%% file: article.tex
\documentclass[12pt]{amsart}
 
 \pdfoutput=1
 
 \usepackage{ifpdf}
\ifpdf
\usepackage[pdftex]{graphicx}
\else
\usepackage[dvips]{graphicx}
\fi
\usepackage{tikz}
 	 \usetikzlibrary{arrows,backgrounds,snakes} 
\usepackage[all]{xy}

\usepackage[pdftex,plainpages=false,hypertexnames=false,pdfpagelabels]{hyperref}
\newcommand{\arxiv}[1]{\href{http://arxiv.org/abs/#1}{\tt arXiv:\nolinkurl{#1}}}
\newcommand{\doi}[1]{\href{http://dx.doi.org/#1}{{\tt DOI:#1}}}

\newcommand{\googlebooks}[1]{(preview at \href{http://books.google.com/books?id=#1}{google books})}
 
 \theoremstyle{plain}
 \newtheorem{Theorem}{Theorem}[section]
 
 \newtheorem{Lemma}[Theorem]{Lemma}
 \newtheorem{Fact}[Theorem]{Fact}
 \newtheorem{DefnThm}[Theorem]{Definition/Theorem}
 \newtheorem{Corollary}[Theorem]{Corollary}

\theoremstyle{definition}
\newtheorem*{Defn}{Definition}
 \newtheorem*{Example*}{Example}
   \newtheorem{Example}{Example}
 \newtheorem*{Notation}{Notation}
 \newtheorem*{Remark}{Remark}
  
  \newcommand{\notetoself}[1]{}
 
 \newcommand{\A}{\mathfrak{A}}
 \newcommand{\C}{\mathbb{C}}
 \newcommand{\id}{\boldsymbol{1}}
 \newcommand{\Hom}[1]{\operatorname{Hom}(#1)}
  \newcommand{\f}[1]{f^{(#1)}}
\newcommand{\gen}[1]{\left\langle #1 \right\rangle}
 \newcommand{\ip}[1]{\left\langle #1 \right\rangle}
 \newcommand{\critical}[2]{\text{critical}_{#1}(#2)}
  \newcommand{\norm}[1]{\left\vert #1 \right\vert}
 \newcommand{\Norm}[1]{\left\Vert #1 \right\Vert}
 
 \newcommand{\proj}[2]{\text{Proj}_{#1}(#2)}
 \newcommand{\PP}{\mathcal{P}}
 \newcommand{\T}{\mathfrak{T}}
 \newcommand{\tr}[1]{\operatorname{tr}(#1)}
 \newcommand{\spn}[1]{\operatorname{span}(#1)}
 \newcommand{\lp}[1]{\mathbf{#1}}
 \newcommand{\pth}[1]{\mathbf{#1}}

 \newcommand{\join}[3]{ \operatorname{Join}_{#1} (#2,#3)}
  \newcommand{\joint}[3]{ \operatorname{Join}'_{#1} (#2,#3)}
\newcommand{\coeff}[2]{\operatorname{Coeff}_{#1}(#2)}

\begin{document}
 \input{TikzStyles}
 
 \title[The Haagerup planar algebra]{A planar algebra construction of the Haagerup subfactor}
 
\author{Emily~Peters}
\email{eep@math.berkeley.edu} 
\address{   Department of Mathematics,   University of California, Berkeley, 94720}

\date{\today}
 
 \begin{abstract} 
Most known examples of subfactors occur in families, coming from algebraic objects such as groups, quantum groups and rational conformal field theories.  The Haagerup subfactor is the smallest index finite-depth subfactor which does not occur in one of these families.  In this paper we construct the planar algebra associated to the Haagerup subfactor, which provides a new proof of the existence of the Haagerup subfactor.  Our technique is to find the Haagerup planar algebra as a singly generated subfactor planar algebra, contained inside of a graph planar algebra.  
 \end{abstract}
 
 \maketitle
 
 \tableofcontents
 
 \section{Introduction}\label{introduction}
 \input{text/introduction}

\section{Background on planar algebras}\label{background}
 \input{text/background}

	 \subsection{Subfactor planar algebras}
 	 \input{text/subfpa}

	\subsection{The planar algebra of a bipartite graph}\label{pabg}
	\input{text/pabg}

	\subsection{Annular Temperley-Lieb modules}\label{ATLmodules}
	\input{text/ATLmodules}

\section{A (potential) generator $T$ for the Haagerup planar algebra}\label{IDgen}
\input{text/IDgen}

\section{Quadratic relations on $T$}\label{Quadratic}
\input{text/Quadratic}
 
	\subsection{Relations among $4$-boxes}
	\input{text/level4}

	\subsection{Interlude:  annular Temperley-Lieb dual bases}\label{dualbasis}
	\input{text/DualBasis}

	\subsection{Relations among $5$- and $6$-boxes}\label{levels56}
	\input{text/levels56}
 
 \section{$T$ generates the Haagerup planar algebra}\label{PATisHaagerup}
\input{text/PATisHaagerup}

\section{The Haagerup planar algebra, and graph planar algebras of other bipartite graphs}\label{othergraphs}
 \input{text/OtherGraphs}

\appendix

 \section{Calculating traces}\label{traces}
 \input{text/Traces}

 \bibliographystyle{plain} 
\bibliography{../../Bibliography/bibliography}

\end{document}

%% file: TikzStyles.tex
\tikzstyle{shaded}=[fill=red!10!blue!20!gray!30!white]
\tikzstyle{shaded line}=[double=red!10!blue!20!gray!30!white, double distance=1.5mm, draw=black]
\tikzstyle{unshaded}=[fill=white]
\tikzstyle{unshaded line}=[double=white, double distance=1.5mm, draw=black]
\tikzstyle{Tbox}=[circle, draw, thick, fill=white, opaque,]
\tikzstyle{empty box}=[circle, draw, thick, fill=white, opaque, inner sep=2mm]
\tikzstyle{background rectangle}= [fill=red!10!blue!20!gray!40!white,rounded corners=2mm] 
\tikzstyle{on}=[very thick, red!50!blue!50!black]
\tikzstyle{off}=[gray]

\tikzstyle{traces}=[scale=.2, inner sep=1mm]
\tikzstyle{quadratic}=[scale=.4, inner sep=1mm, baseline]
\tikzstyle{annular}=[scale=.7, inner sep=1mm, baseline]
\tikzstyle{make triple edge size}= [scale=.4, inner sep=1mm,baseline] 
\tikzstyle{icosahedron network}=[scale=.3, inner sep=1mm, baseline]
\tikzstyle{ATLsix}=[scale=.25, baseline]
\tikzstyle{TL12}=[scale=.15,baseline]
\tikzstyle{PAdefn}=[scale=.7,baseline]
\tikzstyle{TLEG}=[scale=.5,baseline]

%% file: text/introduction.tex
The study of subfactors $N \subset M$ of von Neumann algebras was initiated by Vaughan Jones in \cite{MR696688}, when he defined the index $[M:N]$ of a subfactor and proved it must lie in the set 
$$\{4 \cos^2(\frac{\pi}{n}) | n \geq 3\} \cup [4,\infty],$$
and that all numbers in this set can be realized as the index of a subfactor.  He later defined a finer invariant of subfactors, the principal graph and dual principal graph, which encapsulate the tensor category structure associated to bimodules generated by $M$ as an $N$,$M$-bimodule.  Two natural questions to ask are:  is the principal graph a complete invariant?  Which graphs can be principal graphs?  

Early results in this area are that subfactors with index less than four are classified by their principal graphs, which are exactly the Dynkin diagrams $A_n$, $D_{2n}$, $E_6$ and $E_8$ (see \cite{MR996454} for the broad picture and  \cite{MR1145672, MR1308617, MR1193933, MR1313457} for more details).  Popa \cite{MR1278111} has shown that the finite graphs which occur as principal graphs of index 4 subfactors are exactly the extended Dynkin diagrams, but that the principal graph is not a complete invariant in this case.  However, he goes on to show that amenable subfactors (of which finite-index, finite-depth subfactors are a subset) do have a complete invariant, called the paragroup, which is a structure that includes the principal and dual principal graphs and data on how they fit together.

For subfactors with index more than $4$, the situation is more complicated.  Any index more that $4$ can be realized by a `trivial' subfactor (a subfactor with principal graph $A_\infty$).  But one can still ask which indices more than $4$ can be realized for finite-depth subfactors, i.e. those with finite principal graphs.  The 
first examples started from groups \cite{MR696688} or quantum groups \cite{MR936086}.  The smallest index more than $4$ attained by such examples is $5$.   Taking a combinatorial approach, Goodman, de la Harpe and Jones in \cite{MR999799} constructed a subfactor with finite depth and index $3+\sqrt{3} \approx 4.732$.  
 In \cite{MR1317352}, Haagerup asked what the smallest (finite-depth, irreducible, hyperfinite) subfactor with index more than 4 was.    He showed  such a subfactor cannot have index in the range $(4,\frac{5+\sqrt{13}}{2}\approx 4.303)$, and that if a subfactor of index $\frac{5+\sqrt{13}}{2}$ exists, its principal/dual principal graph pair is the following:

 \input{Diagrams/HandDualH}

In \cite{MR1686551}, Asaeda and Haagerup constructed a hyperfinite subfactor having $(H,H')$ as its principal, dual principal graph pair, and proved its uniqueness.  They first deduced the bimodule structure a Haagerup subfactor must have from the symmetry of $H$, then constructed (using string, or path, algebras) bimodules satisfying these requirements.  From such bimodules one easily produces a subfactor.  Haagerup and Asaeda also used this method to construct and prove uniqueness for a subfactor with index $\frac{5+\sqrt{17}}{2} \approx 4.562$.  These are the smallest ``exotic'' subfactors, that is, subfactors not coming from algebraic structures such as groups, quantum groups, or rational conformal field theories.  

A second construction of the Haagerup subfactor (as a pair of hyperfinite III$_1$ factors) is due to Izumi in \cite{MR1832764}.  He first constructs sectors which form the even part of the desired tensor category (sectors replace bimodules in the III$_1$ setting), then finds a Q-system in order to generate the odd part of the tensor category.
 Izumi's method generalizes to give an infinite family of subfactors, having spoke-shaped principal graphs:  these have $2n+1$ legs with three edges each, which all meet at a single central vertex.  After the Haagerup subfactor, the next smallest of these has index $\frac{7+\sqrt{29}}{2} \approx 6.193$.

A parallel story in the theory of subfactors is the introduction of planar algebras.  
Jones created the planar algebra formalism in \cite{math.QA/9909027} and proved that the tower of relative commutants (also known as the standard invariant, and equivalent to the paragroup and $\lambda$-lattice) of a subfactor is a planar algebra.  Popa has a converse result in \cite{MR1334479};  when translated into planar algebra language, it says that given a planar algebra satisfying certain properties, it is possible to construct a subfactor with that planar algebra as its tower of relative commutants.  More recently, Jones, Shlyakhtenko, Guillonet and Walker \cite{0712.2904, 0807.4146} have provided planar-algebraic proofs of this fact.
These are important results for this paper.  They show that if we can construct a planar algebra with the Haagerup principal graph (from here on, referred to as a {\em Haagerup planar algebra}, this gives the Haagerup subfactor.

In this paper we give an independent proof of the existence of the Haagerup subfactor, by finding a Haagerup planar algebra inside the graph planar algebra of $H$.    
A very fortunate result for us is that if a subfactor planar algebra with principal graph $G$ exists, it is a sub-planar-algebra of the graph planar algebra of $G$ (this result is well-known to experts but is not yet proven in the literature).  
The `graph planar algebra' is the path algebra of a bipartite graph, with a natural planar algebra structure.  Although this planar algebra is in some sense ``too big'' to be a subfactor planar algebra,  it shares many nice properties with subfactor planar algebras and is a convenient place to do calculations.

Working in the graph planar algebra of $H$, we will find a single element which will generate the Haagerup planar algebra.   The main results of this paper are summarized in Theorem \ref{MainThm}, which gives a presentation of the Haagerup planar algebra.  To prove that the planar algebra given by this presentation is the Haagerup planar algebra, we show a small number of relations, then prove that these are sufficient to evaluate any closed tangle.   Thus we have a subfactor planar algebra.  The fact that the planar algebra has principal graph $H$ follows from  Haagerup's result on possible principal graphs with norm $\sqrt{\frac{5+\sqrt{13}}{2}}$.
 
To find this generator inside the graph planar algebra of $H$, we rely on results from \cite{MR1929335, quadratic}.  In these papers Jones considers how to best present a planar algebra.
He classifies annular Temperley-Lieb modules (reproducing the results of \cite{MR1659204}), which provide a convenient way to decompose planar algebras, and gives many relations a planar algebra generator must satisfy.    Examples of subfactor planar algebras for which presentations are known include $A_n$, $D_{2n}$, $E_6$ and $E_8$ (a general form for these presentations is given in \cite{MR1929335}; more combinatorial descriptions of the $D_{2n}$, $E_6$ and $E_8$  planar algebras are found in \cite{d2n}, \cite{ADE}), and Bisch-Haagerup subfactors (\cite{MR1386923}, \cite{0807.4134}).  

This paper is organized as follows:
Section \ref{background} reviews planar algebras, the graph planar algebra of a bipartite graph, and results about annular Temperley-Lieb modules.  In Section \ref{IDgen}, we use these results to identify characteristics of a generator of a Haagerup planar algebra, and then find an element (called $T$) which has these properties in the graph planar algebra of $H$.  In Section \ref{Quadratic}, we make use of annular Temperley-Lieb modules to calculate more relations satisfied by $T$; and in Section \ref{PATisHaagerup}, we show that the planar algebra presented by $T$ and these relations is a subfactor planar algebra with principal graph $H$.  In Section \ref{othergraphs}, we address the question of whether the Haagerup planar algebra appears in graph planar algebras of other bipartite graphs.

The author is grateful to many people for their support for this project.  In particular, I would like to thank Vaughan Jones for suggesting and discussing this problem, and also to thank Stephen Bigelow, Richard Burstein, Scott Morrison and Noah Snyder for interesting discussions, and Jana Comstock for comments on early drafts.  The author was supported in part by NSF Grant DMS0401734 and a fellowship from Soroptimist International.  Part of this work was done while visiting the University of Melbourne.

%% file: text/background.tex
A planar algebra is a family of vector spaces, with a circuit-like structure; elements in the vector spaces can be tied together in various ways to create different outputs.  More precisely, this family of vector spaces has an action by the shaded planar operad (which is a `colored' operad, colored by pairs $\{ n, \pm \}_{n=0,1,2,\ldots}$; see \cite{MR1436912}).

\begin{Defn}
Elements of the shaded planar operad are {\em shaded planar tangles}, which consist of 
\begin{itemize} 
\item
one outer disk $D_o$
\item
$k$ disjoint inner disks $D_i$ 
\item 
A set $S$ of non-intersecting strings between the disks, such that $| S \cap D_o |, |S \cap D_i | \in 2 \mathbb{Z}$
\item
A checkerboard shading on the regions (i.e., connected components of $(D_o \setminus  \bigcup_{i} D_i ) \setminus S$).
\item
A star in a region near the boundary of each disk $D_o, D_i$
\end{itemize}
See Figure \ref{tangleEG} for an example.  We consider two shaded planar tangles equal if they are isotopic.  

\begin{figure}[!ht]
\input{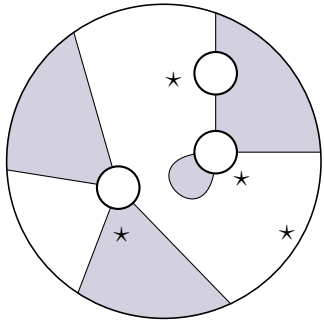}
\caption{A shaded planar tangle}\label{tangleEG}
\end{figure}

By the {\em type} of a disk in a tangle, we mean the pair $(k,\pm)$ determined by half the number of strings on the boundary of the disk ($k$) and whether the star is in an unshaded ($+$) or shaded ($-$) region.  For example, the inner disks in the above figure are of types $(1,+)$, $(2,+)$ and $(2,-)$. 

The operadic structure is given by composition.  Let $A$, $B$ be shaded planar tangles; $A \circ_i B$ exists if the $i^{th}$ inner disk of $A$ has the type as the outer disk of $B$.  In this case, $A \circ_i B$ is the shaded planar tangle built by inserting $B$ in the $i^{th}$ inner disk of $A$ (with the starred regions matching up) and connecting the strings.  For example, 
\input{Diagrams/CompositionEG.tex}
\end{Defn}

\begin{Defn}
  A {\em planar algebra} is a collection of vector spaces $\{ V_{i,\pm} \} _{i=0,1,2,\ldots}$ which is acted on by the shaded planar operad; that is, shaded planar tangles act on tensor products of the $V_{i,\pm}$, in a way that is compatible with composition of tangles.  A shaded planar tangle whose $i^{th}$ inner disk has type $(k_i,\pm_i)$ and outer disk has type $(k_o, \pm_o)$ gives a map $\bigotimes_{i} V_{k_i,\pm_i} \rightarrow V_{k_o, \pm_o}$.  
For example, the first tangle above represents a map $V_{1,+} \otimes V_{2,+} \otimes V_{2,-} \rightarrow V_{3,+}$.
  
Compositional compatibility means composition of tangles and composition of multilinear maps should produce the same result, e.g.  $$(A \circ_1 B) (v_1, \ldots, v_i, w_2, \ldots, w_j) = A( B(v_1, \ldots v_i), w_2, \ldots, w_j).$$  
  
  \end{Defn}
  
  \begin{Defn}
Elements of $V_n$ are called {\em $n$-boxes} and we will sometimes refer to $V_{n}$ as the {\em $n$-box space} of $V$, or as the $n^{th}$ {\em level} of $V$.
\end{Defn}

\begin{Remark}
For $n \geq 1$, $V_{n,+} \simeq V_{n,-}$ (as vector spaces only, not in any algebraic sense), via the shading-changing rotation map
$$\rho^{1/2} = \input{Diagrams/rhohalf}.$$
The isomorphism between $V_{n,+}$ and $V_{n,-}$ is the reason that previous definitions of planar algebras required stars to be in unshaded regions, and only had spaces $V_+$, $V_-$ and $V_i$, $i=1,2,3, \ldots$.  We will drop $+$ and $-$ signs from subscripts when it is either clear, or unimportant, whether we are working in $V_{n,+}$ or $V_{n,-}$.
\end{Remark}
  
\begin{Defn}[The Temperley-Lieb algebra]
The first example of a planar algebra is the Temperley-Lieb algebra with parameter $\delta$ (this algebra was introduced in \cite{MR0498284} and formulated diagrammatically by Kauffman in \cite{MR899057}).  The vector spaces $TL_{i,\pm}$ have a basis (called $B(TL_{i,\pm})$) consisting of non-crossing pairings on $2i$ numbers; these can be drawn as planar tangles with no input disks, $2i$ points on the output disk, and no closed circles (all strings have endpoints on boundary disks).   The number of such pictures is the $i^{th}$ Catalan number $\frac{1}{i+1}\binom{2i}{i}$.

\begin{Example*}
$TL_{3,+}=\{
\input{Diagrams/TLThree/basis}
\}$
\end{Example*}

The action of shaded planar tangles on this basis is straightforward: put the pictures inside the tangle, smooth all strings and throw out closed circles by multiplying the picture by $\delta$ (i.e., if $\tau=\tau' \sqcup \tikz \draw[thick] (0,0) circle (1mm);$ or $\tau=\tau' \sqcup \tikz \filldraw[thick, shaded] (0,0) circle (1mm);$, then $\tau = \delta \tau'$.)  For example,
\input{Diagrams/TLEG}
\end{Defn}

Some shaded planar tangles have particular interpretations in planar algebras.  In the rest of this paper, we will make frequent use of the following.

\begin{Defn} The tangles
$$  	\input{Diagrams/multiplication} \; , \;
 	\input{Diagrams/trace} \; , \text{and} \;
 	\input{Diagrams/inclusion} 
$$
are multiplication, trace and inclusion for $V_{i,+}$; the reverse-shaded tangles are multiplication, trace and inclusion for $V_{i,-}$.
\end{Defn}

\begin{itemize}
	\item
	The multiplication tangle 
 	gives an associative map $$m: V_{k,\pm} \otimes V_{k,\pm} \rightarrow V_{k,\pm}.$$  As usual, we will frequently denote multiplication by putting two elements next to each other, maybe with a dot between them.
 	\item 
	The trace tangle
 	gives a cyclically commutative map $$\operatorname{tr}: V_{k,\pm} \rightarrow V_{0,\pm}.$$  
 	\item 
	The inclusion tangle gives a map $$\iota: V_{k,\pm}  \rightarrow V_{k+1,\pm}$$ which is compatible with multiplication and trace.  
 \end{itemize}
  
\begin{Example*} The identity for multiplication on Temperley-Lieb is $n$ vertical strands, i.e.
$$\id = \begin{tikzpicture}[scale=.6, baseline]	
	\draw (140:2cm) -- (-140:2cm) arc (-140:-110:2cm) -- (110:2cm) arc (110:140:2cm);
	\draw (40:2cm)--(-40:2cm) arc (-40:-70:2cm)--(70:2cm) arc (70:40:2cm);
	
	\node at (0,0) {$\cdots$};
	\node at (180:2cm) [left] {$\star$};
	
	\draw[thick] (0,0) circle (2cm);
\end{tikzpicture}$$
Temperley-Lieb is multiplicatively generated by the {\em Jones projections} 
$$e_i = \frac{1}{\delta} \cdot 
\begin{tikzpicture}[scale=.6, baseline]
	\draw[thick] (0,0) circle (2cm);
	
	\draw (145:2cm) -- (-145:2cm) arc (-145:-110:2cm) -- (110:2cm) arc (110:145:2cm);
	\draw (35:2cm)--(-35:2cm) arc (-35:-70:2cm)--(70:2cm) arc (70:35:2cm);
	
	\draw (80:2cm) .. controls ++(-90:1cm) and ++(-90:1cm) .. (100:2cm);
	\draw (-80:2cm) .. controls ++(90:1cm) and ++(90:1cm) .. (-100:2cm);

	\node at (-1.1,0) {\tiny{$\cdots$}};
	\node at (1.2,0) {\tiny{$\cdots$}};	
	\node at (180:2cm) [left] {$\star$};
\end{tikzpicture}$$
\end{Example*}

%% file: Diagrams/TLThree/basis.tex
\begin{tikzpicture}[TLEG]
	\filldraw[shaded]  (30:1cm) arc (30:90:1cm) arc (30:-150:5mm) arc (150:210:1cm) -- cycle; 
	\filldraw[shaded]  (0,-1) arc (-90:-30:1cm) arc (30:210:5mm);
	\draw[thick] (0,0) circle (1cm);
	\node at (120:1.3cm) {$\star$};
\end{tikzpicture},
\begin{tikzpicture}[TLEG]
	\filldraw[shaded]  (0,1) arc (90:30:1cm) arc (90:270:5mm) arc (-30:-90:1cm) -- cycle; 
	\filldraw[shaded]  (-1,0) arc (180:210:1cm) arc (-90:90:5mm) arc (150:180:1cm);
	\draw[thick] (0,0) circle (1cm);
	\node at (120:1.3cm) {$\star$};
\end{tikzpicture},
\begin{tikzpicture}[TLEG, rotate=180]
	\filldraw[shaded]  (150:1cm) arc (150:90:1cm) arc (-210:-30:5mm) arc (30:-30:1cm)--cycle;
	\filldraw[shaded] (-90:1cm) arc (-90:-150:1cm) arc (150:-30:5mm);
	\draw[thick] (0,0) circle (1cm);
	\node at (-60:1.3cm) {$\star$};
\end{tikzpicture},
\begin{tikzpicture}[TLEG]
	\filldraw[shaded]  (0,-1) arc (-90:-30:1cm) arc (30:210:5mm);
	\filldraw[shaded]  (-1,0) arc (180:210:1cm) arc (-90:90:5mm) arc (150:180:1cm);
	\filldraw[shaded] (90:1cm) arc (90:30:1cm) arc (-30:-210:5mm);
	\draw[thick] (0,0) circle (1cm);
	\node at (120:1.3cm) {$\star$};
\end{tikzpicture},
\begin{tikzpicture}[TLEG]
	\filldraw[shaded]  (90:1cm) arc (30:-150:5mm) arc (150:210:1cm) arc (150:-30:5mm) arc (-90:-30:1cm) arc (-90:-270:5mm) arc (30:90:1cm);
	\draw[thick] (0,0) circle (1cm);
	\node at (120:1.3cm) {$\star$};
\end{tikzpicture}

%% file: text/subfpa.tex
Subfactor planar algebras are planar algebras with additional structure.    They are called `subfactor' because the planar algebra of a subfactor always has these properties.  Furthermore, a planar algebra with these properties is always the standard invariant of some subfactor.  Further details on this connection are found in chapter 4 of \cite{math.QA/9909027}.

The properties that define subfactor planar algebras make them easier to work with than general planar algebras.  In particular, the requirements that each space be finite dimensional and have an inner product means the tools of linear algebra are available to us when we work with subfactor planar algebras.

\begin{Defn}
A {\em subfactor} planar algebra is a planar algebra (over $\mathbb{C}$) which has
\begin{enumerate}
\item\label{spa:inv} Involution:
 a $*$ on each $V_{i,\pm}$ which is compatible with reflection of tangles (so that $\tau^* (v_1^*,v_2^*,\ldots v_n^*)=\tau (v_1,v_2,\ldots v_n)^*$), 
 \item\label{spa:dim} Dimension restrictions:
  $\dim (V_{0,+})=\dim(V_{0,-})=1$ and $\dim{V_{k,\pm}}<\infty$ for all $k$,
   \item\label{spa:sph} Sphericality:
The left trace $\operatorname{tr}_l:V_{1,\pm} \rightarrow V_{0,\mp}\simeq \mathbb{C}$ and the right trace $\operatorname{tr}_r:V_{1,\pm} \rightarrow V_{0,\pm }\simeq \mathbb{C}$ are equal (equivalently, the action of planar tangles is invariant under spherical isotopy),
 \item\label{spa:posdef}  Inner Product:
The bilinear form on $V_{n,\pm}$ given by  $\left< a,b \right>:=\tr{b^* a}$ is positive definite.
\end{enumerate}
\end{Defn}

Note that the one-dimensionality of $V_{0,+}$ and $V_{0,-}$ implies that closed circles must count for a constant, and property (\ref{spa:sph}) implies that shaded and unshaded closed circles count for the same value.

\begin{Notation}
In a subfactor planar algebra, $\delta$ will denote the value of closed circles.  Sometimes we change variables so that 
$$ \delta = [2]_q= (q+q^{-1})$$
in order to use {\em quantum numbers:}
$$[n]_q :=\frac{q^n-q^{-n}}{q-q^{-1}} = q^{n-1} + q^{n-3} + \cdots + q^{-n+1} + q^{-n+3}.$$
We will often write $[n]$ instead of $[n]_q$ in situations where the value of $q$ is known.
\end{Notation}

\begin{Example*}
The Temperley-Lieb planar algebra always meets conditions \ref{spa:inv}, \ref{spa:dim} and \ref{spa:sph}:  The involution is defined by reflection;  $TL_{0,+}=\C \{ \tikz \draw[thick] (0,0) circle (1mm); \}$  and $TL_{0,-}=\C \{ \tikz \draw[thick, shaded] (0,0) circle (1mm); \}$ are both $1$-dimensional; and shaded and unshaded circles both count for $\delta$, hence the planar algebra is spherical.

If $\delta \geq 2$, the bilinear form defined by the trace tangle is positive definite, and so Temperley-Lieb is a subfactor planar algebra.  If $\delta=2 \cos{\frac{\pi}{n}}$ for some $n \geq 3$, the bilinear form is positive semidefinite, and we can form a subfactor planar algebra by quotienting Temperley-Lieb by all $x \in TL$ such that $\tr{x^* x}=0$.
\end{Example*}

In fact, the subfactor Temperley-Lieb planar algebra is an initial object in the category of subfactor planar algebras. 

\begin{Fact}\label{TLinjects}
Let $V$ be a subfactor planar algebra with parameter $\delta$.  Then the map
$$TL(\delta) \hookrightarrow \Hom{\mathbb{C},V} \cong V $$
given by interpreting a Temperley-Lieb diagram as a shaded planar tangle with no inputs 
\begin{itemize}
\item is injective if $\delta \geq 2$;
\item has kernel $\operatorname{Rad}(\left<, \right>)$ if $\delta < 2$.
\end{itemize}
\end{Fact}

It is the algebraic structure of subfactor planar algebras that make them so nice to work with.  This algebraic structure is summarized by its principal graph.  In particular, this graph encodes the decomposition into irreducibles of the inclusion of a minimal idempotent in the next level.  Let's make these notions more precise:

\begin{Defn}
\begin{itemize}
	\item	An {\em idempotent} is an element $p$ of some $V_{k,\pm}$ such that $p=p^*=p p$.  We often draw idempotents in rectangles instead of disks, with an implicit star on the left side.
	\item An idempotent $p \in V_{k,\pm}$ is {\em minimal} if $p \cdot V_{k,\pm} \cdot p$ is one-dimensional.  
	\item Two idempotents $p \in V_{i,\pm}$ and $q \in V_{j,\pm}$ are {\em isomorphic} if there is an element $f$ in $V_{i \rightarrow j,\pm}$ (this is the same space as $V_{(i+j)/2,\pm}$, but drawn with $i$ strings going up and $j$ strings going down) such that $f^*f=p$ and $f f^*=q$:
	$$\input{Diagrams/DefnIsom}$$
\end{itemize}
\end{Defn}

\begin{Example*}[The Jones-Wenzl idempotents  of Temperley-Lieb] A {\em Jones-Wenzl idempotent}, first defined in \cite{MR873400} and denoted $f^{(n)}$, is the unique idempotent in $TL_{n,\pm}$  which is orthogonal to all $TL_{n,\pm}$  basis elements except the identity.  In symbols this says $\f{n}\neq 0$, $\f{n}\f{n}=\f{n}$ and $e_i \f{n}=\f{n} e_i=0$ for all $i \in \{ 1, \ldots , n-1\}$ (since $e_i$ multiplicatively generate $TL$).  These have $\tr{f^{(n)}}=[n+1]$, are minimal, and satisfy {\em Wenzl's relation}:
$$\input{Diagrams/WenzlRelation}$$
\end{Example*}

\begin{Defn}
The {\em principal graph} consists of vertices and edges.
The vertices of the principal graph are the isomorphism classes of idempotents in $V_{k,+}$ for any $k$.  
There are $n$ edges between $p$ and $q$ if $\iota(p) =
\begin{tikzpicture}[baseline]
	\foreach \x in {-3,-1,1,3,5} \draw (\x mm,-4mm)--(\x mm, 4mm);
	\node[rectangle,draw,fill=white] {$\; p \;$};
\end{tikzpicture}
$ contains $n$ copies of $q$, meaning that when $\iota(p)$ is decomposed into minimal idempotents, $n$ of these are isomorphic to $q$.  The dual principal graph is constructed in the same way, for the idempotents in $V_{k,-}$.
\end{Defn}

\begin{Example*}[The principal graph of Temperley-Lieb is $A_n$ or $A_\infty$]
The vertices of the principal graph of Temperley-Lieb are the Jones-Wenzl idempotents.
Wenzl's relation says that if $[n+1]\neq 0$, the projection $\iota (f^{(n)})$ decomposes into minimal projections isomorphic to $f^{(n+1)}$ and $f^{(n-1)}$.  So there is an edge in the principal graph between $f^{(k)}$ and $f^{(k+1)}$ as long as $[k+1] \neq 0$.  

If $\delta>2$, we never have $[k] = 0$, so the principal graph is the Dynkin diagram $A_{\infty}$.
 If $\delta = 2 \cos{\frac{\pi}{n}}$ then $[n]=0$, so by the positive definiteness of the inner product $f^{(n-1)}=0$.  Thus the principal graph has $n-1$ vertices and is the Dynkin diagram $A_{n-1}$.
 \end{Example*}

One more thing worth mentioning about the Jones-Wenzl idempotent is that a reasonably explicit recursion is given, in \cite{morrison}, which calculates the coefficients of $TL$ pictures in the Jones-Wenzl idempotents.  We will make extensive use of these in Section \ref{dualbasis}.

%% file: text/pabg.tex
The planar algebra of a bipartite graph $G$, defined in \cite{MR1865703}, is another example of a planar algebra.  It shares with Temperley-Lieb the property that closed circles count for a constant --- in this case, 
$\delta= \left\Vert G \right\Vert $ (where $\Norm{G}$ is the operator norm of the adjacency matrix of $G$).  It also has an involution, is spherically invariant and has a positive definite inner product.  Usually, however, it is not a subfactor planar algebra, because for most graphs, the zero-box spaces are too big.

\begin{Notation}
All graphs in this paper are simply laced.  Therefore paths or loops can, and will, be entirely described by the vertices they pass through. 
When we concatenate two paths written this way, we have to drop a vertex:
$$a b \ldots c d \sqcup d e \ldots f g = a b \ldots c d e \ldots  f g.$$
When a loop is described in this notation, we might forget to notate the last vertex.

Paths and loops are written in boldface, and indices indicate the position of a vertex in a path or loop: $\pth{p} = p_1 p_2 \ldots p_k$.
\end{Notation}

\begin{Defn}
$PABG(G)_{0,+}$ has the even vertices $U_+$ of $G$ as a basis; $PABG(G)_{0,-}$ has the odd vertices $U_-$ as a basis.  The space $PABG(G)_{i,+}$ has a basis consisting of all based loops of length $2i$ on the graph $G$, based at any even vertex; similarly the space $PABG(G)_{i,-}$ has a basis consisting of all based loops of length $2i$ on the graph $G$, based at any odd vertex.
\end{Defn}

To make this a planar algebra, we need to specify how a shaded planar tangle acts on basis elements.  The idea is to sum over all states of the tangle (a state is an assignment of vertices and edges to regions and strings) which are compatible with the input.  

\begin{Defn} Suppose $\tau$ is a tangle with $k$ input disks.  A {\em state} on $\tau$ is an assignment of even vertices to unshaded regions, odd vertices to shaded regions, and edges to strings in such a way that if an edge is assigned to a string, its endpoints are assigned to the two regions which touch that string.   Given a state $\sigma$ on $\tau$, we can read clockwise around any disk and get a loop on $G$ (the starred region tells us where to base the loop).  Let $\partial_i (\sigma)$ be the loop read from the $i^{\text{th}}$ inner disk, and let $\partial_o (\sigma)$ be the loop read from the outer disk.
\end{Defn}

\begin{Defn}
Define the action of $\tau$ on loops $\lp{p_i}$ by
$$\tau(\lp{p_1},\ldots,\lp{p_k}):= \sum_{\text{states } \sigma \text{ s.t. } \partial_i(\sigma)=\lp{p_i}} c(\tau, \sigma) \cdot \partial_o (\sigma) ,$$
where $c(\tau, \sigma)$ is a number defined below.  Extend this action from the basis of loops to all of $PABG(G)_{i,\pm}$ by linearity.
\end{Defn}

Why do we need a correction factor $c$ in the above equation?
In order for the action of tangles to be isotopy invariant and compatible with composition, we need $c$ to be isotopy invariant and multiplicative (so if $\sigma=\sigma_1 \circ \sigma_2$ is a state on the tangle $\tau=\tau_1 \circ \tau_2$, $c(\tau, \sigma)=c(\tau_1, \sigma_1)c(\tau_2, \sigma_2)$).  Unfortunately, the obvious choice of $c \equiv 1$ will not work because  $c$ has a third role:  it should cause closed circles to count for a constant $\delta$.  

\begin{Defn}
Let $\Lambda$ be the adjacency matrix of $G$.   Its  {\em Perron-Frobenius eigenvalue} is the maximal modulus eigenvalue (which is necessarily real); call this $\delta$.   Let $\lambda$ be the eigenvector with eigenvalue $\delta$, called the {\em Perron-Frobenius eigenvector}.  We normalize $\lambda$ to have norm $\sqrt{2}$, and let $\lambda(a)$ be the entry of $\lambda$ corresponding to vertex $a$.  
\end{Defn}

The normalization above may seem peculiar, but it is chosen because we want $\sum_{a \in U_+} \lambda(a)^2=1$ and $\sum_{a \in U_-} \lambda(a)^2=1$.

\begin{Defn}
To define $c(\sigma, \tau)$, we first put $\tau$ in a `standard form':
isotope $\tau$ so that all strings are smooth, and all its boxes are rectangles, with the starred region on the left, half the strings coming out of the top and the other half coming out the bottom.
Let $E(\tau)$ be the set of maxima and minima on strings of the standard form of $\tau$.  If $t \in E(\tau)$ is a max or min on a string of $\tau$, let 
$\sigma(t_{\text{convex}})$ 
be the vertex assigned by $\sigma$ to the region touching $t$ where the string is convex, and 
${\sigma}(t_{\text{concave}})$
 be the vertex assigned by $\sigma$ to the region touching $t$ where the string is concave.  Then
$$c(\tau, \sigma):=\prod_{t \in E(\tau)} \sqrt{\frac{\lambda({\sigma}(t_{\text{convex}}))}{\lambda({\sigma}(t_{\text{concave}}))}}$$
\end{Defn}

Though hard to state in the abstract, this definition isn't too hard to work with:

\begin{Example}\label{PABGcalc}
Take our bipartite graph to be 
$$H=\input{Diagrams/Hlabelled};$$ 
the Perron-Frobenius data for $H$ is $\delta=\sqrt{\frac{5+\sqrt{13}}{2}}$ and 
\begin{align*}
\lambda(z_i) & = \sqrt{\frac{1}{78} \left(13-3 \sqrt{13}\right)}, & 
\lambda(a_i) & =\sqrt{\frac{1}{78} \left(13-\sqrt{13}\right)}, \\
\lambda(b_i) & =\sqrt{\frac{1}{78} \left(13+3 \sqrt{13}\right)}, & 
 \lambda(c) & =\sqrt{\frac{1}{26} \left(13+\sqrt{13}\right)}. 
\end{align*}
As $H$ is has no multiple edges, we can describe paths on $H$ by referring to which vertices they pass through. 
Here is an example of a planar diagram in standard form acting on a basis element:
\begin{align*}\begin{tikzpicture}[scale=.7,baseline]
	\filldraw[shaded] (0,2)--(0,1) arc (-180:0:5mm) -- (1,2);
	\filldraw[shaded] (0,-1)--(0,-.7) arc (180:0:5mm) -- (1,-1);
	\node at (0,1.2) [rectangle, draw, thick, fill=white, inner sep=2mm] {};
	\draw[thick] (-1,-1) rectangle (2,2);
\end{tikzpicture}( b_0 a_0) & = 
\begin{tikzpicture}[scale=.7,baseline]
	\filldraw[shaded] (0,2)--(0,1) arc (-180:0:5mm) -- (1,2);
	\filldraw[shaded] (0,-1)--(0,-.7) arc (180:0:5mm) -- (1,-1);
	\node at (0,1.2) [rectangle, draw,thick, fill=white, inner sep=2mm] {};
	\draw[thick] (-1,-1) rectangle (2,2);
	\node at (.6,1.2) {{\textcolor{white}{$a_0$}}};
	\node at (-.6,1.2) {{\textcolor{gray}{$b_0$}}};
	\node at (.5,-.6) {{\textcolor{white}{?}}};
\end{tikzpicture}
=
\frac{\lambda(a_0)}{\lambda(b_0)} b_0 a_0 b_0 a_0 + \frac{\sqrt{\lambda(a_0) \lambda(c)}}{\lambda(b_0)} b_0 a_0 c a_0 \\
& =
\sqrt{4-\sqrt{13}} \cdot b_0 a_0 b_0 a_0 + \sqrt{\frac{39(\sqrt{13}-3)}{2} } \cdot b_0 a_0 c a_0 \\
\end{align*}
\end{Example}

Note that the planar algebra of a bipartite graph as defined is not often a subfactor planar algebra; $\dim(PABG(G)_{0,+})$ equals the number of even vertices of $G$, $\dim(PABG(G)_{0,-})$ equals the number of odd vertices of $G$, and these are both $1$ only in rather dull cases.  However, this planar algebra does have an involution $*$ defined on loops by traversing them backwards.  What's more, the planar algebra of a bipartite graph has a genuine trace (that is, a cyclically commutative map $PABG(G)_{i,\pm} \rightarrow \C$):

\begin{Defn} The trace $Z: PABG(G)_{i,\pm} \rightarrow \mathbb{C}$ is defined as the composition 
\begin{equation*}
	\xymatrix{
	PABG(G)_{i,\pm} \ar@{->}[r]^{\operatorname{tr}} & PABG(G)_{0,\pm} \ar@{->}[r]^{Z_0} & \C\qquad\qquad\quad }
\end{equation*}
where $\operatorname{tr}$ is the trace tangle, and $Z_0$ is the linear extension of the map $a \mapsto \lambda(a)^2$.
\end{Defn}

 This trace gives us a positive definite inner product, and by construction (and normalization of $\lambda$) agrees with the trace on $TL(\delta) \subset PABG(G)$.

%% file: text/ATLmodules.tex
We begin by summarizing some of the definitions and theorems of \cite{MR1929335}.

\begin{Defn}
An {\em annular Temperley-Lieb tangle} is any shaded planar tangle having exactly one input disk; for example
$$\input{Diagrams/AnnExample1}\, , \,
\input{Diagrams/AnnExample2}\, , \text{ and }
\input{Diagrams/AnnExample3} \, . $$
We write $ATL_{k,\pm \rightarrow m,\pm}$ for the set of all annular Temperley-Lieb tangles with type $(k,\pm)$ inner disk and type $(m,\pm)$ outer disk.  If $ATL$ is being applied to an element $R$ of type $(k,\pm)$, we may simply write $ATL_{m,\pm}(R)$ for $ATL_{k,\pm \rightarrow m,\pm}(R)$.
\end{Defn}

\begin{Defn} An {\em annular Temperley-Lieb module} is a family of vector spaces $V_{i,\pm}$ which has an action by annular Temperley-Lieb.  This action should be compatible with composition.
\end{Defn}

Note that, in particular, all planar algebras are annular Temperley-Lieb modules.  The idea of annular Temperley-Lieb modules is useful because it is weaker than a planar algebra, and lets us break planar algebras down into smaller, easier-to-understand pieces.

\begin{Defn}
An annular Temperley-Lieb module is {\em irreducible} if it has no proper submodules.  This is equivalent to being indecomposable.
\end{Defn}

Irreducible annular Temperley-Lieb modules were studied first in \cite{MR1659204}, and recast in the language of planar algebras in \cite{MR1929335}.  If $\delta>2$, they are classified by an eigenvector and its eigenvalue.

\begin{DefnThm}[\cite{MR1659204}, \cite{MR1929335}]
If $V$ is an irreducible annular Temperley-Lieb (for $\delta>2$) module, it has a {\em low weight space:}  some $k$ such that $V_{i,\pm} = 0$ if $i<k$ and $V_{k,+} \neq 0$ or $V_{k,-} \neq 0$.  
Then $\dim (V_{k,\pm})\leq 1$, and any non-zero element $T \in V_{k,\pm}$ generates all of $V$ as an annular Temperley-Lieb module. In particular, we may choose $T=T^*$.

If $V$ has low weight $(0,+)$,  $T$ is an eigenvector of the double-circle operator
$$\sigma
=\input{Diagrams/sigma},$$
and the eigenvalue $\mu^2$ of $T$ can be any number in $[0,\delta^2]$.  We write $V^{(0,+),\mu}$ for this module.

If $V$ has low weight $(0,-)$,  $T$ is an eigenvector of the double-circle operator, $\sigma$ from above with its shading reversed.
The eigenvalue $\mu^2$ of $T$ can be any number in $[0,\delta^2]$.  We write $V^{(0,-),\mu}$ for this module.

If $V$ has low weight $k>0$, $T$ is an eigenvector of the rotation operator 
$$\rho
=\input{Diagrams/rho},$$ 
and the eigenvalue $\zeta$ of $T$ can be any $k^{th}$ root of unity.  We write $V^{k,\zeta}$ for this module.
\end{DefnThm}

\begin{Remark}
For $V^{k,\zeta}$ the condition that the generator $T$ is a low weight element (i.e., if $m<k$ 
then $ATL_{m,\pm} (T)=0$) is equivalent to the easier-to-check condition that
all of the capping-off operators in $ATL_{k,\pm \rightarrow k-1,\pm}$ give $0$:  if we define
$$\epsilon_1=\input{Diagrams/epsilon1}\, , 
\epsilon_2=\input{Diagrams/epsilon2}\, , \ldots,
\epsilon_{2k}=\input{Diagrams/epsilon2k}$$
then  $\epsilon_i(T)=0$ for all $i$.  
\end{Remark}

 To decompose any annular Temperley-Lieb module (for instance, a planar algebra) into irreducible annular Temperley-Lieb modules, we count dimensions and pay attention to the action of $\rho$, the rotation operator.  

\begin{Theorem}\label{IrredDims} If $k>0$, for $m > k$
$$\dim (V_m^{k,\zeta}) = \binom{2m}{m-k};$$
If $k=0$, the dimensions depend on $\mu$:
\begin{align*}
\dim (V_{0,\pm}^{(0,\pm),0}) = 1, \; &
\dim (V_{0,\mp}^{(0,\pm),0}) = 0, &&
\dim (V_m^{(0,\pm),0}) = \frac{1}{2} \binom{2m}{m}, \\
\dim (V_{0,\pm}^{0,\delta})&=1  && \dim (V_m^{0,\delta})= \frac{1}{m+1} \binom{2m}{m}  , \\
\dim (V_{0,\pm}^{0,\mu})&=1  && \dim (V_m^{0,\mu})  = \binom{2m}{m} &
\text{ if } \mu \in (0,\delta)
\end{align*}
\end{Theorem}

To see how we can make use of this information, we are going to deduce a number of relations on the (so far hypothetical) Haagerup planar algebra.

\begin{Example}[Decomposition of a Haagerup planar algebra into irreducibles]\label{HaagerupDecomp}
Consider a (hypothetical) planar algebra $P$ having principal graph $H$.  It will have $\delta=\sqrt{\frac{5+\sqrt{13}}{2}}$, and from the principal graph we can reconstruct the Bratteli diagram and sequence of dimensions:

\begin{figure}[!ht]
\input{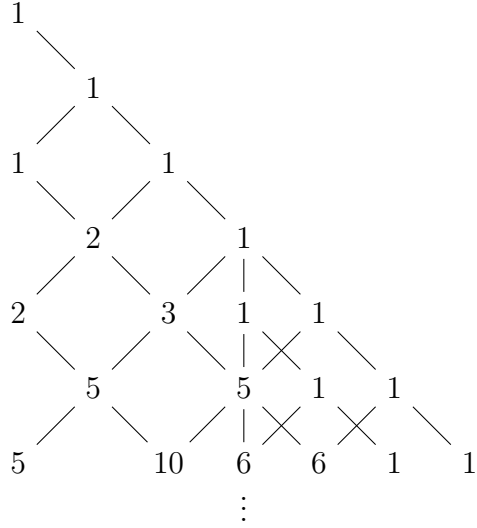}
\caption{The Bratteli diagram of a Haagerup planar algebra.  
The Bratelli diagram is built from the principal graph; 
the number at a vertex of the $k^{th}$ row represent the multiplicity of the corresponding minimal idempotent in $P_k$.  Thus, the sum of the squares of the numbers in the $k^{th}$ row is the dimension of $P_k$.}
\end{figure}

\begin{multline*}
\dim (P_{0,\pm})=1 , \, \dim (P_1)=1 , \, \dim (P_2)= 2, \, \dim (P_3)= 5, \\
 \dim (P_4)= 15, \, \dim (P_5)=52 , \, \dim (P_6)=199 , \ldots
\end{multline*}

From Fact \ref{TLinjects}, $P$ contains a copy of Temperley-Lieb; and Temperley-Lieb is an irreducible annular Temperley-Lieb module: $TL \simeq V^{0,\delta}$.
The sequence of dimensions of $TL$ is 
\begin{multline*}
\dim (TL_{0,\pm})=1 , \, \dim (TL_1)=1 , \, \dim (TL_2)= 2, \, \dim (TL_3)= 5, \\
 \dim (TL_4)= 14, \, \dim (TL_5)=42 , \, \dim (TL_6)=132 , \ldots
\end{multline*}

These sequences differ for the first time at level 4, meaning that $P$ contains no other low weight $0$ modules, and no low weight $1$, $2$ or $3$ modules, but $P$ does contains one low weight 4 module:
$$\dim(V_4^{4,\zeta})=1 , \, \dim(V_5^{4,\zeta})=10 , \, \dim(V_6^{4,\zeta})=66 , \ldots$$
From  \cite[Thm. 6.1.2]{quadratic}, we can deduce that in this case, $\zeta=-1$.

 By counting dimensions we see that $P_5=TL_5 \oplus V_5^{4,-1}$.  Then $\dim(P_6)=199$ and $\dim(TL_6 \oplus V_6^{4,\zeta})=132+66$, so $P$ also contains a lowest weight 6 module:
$$P = (TL\simeq V^{0,\delta}) \oplus V^{4,-1} \oplus V^{6, \zeta} \oplus \cdots$$
\begin{flushright}
\qedsymbol
\end{flushright}
\end{Example}

This decomposition of $P$ up to level $6$ contains information about a generators-and-relations presentation of $P$.  Since $TL$ is generated (as a planar algebra) by the empty set, the first generator of $P$ must have low weight $4$.  Call this generator $T$ and pick it to have sign $+$ (i.e., its star is in an unshaded region).  Now, any diagram with some  $T$s inside it, and only eight or ten strands on the outer boundary, is contained in $ATL_{4}(T) \oplus TL_4$ or $ATL_{5}(T) \oplus TL_5$.  So $P$ must have many ``quadratic'' (and higher-order) relations.  For instance, if $R$ is a diagram with two $T$s connected by four strands (for instance, $T^2$), then
$$R \in ATL_{4}(T) \oplus TL_4$$
or, if $S$ is a diagram with two $T$s connected by three strands,
$$S \in ATL_{5}(T) \oplus TL_5.$$
Since there is a single low weight 6 element, there can be only one new element at level 6.  Call this $N$.  We then expect to see a number of quadratic and higher-order relations saying that any diagram with some $T$s, and twelve strands on the outer boundary, is in $ATL_{6}(T) \oplus N \oplus TL_6$.
In Section \ref{Quadratic}, we will see exactly these kinds of relations appearing.

%% file: Diagrams/sigma.tex
\begin{tikzpicture}[scale=.5,baseline]
	\draw[shaded] (0,0) circle (20mm);
	\draw[unshaded] (0,0) circle (10mm);
	\draw[thick] (0,0) circle (5mm);
	\draw[thick] (0,0) circle (25mm);
\end{tikzpicture}

%% file: text/IDgen.tex
We already know (from Example \ref{HaagerupDecomp}) a fair bit about what a presentation of the Haagerup planar algebra should look like.  The main result of this paper is that the following is a presentation of the Haagerup planar algebra.  

\begin{Theorem}\label{MainThm}
  There is a unique element $T \in PABG(H)_{4,+}$ satisfying the following relations, and the planar algebra generated by $T$ is a subfactor planar algebra with principal graph $H$.
  \begin{enumerate}
  \item $T=T^*$,
  \item $\input{Diagrams/AnnRelns/Cap},$
  \item \hspace{.5in}  $\input{Diagrams/AnnRelns/Rotate}$,
  \item $\input{Diagrams/QuadRelns/4}$ 
  \item \hspace{.5in}  $\input{Diagrams/QuadRelns/5simplified}$
  \item $\input{Diagrams/QuadRelns/6withN}$
  \end{enumerate}
\end{Theorem}

This theorem will be proved over the next three sections.  In this section, we find the elements of $PABG(H)_{4,+}$ satisfying (1), (2) and (3) above.  We then find that relation (4) rules out all of these elements except one, and call the remaining element $T$.
In Section \ref{Quadratic}, we prove relations (4), (5) and (6), and in Section \ref{PATisHaagerup} we will show that $T$ does, in fact, generate a subfactor planar algebra with principal graph $H$.  

To work in the planar algebra of the bipartite graph $H$, we use the notation of Example \ref{PABGcalc} of Section \ref{pabg}.  Let $X^{4,-1}$ be the subspace of $PABG(H)_{4,+}$ consisting of elements satisfying the relations
$$\epsilon_i(T)=0 \text{ for }i=1,\ldots,8, \quad \rho(T)=-T.$$

\begin{Theorem}
$X^{4,-1}$ is four-dimensional.
\end{Theorem}

\begin{proof}
To see this we decompose $PABG(H)$ into irreducible annular Temperley-Lieb modules, using the methods of \cite{MR1929335} as outlined in Section \ref{ATLmodules}.

Let $\Lambda$ be the even-odd adjacency matrix of $H$.  The eigenvalues of $\Lambda \Lambda^t$ and $\Lambda^t \Lambda$  are 
$ \{ \frac{5+\sqrt{13}}{2}, 2, 2 , \frac{5-\sqrt{13}}{2}, 0 , 0 \}$
and
$ \{ \frac{5+\sqrt{13}}{2}, 2, 2 , \frac{5-\sqrt{13}}{2} \}$,
so 
$$PABG(H)_{0,+}  = V^{0,\delta}_+ \oplus 3 V^{0,\mu}_+ \oplus 2 V^{(0,+),0}_+$$
and 
$$PABG(H)_{0,-}  = V^{0,\delta}_- \oplus 3 V^{0,\mu}_- \oplus 2 V^{(0,+),0}_-.$$

We then know that 
$$PABG(H)_{1}  \supset V^{0,\delta}_1 \oplus 3 V^{0,\mu}_1 \oplus 2 V^{(0,+),0}_1$$
 From counting length-2 loops on $H$ based at even vertices, we know that $PABG(H)_1$ is 9-dimensional.  Theorem \ref{IrredDims} tells us that $V^{0,\delta}_1 \oplus 3 V^{0,\mu}_1 \oplus 2 V^{(0,+),0}_1$ is also 9-dimensional, so
$$PABG(H)_{1}  = V^{0,\delta}_1 \oplus 3 V^{0,\mu}_1 \oplus 2 V^{(0,+),0}_1.$$

Next we compute
$$PABG(H)_{2}  \supset V^{0,\delta}_2 \oplus 3 V^{0,\mu}_2 \oplus 2 V^{(0,+),0}_2$$
and counting length-4 loops on $H$ based at even vertices tells us $\dim (PABG(H)_2)=27$, while Theorem \ref{IrredDims} tells us that $\dim( V^{0,\delta}_2 \oplus 3 V^{0,\mu}_2 \oplus 2 V^{(0,+),0}_2 )=26$.  Therefore $PABG(H)$ must contain a copy of $V^{2,\omega}$. 

To decide whether $\omega=1$ or $-1$, we consider the action of $\rho$ on $PABG(H)_2$ and on the irreducible modules.  On the (even-based) loop basis of $PABG(H)_2$, $\rho$ has $9$ fixed points, and $9$ orbits of two elements each.  If $\rho(x)=x$ then $\rho$ has eigenvalue $1$; if $\rho(x)=y$ and $\rho(y)=x$ then $x+y$ and $x-y$ are eigenvectors of $\rho$, with eigenvalues $1$ and $-1$.  So $\rho$ on $PABG(H)_2$ has eigenvalue $1$ with multiplicity $18$, and eigenvalue $-1$ with multiplicity $9$.  

To figure out the action of $\rho$ on 
$V^{0,\delta}_2 \oplus 3 V^{0,\mu}_2 \oplus 2 V^{(0,+),0}_2$, we pick a basis on which $\rho$ acts by permutation.  For $V^{0,\delta}$ this is just Temperley-Lieb pictures with four boundary points.  On $V^{0,\mu}$, the basis is $TL_2$ pictures with the generator $T$ of $V^{0,\mu}$ in one of the four regions.  On $V^{(0,+),0}$, the basis is $TL_2$ pictures with the generator $T$ in one of the two unshaded regions.  On these bases, $\rho$ has $10$ fixed points and $8$ orbits with two elements each, and therefore $\rho$ has eigenvalue $1$ with multiplicity $18$, and eigenvalue $-1$ with multiplicity $8$.  

Therefore $V^{2,-1} \subset PABG(H)$, and 
$$PABG(H)_2  = V^{0,\delta}_2 \oplus 3 V^{0,\mu}_2 \oplus 2 V^{(0,+),0}_2 \oplus V^{2,-1}_2.$$

This implies
$$PABG(H)_3  \supset V^{0,\delta}_3 \oplus 3 V^{0,\mu}_3 \oplus 2 V^{(0,+),0}_3 \oplus V^{2,-1}_3$$
and, by counting the dimensions of the left and right spaces above, we find that $PABG(H)$ contains five modules of the form $V^{3,\zeta}$.

To find the multiplicities of $V^{3,1}$, $V^{3,\omega}$, and $V^{3,\omega^2}$ in $PABG(H)$, we compute the action of $\rho$ on $PABG(H)_3$ and on its known submodules.  By counting fixed points and orbits on the loop basis, we find $\rho$ has eigenvalue $1$ with multiplicity $38$, and eigenvalues $\omega$ and $\omega^2$ with multiplicity $29$.

On $V^{0,\delta}_3 \oplus 3 V^{0,\mu}_3 \oplus 2 V^{(0,+),0}_3 \oplus V^{2,-1}_3$, we pick a basis on which $\rho$ acts by permutation.  The modules $V^{0,?}_3$ have a basis of $TL_3$ pictures, or $TL_3$ pictures with a generator inserted in all, or in all unshaded, regions (for the $\mu=\delta$, $\delta > \mu >0$, and $\mu=0$ cases respectively).   The module $V^{2,-1}_3$ has a basis of annular tangles with $4$ strands connecting the generator to the outer boundary, and two points on the outer boundary connected to each other.  On these bases, $\rho$ has eigenvalue $1$ with multiplicity $37$, and eigenvalues $\omega$ and $\omega^2$ with multiplicity $27$.

Therefore $V^{3,1} \oplus 2 V^{3,\omega} \oplus 2 V^{3,\omega^2} \subset PABG(H)$, and 
$$PABG(H)_3  = V^{0,\delta}_3 \oplus 3 V^{0,\mu}_3 \oplus 2 V^{(0,+),0}_3 \oplus V^{2,-1}_3 
\oplus V^{3,1}_3\oplus 2 V^{3,\omega}_3 \oplus 2 V^{3,\omega^2}_3.
$$
Then we know
$$PABG(H)_4  \supset V^{0,\delta}_4 \oplus 3 V^{0,\mu}_4 \oplus 2 V^{(0,+),0}_4 \oplus V^{2,-1}_4 
\oplus V^{3,1}_4\oplus 2 V^{3,\omega}_4 \oplus 2 V^{3,\omega^2}_4.
$$ and a dimension count implies that $PABG(H)$ contains $13$ modules of the form $V^{4,\zeta}$.  

Again, we analyze the action of $\rho$ on the basis of $PABG(H)$ to determine that $\rho$ has eigenvalues $\{1, -1, i, -1\}$ with multiplicities $\{ 105, 96, 87, 87 \}$.  If we analyze the action of $\rho$ on the bases of the irreducible modules, we get that it has eigenvalues $\{1, -1, i, -1\}$ with multiplicities $\{ 102, 92, 84, 84 \}$.  (As above, $\rho$ acts by permuting the natural bases of $V^{0,?}$ and $V^{3,?}$.  However one must be a bit careful on $V^{2,-1}$, where $\rho$ has $6$ orbits with four elements, and $2$ not-quite-orbits, of the form $x \rightarrow y \rightarrow -x$.  These not-quite-orbits give eigenvectors $x - i y$ and $x + i y$ which have eigenvalues $i$ and $-i$.)

Therefore
\begin{multline*}
PABG(H)_4  = V^{0,\delta}_4 \oplus 3 V^{0,\mu}_4 \oplus 2 V^{(0,+),0}_4 \oplus V^{2,-1}_4 
\oplus V^{3,1}_4\oplus 2 V^{3,\omega}_4 \oplus 2 V^{3,\omega^2}_4 \\
\oplus 3 V^{4,1}_4 \oplus 4 V^{4,-1}_4 \oplus 3 V^{4,i}_4 \oplus 3 V^{4,-i}_4
\end{multline*}

And, as desired, we see that the subspace of $PABG(H)_4$ of low weight 4 elements with $\rho$-eigenvalue $-1$ is four-dimensional.  
\end{proof}

Of course, for the purpose of doing computations, it would be better if we could explicitly write $X^{4,-1}$ as a four-dimensional subspace of $PABG(H)$.  The symmetries of $H$, and the requirements on $X^{4,-1}$ (such as $\epsilon_i(X^{4,-1})=0$) mean that the space $X^{4,-1}$ is very symmetric.  In order to highlight these symmetries, and shorten the description, we introduce operators (some of which are not planar) on $PABG(H)_4$.  

\begin{Notation}

\begin{itemize}
\item The rotation operator $\rho$ has already been introduced:
	$$\rho=\input{Diagrams/rho4}$$
\item
The operator $\alpha$ is a symmetry of $H$ which permutes legs; Specifically, it sends $z_j \mapsto z_{j+1 \mod 3}$,  $a_j \mapsto a_{j+1 \mod 3}$,  and  $b_j \mapsto b_{j+1 \mod 3}$.  
\item
The operator $f_i$ on paths ``flips''  vertex $(i+1)$; if $v_{i}=v_{i+2}$, and $v_{i}$ is adjacent to only two vertices ($v_{i+1}$ and another, call it $w$), then 
	$$f_i (v_1 \ldots v_{i} v_{i+1} v_{i+2} \ldots v_8) = 
	- \sqrt{\frac{\lambda(v_{i+1})}{\lambda(w)}} v_1 \ldots v_{i} w v_{i+2} \ldots v_8.$$  
\end{itemize}
\end{Notation}

When paths differ in only one position, and there are exactly two vertices which can go in this position, the condition $\epsilon_i(X^{4,-1})=0$ forces the coefficient of one path to be a specific multiple of the coefficient of the other.  The purpose of the operators $f_i$ is to concisely express this relation.

\begin{Theorem} Let $t_0$, $t_1$, $t_2$ and $s$ be free variables.  If $X \subset PABG(H)_{4,+}$ is defined by
\begin{align*}
X=& \big(1-\rho + \rho^2 - \rho^3 \big) \cdot \\
&\Big(\big( t_0  + t_1 \cdot \alpha + t_2 \cdot \alpha^2 \big)\big(1+f_1 + f_3+ f_1f_3 +f_2f_1f_3 \big) 
\big(b_0cb_0cb_0cb_1c - b_0cb_0cb_0cb_2c \big)\\
&+  \big( s + s \cdot \alpha + s \cdot \alpha^2 \big) \big( 1+f_1 + f_5 + f_1f_5 \big)  
	\big(b_0cb_0cb_1cb_1c \big)  \\
&+ \big((t_0-s) + (t_1-s) \cdot \alpha + (t_2 -s) \cdot \alpha^2 \big) \big(1+f_1 \big) 
	\big(b_0cb_0cb_1cb_2c \big) \\
&+ \big((-t_0-s) + (-t_1-s) \cdot \alpha + (-t_2 -s) \cdot \alpha^2 \big) \big(1+f_1 \big) 
	\big(b_0cb_0cb_2cb_1c \big) \\
& + \big((t_1-t_2) + (t_2-t_0) \cdot \alpha +  (t_0-t_1)\cdot \alpha^2 \big) \big(b_0cb_1cb_0cb_2c \big) \\
 & +\frac{1}{2} \big( (t_2 -t_0 - t_1)  + (t_0 - t_1 - t_2) \cdot \alpha +  ( t_1 - t_2 - t_0) \cdot \alpha^2 \big)
 	\big(b_0cb_1cb_0cb_1c \big)\Big)
\end{align*}
then $X=X^{4,-1}$.
\end{Theorem}

\begin{proof}  
We must check that for all  $T \in X$,
$$\epsilon_i(T)=0 \text{ for }i=1,\ldots,8, \text{ and } \rho(T)=-T.$$

The fact that $X$ is a product of $(1-\rho + \rho^2 - \rho^3)$ means that everything in $X$ has rotational eigenvalue $-1$.

To prove that $\epsilon_i(X)=0$ for all $i$, note that the only possible states on the diagram 
$$\epsilon_i = \input{Diagrams/epsilonInLevelFour}$$ 
must have, along their inner boundary, loops $\lp p$ such that $p_i = p_{i+2}$.  
(By $p_k$ we mean the $k^{th}$ vertex of $\lp p$.)
Then $\epsilon_i(\lp p)$ is (some multiple of) the length-six loop which comes from deleting vertices $p_{i}$ and $p_{i+1}$:  $p_1 \ldots p_i p_{i+3} \ldots p_8$.  

So the loops in $X$ that could contribute a non-zero summand to  $\epsilon_i(X)$ 
are terms where the $i$ and $i+2$ vertices are identical.  
We need the contributions of such loops to cancel each other out;
requiring the coefficient of each 6-loop to be 0 gives a relation among 
those 8-loops which give the same 6-loop.  So we suppose $\lp p$ is a 
6-loop and consider three cases.

First, if vertex $i$ of $\lp p$ is $z_j$,  there is a single 8-loop $\lp q$ such 
that $\epsilon_i(\lp q)$ is a multiple of  $\lp p$, and so the coefficient of 
$\lp q$ in $X$ must be zero.  Note that every loop with non-zero 
coefficient in $X$ passes through vertex $c$ (at least) twice, hence 
cannot pass through $z_j$ twice also.  

Second, if vertex $i$ of $\lp p$ is  $a_j$ or $b_j$, there are two 8-loops $\lp q$ and $\lp r$ such that 
$\epsilon_i(\lp q)$ and $\epsilon_i(\lp r)$ are multiples of $\lp p$.  In this case, if the coefficient in $X$ of $\lp q$ is $x$ and the coefficient in $X$ of $\lp r$ is $-x\sqrt{\frac{\lambda(q_{i+1})}{\lambda(r_{i+1})}}$ if $i=1,2,3,5,6,7$ or $-x\frac{\lambda(q_{i+1})}{\lambda(r_{i+1})}$ if $i=4,8$, then the coefficient of $\lp p$ in $\epsilon_i(X)$ is $0$.  By the definition of $f_i$, $X$  has this property for all $i$ and loops whose $i$ and $i+2$ vertices are $a_j$ or  $b_j$.

Finally, if vertex $i$ is $c$, there are three 8-loops $\lp q$, $\lp r$ and $\lp s$ such that 
$\epsilon_i(\lp q)$, $\epsilon_i(\lp r)$ and $\epsilon_i(\lp s)$ are multiples of  $\lp p$.  One can check that in this case, if the coefficient in $X$ of $\lp q$ is $x$, the coefficient in $X$ of $\lp r$ is $y$ and the coefficient in $X$ of $\lp s$ is 
$-x\sqrt{\frac{\lambda(q_{i+1})}{\lambda(s_{i+1})}} - 
y\sqrt{\frac{\lambda(r_{i+1})}{\lambda(s_{i+1})}}$, 
then the coefficient of $\lp p$ in $\epsilon_i(X)$ is $0$.  In this case, all vertices adjacent to $c$, namely $b_0$, $b_1$ and $b_2$, have the same Perron-Frobenius eigenvector entries, so the coefficient in $X$ of $\lp s$ is 
$-x - y$.  By inspection, $X$ as stated in the theorem has this property for all $i$ and loops whose $i$ and $i+2$ vertices are $c$.

Therefore, $\epsilon_i(X)=0$, and $X=X^{4,-1}$.
\end{proof}

When we add the requirement that for all $T \in X$, $T^*=T$, we get that $t_0, t_1, t_2 \in \mathbb{R}$ and $s \in i\mathbb{R}$. 

Which, if any, element in $X^{4,-1}$ could generate a Haagerup planar algebra?

\begin{Theorem} The only elements of $X^{4,-1} \subset PABGH_4(H)$ which could generate a Haagerup planar algebra are multiples of
\begin{align*}
T=& t(1-\rho + \rho^2 - \rho^3) \cdot ( 1  +   \alpha +   \alpha^2) \cdot \\
&((1+f_2 + f_4+ f_2f_4 +f_3f_2f_4) 
(b_0cb_0cb_0cb_1c - b_0cb_0cb_0cb_2c)\\
&+ i \sqrt{ \frac{1+\sqrt{13}}{2}} \cdot ( 1+f_2 + f_6 + f_2f_6)  (b_0cb_0cb_1cb_1c)  \\
&+(1- i \sqrt{ \frac{1+\sqrt{13}}{2}}) \cdot (1+f_2) (b_0cb_0cb_1cb_2c) \\
&+(-1- i \sqrt{ \frac{1+\sqrt{13}}{2}})\cdot (1+f_2) (b_0cb_0cb_2cb_1c) \\
 & -\frac{1}{2}
 (b_0cb_1cb_0cb_1c)) \\
\end{align*}
\end{Theorem}

\begin{proof}  First note that $T \in X^{4,-1}$.  

Since $\dim (P_{0,+})=\dim (TL_{0,+})=1$ for a (still hypothetical) Haagerup planar algebra $P$,  we look for $T \in X$ such that $\tr{T^2} \in TL_{0,+} \subset PABG(H)_{0,+}$.  The single basis element of $TL_+$, when interpreted as an element of $PABG(H)_+$, is $z_0+z_1+z_2+b_0+b_1+b_2$.  Thus,
 the coefficients in $\tr{T^2}$ of $z_0$, $z_1$ and $z_2$ must all be equal, which quickly gives us that 
 $t_0=t_1=t_2$ (and from now on we will denote $t_0$, $t_1$ and $t_2$ by $t$).
 
 To find the relation between $s$ and $t$, we use the fact that the coefficients of $z_0$ and $b_0$ are equal.  Let $\coeff{x}{y}$ denote the coefficient in $x$ of $y$.  Then we calculate that 
 
\begin{align*}
\coeff{\tr{T^2}}{b_0}=&
t ^2(8
+ 16 ( \frac{1+\sqrt{13}}{2}  )
+8 ( \frac{1+\sqrt{13}}{2}  )^2
+6 ( \frac{1+\sqrt{13}}{2}  )^2 ( \frac{3+\sqrt{13}}{2}  ))\\
&+
\vert s \vert^2(4
+8 ( \frac{1+\sqrt{13}}{2}  )
+4 ( \frac{1+\sqrt{13}}{2}  )^2)\\
&+
(t^2+ \vert s \vert^2)(8
+8 ( \frac{1+\sqrt{13}}{2}  )) +2t^2\\
=&
t^2(109+31\sqrt{13}) +\\
&\vert s \vert ^2(34 + 10\sqrt{13})
\end{align*}

 and
 
 \begin{align*}
\coeff{\tr{T^2}}{z_0}=&
t^2(2 ( \frac{1+\sqrt{13}}{2}  )^2 ( \frac{3+\sqrt{13}}{2}  )^3) \\
= & t^2(191 + 53 \sqrt{13})
\end{align*}

from which it follows that, as $s$ is purely imaginary, $$s=\pm i t \sqrt{ \frac{1 + \sqrt{13}}{2}}.$$
 This tells us what $T$ is, up to a choice of normalization.  In order to be consistent with \cite{quadratic}, we choose $t=\frac{1}{3}(4-\sqrt{13})$ so that  $\ip{T,T}=Z(T^2)=[5]=3+\sqrt{13}$.
\end{proof}

%% file: Diagrams/AnnRelns/Cap.tex
\begin{tikzpicture}[annular]
	\clip (0,0) circle (2cm);

	\filldraw[shaded] (0,0)--(158:1cm) .. controls ++(158:1cm) and ++(112:1cm) .. (112:1cm) -- (0,0);
	\filldraw[shaded] (0,0)--(-158:2cm) arc (-158:-112:2cm) -- (0,0);
	\filldraw[shaded] (0,0)--(22:2cm) arc (22:68:2cm) -- (0,0);
	\filldraw[shaded] (0,0)--(-22:2cm) arc (-22:-68:2cm) -- (0,0);
		
	\draw[ultra thick] (0,0) circle (2cm);
	
	\node at (0,0)  [Tbox] (T) {$T$};
	\node at (T.180) [left] {$\star$};
	\node at (180:2cm) [right] {$\star$};
\end{tikzpicture}
=
\begin{tikzpicture}[annular]
	\clip (0,0) circle (2cm);

	\filldraw[shaded] (0,0)--(158:2cm) arc (158:22:2cm) -- (0,0) -- (68:1cm) .. controls ++(68:1cm) and ++(112:1cm) .. (112:1cm) -- (0,0);
	\filldraw[shaded] (0,0)--(-158:2cm) arc (-158:-112:2cm) -- (0,0);
	\filldraw[shaded] (0,0)--(-22:2cm) arc (-22:-68:2cm) -- (0,0);
		
	\draw[ultra thick] (0,0) circle (2cm);
		
	\node at (0,0)  [Tbox] (T) {$T$};
	\node at (T.180) [left] {$\star$};
	\node at (180:2cm) [right] {$\star$};
\end{tikzpicture}
=0

%% file: Diagrams/AnnRelns/Rotate.tex
\begin{tikzpicture}[annular]
	\clip (0,0) circle (2cm);

	\filldraw[shaded] (0,0)--(158:2cm) arc (158:112:2cm) -- (0,0);
	\filldraw[shaded] (0,0)--(-158:2cm) arc (-158:-112:2cm) -- (0,0);
	\filldraw[shaded] (0,0)--(22:2cm) arc (22:68:2cm) -- (0,0);
	\filldraw[shaded] (0,0)--(-22:2cm) arc (-22:-68:2cm) -- (0,0);
		
	\draw[ultra thick] (0,0) circle (2cm);
		
	\node at (0,0)  [Tbox] (T) {$T$};
	\node at (T.180) [left] {$\star$};
	\node at (180:2cm) [right] {$\star$};
\end{tikzpicture}
=
-
\begin{tikzpicture}[annular]
	\clip (0,0) circle (2cm);

	\filldraw[shaded] (0,0)--(158:2cm) arc (158:112:2cm) -- (0,0);
	\filldraw[shaded] (0,0)--(-158:2cm) arc (-158:-112:2cm) -- (0,0);
	\filldraw[shaded] (0,0)--(22:2cm) arc (22:68:2cm) -- (0,0);
	\filldraw[shaded] (0,0)--(-22:2cm) arc (-22:-68:2cm) -- (0,0);
		
	\draw[ultra thick] (0,0) circle (2cm);

	\node at (0,0)  [Tbox] (T) {$T$};
	\node at (T.270) [below] {$\star$};
	\node at (180:2cm) [right] {$\star$};
\end{tikzpicture}

%% file: Diagrams/epsilonInLevelFour.tex
\begin{tikzpicture}[annular]
	\clip (0,0) circle (2cm);

	\filldraw[shaded] (0,0) .. controls ++(170:2cm) and ++(100:2cm) .. (0,0);
	\filldraw[shaded] (-158:4cm)--(0,0)--(-112:4cm);
	\draw[shaded] (68:4cm)--(0,0)--(22:4cm);
	\draw[shaded] (-68:4cm)--(0,0)--(-22:4cm);
		
	\draw[ultra thick] (0,0) circle (2cm);
	
	\node at (0,0)  [empty box] (T) {};
\end{tikzpicture}

%% file: text/Quadratic.tex
A shaded planar tangle is called {\em quadratic} if it has two inner disks.  In this section, we prove some important quadratic skein relations involving $T$, using linear algebra and the following six facts.
   
\begin{Lemma}\label{powersofT}
Recall that $\rho^{1/2}$ is the shading-changing rotation morphism 
$$PABG(H)_{4,+} \rightarrow PABG(H)_{4,-}.$$
 \begin{itemize}
 \item $Z(T^2)=3+\sqrt{13}$
 \item $Z(T^3)=0$
 \item $Z(T^4)=3+\sqrt{13}$
 \item $Z((\rho^{1/2} T)^2)=-3-\sqrt{13}$
 \item $Z((\rho^{1/2} T)^3)=i\sqrt{\frac{8(3+\sqrt{13})}{3}}$
 \item $Z((\rho^{1/2}T)^4)=\frac{17+3\sqrt{13}}{3}$ 
 \end{itemize}
 \end{Lemma}
 
\begin{proof}
 The proofs of these facts are straightforward once the above expression for $T$ is written out explicitly and converted into a sum of matrices, as is done in appendix \ref{traces}.  
\end{proof}

 This small number of facts enables us to calculate a great deal about $T$.  In this section we will prove the following theorem about the existence of relations at levels 4, 5 and 6.
For our purposes, it is not the exact values of these relations that are important, but simply that there {\em is} a relation between, for instance, two Ts joined by three strands and diagrams containing at most one T.   For example, in the proof of Theorem \ref{IsSubfactor}, we need to know that there is a relation allowing us to replace two $T$s connected with three strings by a sum of diagrams involving one or no $T$s; but the coefficients in this relation are irrelevant to the proof.

 \begin{Theorem}\label{skeinrelations}
 The following quadratic skein relations hold in $PABG(H)_4$:
 \begin{enumerate}
\item[(R4), (R4')]  $\input{Diagrams/QuadRelns/4}  $
\item[(R5)]  $\input{Diagrams/QuadRelns/5simplified} $
\item[(R6)]  $\input{Diagrams/QuadRelns/6withN} $
\end{enumerate}
 \end{Theorem}

As mentioned above, the statements in Theorem \ref{skeinrelations} are the right precision for the rest of this paper.  However,  in order to prove them, we need more explicit versions where the projections of these elements onto the relevant subspaces has been calculated.  So we prove Theorem \ref{skeinrelations} as a corollary of 

\begin{Theorem}
	\begin{enumerate}
	 \item[(R4)] $T^2=\f{4}$
	 \item[(R4')] $\joint{4}{T}{T}= - \rho^{-1/2} \f{4} - i \sqrt{\frac{8}{3[5]}} T $
	 \item[(R5)] $\join{3}{T}{T}=\frac{[5]}{[6]}\f{5}  - i\sqrt{\frac{8(3+\sqrt{13})}{3}}\hat{w}_{10} $
	 \item[(R6)] $\joint{2}{T}{T} = - \frac{[5]}{[7]} \rho^{- 1/2}\f{6} + \frac{ \sqrt{ 4-\sqrt{13} }}{2} N + i\sqrt{\frac{8(3+
	 	\sqrt{13})}{3}}( - \hat{w}_{4,5} +   \hat{w}_{10,11})$
	 \end{enumerate}
\end{Theorem}

We will define later $\hat{w}_i$ and $\hat{w}_{i,j}$ (but before we prove (R5) and (R6)); the rest of the notation is as follows.

\begin{Notation}
Let $\join{n}{A}{B}$ be the tangle with an $A$ box joined to a $B$ box by the $n$ strings counterclockwise of $A$'s star and the $n$ strings clockwise of $B$'s star.  There is also a `twisted' version of this, $\joint{n}{A}{B}=\rho^{-1/2}\join{n}{\rho^{1/2}A}{\rho^{1/2}B}$:

$$ \join{n}{A}{B} = 
\input{Diagrams/join},  \quad
\joint{n}{A}{B} =
\input{Diagrams/jointwist} 
$$

Let $N$ be the part of $\join{2}{T}{T}$ which is new, i.e.  
$$
N
=
(1-\text{Proj}_{TL}-\text{Proj}_{ATL_6(T)})
\left(
\input{Diagrams/QuadRelns/New}
\right).
$$

\end{Notation}

%% file: text/level4.tex
 Although (R4) and (R4') follow from (R5), it is useful to prove them independently; they will give insight into how the proofs of (R5) and (R6) work.  All of these relations will be proved in a similar way, via Bessel's inequality.  Bessel's inequality says that for any $v \in V$, $W \subset V$, we have $\left\Vert v \right\Vert^2 \geq \left\Vert \proj{W}{v} \right\Vert^2$, with equality if and only if $v \in W$.
 
\begin{proof}[Proof of (R4)] 
To show 
$$T^2=\proj{TL}{T^2}=f^{(4)},$$ 
we will first show 
$$ \left\Vert T^2 \right\Vert^2 = 3+\sqrt{13}$$
and then
$$ \left\Vert \proj{TL}{T^2} \right\Vert^2 =  \left\Vert f^{(4)} \right\Vert^2 = 3+\sqrt{13}.$$

First, from \ref{powersofT}, we know
$$\Norm{T^2}^2= \ip{T^2,T^2}=Z(T^4)=3+\sqrt{13}.$$
Next, recall that $B(TL_4)$ is the usual picture basis of $TL_4$, and compute
\begin{align*}
 \left\Vert \proj{TL}{T^2} \right\Vert^2 
 	& =\Norm{\sum_{\beta \in B(TL_4)} \ip{T^2,\beta} \hat{\beta}}^2 
	 =\Norm{ \ip{T^2,\id} \frac{f^{(4)}}{[5]} }^2 \\
 	& = \left\Vert  \tr{T^2} \frac{f^{(4)}}{[5]}  \right\Vert^2 \\
 	& = \left\Vert  f^{(4)}      \right\Vert^2 = \ip{f^{(4)},f^{(4)}} = 3+\sqrt{13}
\end{align*}
\end{proof}
 
The proof of (R4') is analogous and only slightly more complicated.
 
 \begin{proof} [Proof of (R4')]
 To show
 $$\joint{4}{T}{T} = \proj{\spn{TL_4, T}}{\joint{4}{T}{T}}=  -\rho^{-1/2} \f{4} - i\sqrt{\frac{8}{3[5]}} T,$$
 we show
 $$\Norm{\joint{4}{T}{T}}^2 = \frac{17+3\sqrt{13}}{3}$$
 and
 \begin{align*}
 \Norm{\proj{\spn{T,TL_4}}{\joint{4}{T}{T}}}^2 & = \Norm{ -\rho^{-1/2} \f{4} - i \sqrt{\frac{8}{3[5]}} T }^2 =\frac{17+3\sqrt{13}}{3}.
 \end{align*}
 
 First, 
 \begin{align*}
 \Norm{\joint{4}{T}{T}}^2 & =\ip{\joint{4}{T}{T},\joint{4}{T}{T}}=Z((\rho^{1/2}T)^4)=\frac{17+3\sqrt{13}}{3}
 \end{align*}
 Next,
  \begin{align*}
 \Norm{\proj{\spn{T,TL_4}}{\joint{4}{T}{T}}}^2 & = \Norm{\proj{TL}{\joint{4}{T}{T}} + \proj{T}{\joint{4}{T}{T}}}^2 \\
 & = \Norm{  \sum_{\beta \in B(TL_4)} \ip{\joint{4}{T}{T}, \beta} \hat{\beta} + \frac{\ip{\joint{4}{T}{T},T}}{\ip{T,T}} T } ^2 \\
 & = \Norm{ \ip{\joint{4}{T}{T}, \input{Diagrams/e1e2e3e4}} \widehat{\input{Diagrams/e1e2e3e4}} }^2 
 + \Norm{ \frac{ -Z((\rho^{1/2}T)^3)}{Z(T^2)}T }^2 \\
 & = \Norm{-[5] \frac{\rho^{-1/2}{f^{(4)}}}{[5]}}^2 + \Norm{ -i \sqrt{ \frac{8}{3[5]}} T }^2 \\
 & = [5]+\frac{8}{3[5]}[5] \\
 & =\frac{17+3\sqrt{13}}{3}
 \end{align*}
  \end{proof}

%% file: Diagrams/e1e2e3e4.tex
\begin{tikzpicture}[scale=.2]
 	\draw (1,1)--(3,-1);
	\draw (2,1)--(4,-1);
	\draw (1,-1) arc (180:0:.5cm);
	\draw (3,1) arc (-180:0:.5cm);
 \end{tikzpicture}

%% file: text/DualBasis.tex
 In the above proofs of (R4) and (R4'), calculating projections onto $TL_4$ was easy, even though $TL_4$ is 14 dimensional, and its natural basis (of non-crossing pairings) is not orthogonal.  This was because $\join{4}{T}{T}$ and $\joint{4}{T}{T}$ each had non-zero inner product with only one $TL$ basis elements, whose dual was easily expressed in terms of $f^{(4)}$.  

In order to prove (R5) and (R6), we need to be able to calculate projections onto $ATL_5(T)$ and $ATL_6(T)$.  The natural bases of these consists of tangles with one $T$ and some caps on the outer boundary.  The simplest way to calculate projections onto these spaces is to learn more about the dual basis to the natural basis, and that is the purpose of this section.  
  
In the usual basis of $TL$, the Jones-Wenzl projection is dual to the identity.  The question of explicitly writing $f^{(n)}$ as a sum of basis elements was solved in \cite{MR1446615} and later, in a more pictorally natural form, in \cite{morrison}.   Duals to other Temperley-Lieb elements can sometimes be built up from Jones-Wenzl idempotents, and we will do something similar in  $ATL_5(T)$ and $ATL_6(T)$.

\begin{Defn}
By $w_{i_1, i_2, \ldots , i_n}(T)$, we denote the element of $ATL_{4+n}(T)$ with caps originating from the $i_1, i_2, \ldots , i_n$ external boundary points (counted clockwise from the external star).  Since $T$ is a rotational eigenvector, the position of the internal star affects $w_{i_1, i_2, \ldots , i_n}(T)$ only us to a choice of sign.  When possible, it is convenient to place the star on $T$ is in the same region as the star on the outer boundary.  Otherwise, we reserve the right to specify the position of the star later.
\end{Defn}

Let $\mathcal{S}$ be the set of $n$-element subsets of $\{1, \ldots, 2n+8\}$.  Then $\{ w_I(T) \}_{I \in \mathcal{S}}$ is a basis for $ATL_{4+n}(T)$.  For example, let $B(ATL_5(T))$ be the basis for $ATL_5(T)$ given by
$$\input{Diagrams/ATLFiveEGs/basis}$$
and let $B(ATL_6(T))$ be the basis for $ATL_6(T)$ given by
$$\input{Diagrams/ATLSixEGs/nnplusone}$$
$$\input{Diagrams/ATLSixEGs/nnplustwo}$$
\vspace{6pt}
$$\hdots$$
\vspace{6pt}
$$\input{Diagrams/ATLSixEGs/nnplussix}.$$
\vspace{6pt}

$ATL_5(T)$ is simple enough that we can guess the dual basis to $B(ATL_5(T))$.

\begin{Defn} The element $\alpha_i(T) \in ATL_5(T)$ is $w_i(T)$ with an $f^{(10)}$ wrapped around it, in the only possible non-zero way.  For example, $\alpha_1(T)$ is
$$\input{Diagrams/ATLFiveEGs/SampleDual}$$
\end{Defn}

\begin{Lemma}
Up to renormalization, $\alpha_i(T)$ is dual to $w_i(T)$:
$$\hat{w}_i(T) = \frac{\alpha_i(T)}{\ip{w_i(T), \alpha_i(T)}}$$
\end{Lemma}

\begin{proof}
Any other $w_j(T)$ gives zero when paired with $\alpha_i(T)$, as its boundary cap lands on the Jones-Wenzl idempotent.  
\end{proof}

It would be nice if we could guess the basis dual to $B(ATL_6(T))$. 

\begin{Defn} The element $\alpha_{i,i+1}(T) \in ATL_{6}(T)$ is  $w_{i,i+1}(T)$ with an $\f{12}$ wrapped around it, and $\alpha_{i,j}(T)$ is $w_{i,j}(T)$ with two Jones-Wenzl projections wrapped around it, in the only possible non-zero way.  For example, $\alpha_{4,5}(T)$ and $\alpha_{6,12}(T)$ are 
$$\input{Diagrams/ATLSixEGs/SampleDual}$$
\end{Defn}

\begin{Lemma}
Up to renormalization, $\alpha_{i,i+1}(T)$ is dual to $w_{i,i+1}(T)$:
$$ \hat{w}_{i,i+1}(T) = \frac{\alpha_{i,i+1}(T)}{\ip{w_{i,i+1}(T), \alpha_{i,i+1}(T)}}.$$
\end{Lemma}

\begin{proof}
Any $w_J(T)$ ($J \neq \{i,i+1\}$) gives zero when paired with $\alpha_{i,i+1}(T)$, as one of its boundary caps lands on the Jones-Wenzl idempotent.  
\end{proof}

For $w_{i,j}(T)$ with $j \neq i +1$, unfortunately $\alpha_{i,j}(T)$ has non-zero inner product with not just $w_{i,j}(T)$ but also $w_{i-1,i}(T)$ and $w_{j-1,j}(T)$.  Fortunately, we already know the duals to this second kind of basis element.  So to find $\hat{w}_{i,j}(T)$, we just subtract the part of $\alpha_{i,j}(T)$ that comes from 
$\hat{w}_{i-1,1}(T)$ and $\hat{w}_{j-1,j}(T)$
\begin{Lemma}
$$	\hat{w}_{i,j}(T) =\frac{\alpha_{i,j}(T) - \ip{\alpha_{i,j}(T),w_{i-1,i}(T)} \hat{w}_{i-1,i}(T)
	- \ip{\alpha_{i,j}(T),w_{j-1,j}(T)} \hat{w}_{j-1,j}(T)}{\ip{\alpha_{i,j}(T),w_{i,j}(T)}}.
$$
\end{Lemma}

Unfortunately, calculating the inner products required to correctly normalize is a somewhat involved process, as it involves knowing coefficients of various $TL$ basis elements in Jones-Wenzl idempotents.  We work out explicit coefficients in the above formulae only for the the cases we will need in order to prove relations (R5) and (R6).

\begin{Lemma}
\begin{enumerate}
	\item $\hat{w}_i(T)=\frac{1}{\sqrt{2(5+\sqrt{13})}} \cdot \alpha_i(T)$
	\item $\Norm{\hat{w}_i(T)}^2=\frac{1}{\sqrt{2(5+\sqrt{13})}}$
\end{enumerate}
\end{Lemma}

\begin{proof}
\begin{enumerate}
\item
Since
$$\hat{w}_{i}(T)=\frac{ \alpha_{i}(T)}{\ip{\alpha_{i}(T), w_{i}(T)}},$$ 
we need to compute 
$\ip{\alpha_{i}(T), w_{i}(T)}=$
$$\input{Diagrams/InnerProducts/ATLfiveAlphaAlpha}$$

(We've drawn the tangle for $\ip{\alpha_6(T),w_6(T)}$, but moving both stars together and/or switching the shading won't change the value of the tangle; so $\ip{\alpha_i(T),w_i(T)}=\ip{\alpha_6(T),w_6(T)}$.)

The $\beta$ which give a non-zero number when inserted in the above tangle are 
$$\input{Diagrams/TLTen/NonZeroBeta}.$$  
The following table gives their coefficients in $f^{(10)}$, and the value of the tangle with $\beta$ inserted.

\begin{center}
\begin{tabular}{c | c | c }
	$\beta \in B(TL_{10})$ & $\coeff{f^{(10)}}{\beta}$ & value of tangle \\
	\hline
	\rule[-3mm]{0mm}{7mm} \input{Diagrams/TLTen/One} & $1$ & $[2][5]$ \\
	\rule[-3mm]{0mm}{6mm} \input{Diagrams/TLTen/Two} & $-\frac{[9]}{[10]}$ & $[5]$ \\
	\rule[-3mm]{0mm}{6mm} \input{Diagrams/TLTen/Three} & $-\frac{[9]}{[10]}$ & $[5]$ \\
	\rule[-3mm]{0mm}{6mm} \input{Diagrams/TLTen/Four} & $-\frac{1}{[10]}$ & $-[5]$ \\
	\rule[-4mm]{0mm}{7mm} \input{Diagrams/TLTen/Five} & $-\frac{1}{[10]}$ & $-[5]$ \\
\end{tabular}
\end{center}

\vspace{6pt}
\noindent
Summing the products of the coefficients and values of the tangle, we get that $\ip{\alpha_{i}(T), w_{i}(T)}=\sqrt{2(5+\sqrt{13})}$.
\item 
Because $f^{(10)} \cdot f^{(10)}=f^{(10)}$, we have $\ip{\alpha_i(T),\alpha_i(T)}=\ip{\alpha_i(T),w_i(T)}=\sqrt{2(5+\sqrt{13})}$.  So
\begin{align*}
	\ip{\hat{w}_i(T),\hat{w}_i(T)} & = \frac{1}{2(5+\sqrt{13})} \ip{\alpha_i(T),\alpha_i(T)} \\ 
	& = \frac{1}{\sqrt{2(5+\sqrt{13})}}
\end{align*}
\end{enumerate}
\end{proof}

The following facts about $ATL_6(T)$ can be proved in a similar, although more involved, fashion.  The casual reader may want to take these on faith and go right to Section \ref{levels56} (or even Section \ref{PATisHaagerup}).

\begin{Lemma}\label{dualbasis}
\begin{enumerate}
\item $ \hat{w}_{i,i+1}(T)=\frac{19+5\sqrt{13}}{144} \cdot \alpha_{i,i+1}(T)$
\item\label{normiiplusone}  $ \| \hat{w}_{i,i+1}(T) \|^2=\frac{19+5\sqrt{13}}{144} $
\item $\hat{w}_{i,i+6}(T)=\frac{11+\sqrt{13}}{36} \cdot \alpha_{i,i+6} $
\item\label{normiiplussix} $\| \hat{w}_{i,i+6}(T) \|^2 = \frac{11+\sqrt{13}}{36}$
\item\label{decompiiplusone} $\ip{\hat{w}_{4,5}(T), \hat{w}_{10,11}(T)} = \ip{\hat{w}_{10,11}(T),\hat{w}_{4,5}(T)} = \frac{17-5\sqrt{13}}{144}$
\item\label{decompiiplussix} $\ip{\hat{w}_{4,5}(T),\hat{w}_{6,12}(T)} = \ip{\hat{w}_{10,11}(T),\hat{w}_{6,12}(T)} = \frac{-3}{[6](11-\sqrt{13})}
$
\end{enumerate}
\end{Lemma}

\begin{proof}
\begin{enumerate}

\item Since
$$\hat{w}_{i,i+1}(T)=\frac{ \alpha_{i,i+1}(T)}{\ip{\alpha_{i,i+1}(T), w_{i,i+1}(T)}},$$ 
we need to compute 
$\ip{\alpha_{i,i+1}(T), w_{i,i+1}(T)}=$
$$\input{Diagrams/InnerProducts/AlphaOneAlphaOne}$$

(Again, the stars and shading can be changed to make this compute any $\ip{\alpha_{i,i+1}(T),\alpha_{i,i+1}(T)}$ without changing the result.)

Table \ref{TLTwelve} is a list of the $\beta$ in the basis of $TL_{12}$ that don't cap off either $T$; the coefficient of $\beta$ in $f^{(12)}$ (calculated as described in \cite{morrison}); and the value of the above tangle with $\beta$ inserted.  When we sum (over $\beta \in B(TL_{12})$) the product of the coefficient and the value of the tangle, we get that $\ip{\alpha_{i,i+1}(T), w_{i,i+1}(T)}=76-20\sqrt{13}$.  Thus $\hat{w}_{i,i+1}(T)=\frac{19+5\sqrt{13}}{144} \cdot \alpha_{i,i+1}(T)$.

\input{Tables/TLAlphaOneAlphaOne}

\item Because $f^{(12)} \cdot f^{(12)} = f^{(12)}$, we have $ \ip{\alpha_{i,i+1}(T), \alpha_{i,i+1}(T)} = \ip{\alpha_{i,i+1}(T), w_{i,i+1}(T)}=76 - 20 \sqrt{13}$.  So
\begin{align*}
	\ip{\hat{w}_{i,i+1}(T), \hat{w}_{i,i+1}(T)} & = (\frac{19-5\sqrt{13}}{144})^2 \ip{\alpha_{i,i+1}(T), \alpha_{i,i+1}(T)} \\
	& = \frac{19-5\sqrt{13}}{144}.
\end{align*}

\item Before proving this, we should note that much of this section is applicable in the case where $T$ is replaced by some other generator of an irreducible $TL$ module.  However, the fact that $\hat{w}_{i,i+6}$ is a multiple of $\alpha_{i,i+6}$ is specific to the case where $\rho(T)=-T$.

 Recall that
\begin{multline*}
	\hat{w}_{i,i+6}(T) =\frac{1}{\ip{w_{i,i+6}(T),\alpha_{i,i+6}(T)}} \cdot ( \alpha_{i,i+6}(T) - \ip{\alpha_{i,i+6}(T),w_{i-1,i}(T)} \hat{w}_{i-1,i}(T) \\
	- \ip{\alpha_{i,i+6}(T),w_{i+5,i+6}(T)} \hat{w}_{i+5,i+6}(T)),
\end{multline*}
and in this case, we will show $\ip{\alpha_{i,i+6}(T),w_{i-1,i}(T)}=0$ and $\ip{\alpha_{i,i+6}(T),w_{i+5,i+6}(T)} =0$ so
$$ \hat{w}_{i,i+6}(T)=\frac{11+\sqrt{13}}{36} \cdot \alpha_{i,i+6}. $$

To prove this, we first calculate that (up to a sign depending on the position of the stars, which we have omitted)
$\ip{\alpha_{i,i+6}(T), w_{i-1,i}(T)}=$
$$\input{Diagrams/InnerProducts/AlphaOneAlphaSix}$$
and we see that the only elements of the basis of $TL_6$ which don't cap off a $T$ when inserted for $\beta_1$ are 
$$\input{Diagrams/TLsix/One} \, , \, \input{Diagrams/TLsix/Two}, \input{Diagrams/TLsix/Three} \, , \, \input{Diagrams/TLsix/Four} \, , \, \input{Diagrams/TLsix/Five} \, , \text{ and } \input{Diagrams/TLsix/Six}$$
 (which, respectively, have coefficients $1$, $\frac{[5]}{[6]}$,  $\frac{[4]}{[6]}$,  $\frac{[1]}{[6]}$,  $\frac{[2]}{[6]}$, and  $\frac{[3]+[5]}{[5][6]}$ in $f^{(6)}$).  Vertical reflections of these are the only basis elements which don't cap off a $T$ when inserted for $\beta_2$ (and the reflections have the same coefficients).  

Table \ref{TLSixFirst} gives the value of $\prod_{i} \operatorname{Coeff}_{f^{(6)}} (\beta_i)$ times the value of the   tangle with $\beta_1$ and $\beta_2$ inserted.  Since these sum to $0$, we get $\ip{\alpha_{i,i+6}(T), w_{i-1,i}(T)}=\pm 0=0$.

\input{Tables/TLAlphaOneAlphaSix}

Similarly, we find $\ip{\alpha_{i,i+6}(T), w_{i+5,i+6}(T)}=0$.  

Now, we need to do a similar calculation to find $\ip{w_{i,i+6}(T),\alpha_{i,i+6}(T)}=$
$$\input{Diagrams/InnerProducts/AlphaSixAlphaSix}$$
(Again, the stars and shading can be changed to make this compute any $\ip{w_{i,i+6}(T),\alpha_{i,i+6}(T)}$ without changing the result.)

The only elements of the basis of $TL_6$ which don't cap off a $T$ when inserted for $\beta_1$ or $\beta_2$ are 
$$\input{Diagrams/TLsix/One} \, , \, \input{Diagrams/TLsix/Two} \, , \, \input{Diagrams/TLsix/FourReflect} \, , \, \input{Diagrams/TLsix/Four} \, , \, \input{Diagrams/TLsix/TwoReflect} \, , \text{ and }\input{Diagrams/TLsix/Seven}$$
 (which, respectively, have coefficients $1$, $\frac{[5]}{[6]}$,  $\frac{[1]}{[6]}$,  $\frac{[1]}{[6]}$,  $\frac{[5]}{[6]}$, and  $\frac{[2]+[4][5]}{[5][6]}$ in $f^{(6)}$).

\input{Tables/TLAlphaSixAlphaSix}

Table \ref{TLSixSecond} gives the value of $\prod_{i} \operatorname{Coeff}_{f^{(6)}} (\beta_i)$ times the inner product diagram with $\beta_1$ and $\beta_2$ inserted.  When we sum these, we see $\ip{w_{i,i+6}(T),\alpha_{i,i+6}(T)}=\frac{11-\sqrt{13}}{3}$.

Therefore 
$$	\hat{w}_{i,i+6}(T)  = \frac{3}{11-\sqrt{13}} \cdot \alpha_{i,i+6}(T) = \frac{11+\sqrt{13}}{36} \cdot \alpha_{i,i+6}(T) 
$$

\item
Since the tangle for $\ip{ \alpha_{i,i+6}(T), \alpha_{i,i+6}(T)}$ has the Jones-Wenzl idempotents next to each other, 
$\ip{ \alpha_{i,i+6}(T), \alpha_{i,i+6}(T)} = \ip{ \alpha_{i,i+6}(T), w_{i,i+6}(T)} $.  Therefore
\begin{align*}\| \hat{w}_{i,i+6}(T) \|^2 & =  (\frac{11+\sqrt{13}}{36} )^2 \ip{ \alpha_{i,i+6}(T), \alpha_{i,i+6}(T)} \\
& =  (\frac{11+\sqrt{13}}{36} )^2 \ip{ \alpha_{i,i+6}(T), w_{i,i+6}(T)} \\
&= \frac{11+\sqrt{13}}{36} 
\end{align*}

\item Our goal is to show
$$\ip{\hat{w}_{4,5}(T), \hat{w}_{10,11}(T)} = \frac{17-5\sqrt{13}}{144}.$$
We have
$$\ip{\hat{w}_{4,5}(T), \hat{w}_{10,11}(T)} =  \sum_I \ip{\hat{w}_{4,5}(T) , \coeff{\hat{w}_{10,11}}{ w_I} w_I} =\coeff{\hat{w}_{10,11}}{ w_{4,5}}$$
and so we want to know which TL pictures can be inserted in the rectangle below to produce a multiple of $w_{4,5}$, what that multiple is, and what their coefficients in $f^{(12)}$ are:
$$\input{Diagrams/InnerProducts/AlphaFourFiveAlphaTenEleven}$$
One can check that the only pictures which don't give $0$ when inserted into the above tangle (and have the bottom caps in the required positions) are:
$$\input{Diagrams/TLTwelve/PicsGivingFiveSix/One} \, , \,
\input{Diagrams/TLTwelve/PicsGivingFiveSix/OneReflect}\, , \,
\input{Diagrams/TLTwelve/PicsGivingFiveSix/Two} \, , \,
\input{Diagrams/TLTwelve/PicsGivingFiveSix/Three}\, , \text{ and } 
\input{Diagrams/TLTwelve/PicsGivingFiveSix/ThreeReflect} \, . $$
In $f^{(12)}$, these have coefficients
$\frac{[5][6]}{[11][12]}$, $\frac{[5][6]}{[11][12]}$, $\frac{[6]^2}{[11][12]}$, $\frac{[5][6]( [2] + [10] + [12])}{[10][11][12]}$ and $\frac{[5][6]( [2] + [10] + [12] )}{[10][11][12]}$ respectively, and when inserted in the above tangle, give $w_{4,5}$ times $1$, $1$, $- \delta$, $-1$, and $-1$.  

We add up the product of all these, multiply by $\frac{19+5\sqrt{13}}{144}$ (the multiple of $\alpha_{i,i+1}$ which gives $\hat{w}_{i,i+1}$), and get

\begin{align*}
\ip{\hat{w}_{4,5}(T), \hat{w}_{10,11}(T)} & = 
\frac{19+5\sqrt{13}}{144} \cdot \frac{[6]}{[11][12]} \left( 2[5]+[2][6]-2 \frac{[5]}{[10]}([2] + [10] + [12] ) \right) \\
& = \frac{17-5\sqrt{13}}{144}
\end{align*}

\item

We want to find
$$ \ip{\hat{w}_{4,5}(T), \hat{w}_{6,12}(T)} = \coeff{\hat{w}_{6,12}}{w_{4,5}}$$ 
and so we want to know which TL pictures can be inserted in the rectangles below to produce a multiple of $w_{4,5}$, what that multiple is, and what their coefficients in $f^{(6)}$ are:
$$\input{Diagrams/InnerProducts/AlphaFourFiveAlphaSixTwelve}$$
One can check that the only pair of pictures which doesn't cause the above to evaluate to zero is
$$\input{Diagrams/TLsix/Four} \text{ and } \input{Diagrams/TLsix/One}\, .$$
These have coefficients $\frac{1}{[6]}$ and $1$ in $f^{(6)}$, and when inserted, give $w_{4,5}$.

Multiplying $\frac{1}{[6]}$ by $\frac{11+\sqrt{13}}{36}$ (the multiple of $\alpha_{i,i+6}$ which gives $\hat{w}_{i,i+6}$) gives

$$\ip{\hat{w}_{4,5}(T), \hat{w}_{6,12}(T)} = \sqrt{\frac{5+\sqrt{13}}{2}} \cdot \frac{22-2\sqrt{13}}{9}$$

\end{enumerate}
\end{proof}

%% file: Diagrams/ATLFiveEGs/basis.tex
w_{1}=
\begin{tikzpicture}[ATLsix]
	\clip (0,0) circle (5cm);
	\draw[shaded] (162:6cm) -- (162:5cm) .. controls (162:3cm) and (126:3cm) .. (126:5cm) -- (126:6cm);
	\draw[shaded] (0,0) -- (90:6cm)--(54:6cm)--(0,0);
	\draw[shaded] (0,0) -- (-18:6cm)--(18:6cm)--(0,0);
	\draw[shaded] (0,0) -- (-90:6cm)--(-54:6cm)--(0,0);
	\draw[shaded] (0,0) -- (-162:6cm)--(-126:6cm)--(0,0);
	\node at (0,0) [Tbox] (T) {$T$};
	\node at (T.180) [left] {$\star$};
	\node at (180:5cm) [right] {$\star$};
	\draw [ultra thick] (0,0) circle (5cm);
\end{tikzpicture}
\; , \;
w_{2}=
\begin{tikzpicture}[ATLsix]
	\clip (0,0) circle (5cm);
	\draw[shaded] (0,0)--(162:6cm) -- (126:6cm) -- (126:5cm) .. controls (126:3cm) and (90:3cm) .. (90:5cm) -- (90:6cm) -- (54:6cm) --(0,0);
	\draw[shaded] (0,0) -- (-18:6cm)--(18:6cm)--(0,0);
	\draw[shaded] (0,0) -- (-90:6cm)--(-54:6cm)--(0,0);
	\draw[shaded] (0,0) -- (-162:6cm)--(-126:6cm)--(0,0);
	\node at (0,0) [Tbox] (T) {$T$};
	\node at (T.180) [left] {$\star$};
	\node at (180:5cm) [right] {$\star$};
	\draw [ultra thick] (0,0) circle (5cm);
\end{tikzpicture}
\; , \ldots , \;
w_{10}=
\begin{tikzpicture}[ATLsix, rotate=72]
	\clip (0,0) circle (5cm);
	\draw[shaded] (0,0)--(162:6cm) -- (126:6cm) -- (126:5cm) .. controls (126:3cm) and (90:3cm) .. (90:5cm) -- (90:6cm) -- (54:6cm) --(0,0);
	\draw[shaded] (0,0) -- (-18:6cm)--(18:6cm)--(0,0);
	\draw[shaded] (0,0) -- (-90:6cm)--(-54:6cm)--(0,0);
	\draw[shaded] (0,0) -- (-162:6cm)--(-126:6cm)--(0,0);
	\node at (0,0) [Tbox] (T) {$T$};
	\node at (T.-108) [below] {$\star$};
	\node at (108:5cm) [right] {$\star$};
	\draw [ultra thick] (0,0) circle (5cm);
\end{tikzpicture}

%% file: Diagrams/ATLSixEGs/nnplusone.tex
w_{1,2}=
\begin{tikzpicture}[ATLsix]
	\clip (0,0) circle (5cm);
	
	\draw[shaded] (165:5cm) .. controls (165:3cm) and (75:3cm) .. (75:5cm) arc (75:105:5cm) .. controls (105:3cm) and (135:3cm) .. (135:5cm) arc (135:165:5cm);

	\draw[shaded] (0,0) -- (45:6cm)--(15:6cm)--(0,0);
	\draw[shaded] (0,0) -- (-45:6cm)--(-15:6cm)--(0,0);
	\draw[shaded] (0,0) -- (-105:6cm)--(-75:6cm)--(0,0);
	\draw[shaded] (0,0) -- (-165:6cm)--(-135:6cm)--(0,0);

	\node at (0,0) [Tbox] (T) {$T$};
	\node at (T.180) [left] {$\star$};
	\node at (180:5cm) [right] {$\star$};
	
	\draw [ultra thick] (0,0) circle (5cm);
\end{tikzpicture}
\; , \;
w_{2,3}=
\begin{tikzpicture}[ATLsix]
	\clip (0,0) circle (5cm);
	
	\begin{pgfonlayer}{background} 
		\fill[shaded] (0,0) circle (5cm);
	\end{pgfonlayer} 
	
	\begin{scope}[rotate=-30]
		\draw[unshaded] (165:5cm) .. controls (165:3cm) and (75:3cm) .. (75:5cm) arc (75:105:5cm) .. controls (105:3cm) and (135:3cm) .. (135:5cm) arc (135:165:5cm);

		\draw[unshaded] (0,0) -- (45:6cm)--(15:6cm)--(0,0);
		\draw[unshaded] (0,0) -- (-45:6cm)--(-15:6cm)--(0,0);
		\draw[unshaded] (0,0) -- (-105:6cm)--(-75:6cm)--(0,0);
		\draw[unshaded] (0,0) -- (-165:6cm)--(-135:6cm)--(0,0);
	\end{scope}

	\node at (0,0) [Tbox] (T) {$T$};
	\node at (T.180) [left] {$\star$};
	\node at (180:5cm) [right] {$\star$};
	
	\draw [ultra thick] (0,0) circle (5cm);
\end{tikzpicture}
\; , \ldots , \;
w_{12,1}=
\begin{tikzpicture}[ATLsix]
	\clip (0,0) circle (5cm);
	
	\begin{pgfonlayer}{background} 
		\fill[shaded] (0,0) circle (5cm);
	\end{pgfonlayer} 
	
	\begin{scope}[rotate=30]
		\draw[unshaded] (165:5cm) .. controls (165:3cm) and (75:3cm) .. (75:5cm) arc (75:105:5cm) .. controls (105:3cm) and (135:3cm) .. (135:5cm) arc (135:165:5cm);

		\draw[unshaded] (0,0) -- (45:6cm)--(15:6cm)--(0,0);
		\draw[unshaded] (0,0) -- (-45:6cm)--(-15:6cm)--(0,0);
		\draw[unshaded] (0,0) -- (-105:6cm)--(-75:6cm)--(0,0);
		\draw[unshaded] (0,0) -- (-165:6cm)--(-135:6cm)--(0,0);
	\end{scope}

	\node at (0,0) [Tbox] (T) {$T$};
	\node at (T.-70) [below right] {$\star$};
	\node at (180:5cm) [right] {$\star$};
	
	\draw [ultra thick] (0,0) circle (5cm);
\end{tikzpicture}

%% file: Diagrams/ATLSixEGs/nnplustwo.tex
w_{1,3}=
\begin{tikzpicture}[ATLsix]
	\clip (0,0) circle (5cm);
	
	\draw[shaded] (165:5cm) .. controls (165:3cm) and (135:3cm) .. (135:5cm) arc (135:165:5cm);
	\draw[shaded] (105:5cm) .. controls (105:3cm) and (75:3cm) .. (75:5cm) arc (75:105:5cm);

	\draw[shaded] (0,0) -- (45:6cm)--(15:6cm)--(0,0);
	\draw[shaded] (0,0) -- (-45:6cm)--(-15:6cm)--(0,0);
	\draw[shaded] (0,0) -- (-105:6cm)--(-75:6cm)--(0,0);
	\draw[shaded] (0,0) -- (-165:6cm)--(-135:6cm)--(0,0);

	\node at (0,0) [Tbox] (T) {$T$};
	\node at (T.180) [left] {$\star$};
	\node at (180:5cm) [right] {$\star$};
	
	\draw [ultra thick] (0,0) circle (5cm);
\end{tikzpicture}
\; , \;
w_{2,4}=
\begin{tikzpicture}[ATLsix]
	\clip (0,0) circle (5cm);
	
	\begin{pgfonlayer}{background} 
		\fill[shaded] (0,0) circle (5cm);
	\end{pgfonlayer} 
	
	\begin{scope}[rotate=-30]
		\draw[unshaded] (165:5cm) .. controls (165:3cm) and (135:3cm) .. (135:5cm) arc (135:165:5cm);
		\draw[unshaded] (105:5cm) .. controls (105:3cm) and (75:3cm) .. (75:5cm) arc (75:105:5cm);

		\draw[unshaded] (0,0) -- (45:6cm)--(15:6cm)--(0,0);
		\draw[unshaded] (0,0) -- (-45:6cm)--(-15:6cm)--(0,0);
		\draw[unshaded] (0,0) -- (-105:6cm)--(-75:6cm)--(0,0);
		\draw[unshaded] (0,0) -- (-165:6cm)--(-135:6cm)--(0,0);
	\end{scope}

	\node at (0,0) [Tbox] (T) {$T$};
	\node at (T.180) [left] {$\star$};
	\node at (180:5cm) [right] {$\star$};
	
	\draw [ultra thick] (0,0) circle (5cm);
\end{tikzpicture}
\; , \ldots , \;
w_{12,2}=
\begin{tikzpicture}[ATLsix]
	\clip (0,0) circle (5cm);
	
	\begin{pgfonlayer}{background} 
		\fill[shaded] (0,0) circle (5cm);
	\end{pgfonlayer} 
	
	\begin{scope}[rotate=30]
		\draw[unshaded] (165:5cm) .. controls (165:3cm) and (135:3cm) .. (135:5cm) arc (135:165:5cm);
		\draw[unshaded] (105:5cm) .. controls (105:3cm) and (75:3cm) .. (75:5cm) arc (75:105:5cm);

		\draw[unshaded] (0,0) -- (45:6cm)--(15:6cm)--(0,0);
		\draw[unshaded] (0,0) -- (-45:6cm)--(-15:6cm)--(0,0);
		\draw[unshaded] (0,0) -- (-105:6cm)--(-75:6cm)--(0,0);
		\draw[unshaded] (0,0) -- (-165:6cm)--(-135:6cm)--(0,0);
	\end{scope}

	\node at (0,0) [Tbox] (T) {$T$};
	\node at (T.-70) [below right] {$\star$};
	\node at (180:5cm) [right] {$\star$};
	
	\draw [ultra thick] (0,0) circle (5cm);
\end{tikzpicture}

%% file: Diagrams/ATLSixEGs/nnplussix.tex
w_{1,7}=
\begin{tikzpicture}[ATLsix]
	\clip (0,0) circle (5cm);
	
	\draw[shaded] (165:5cm) .. controls (165:3cm) and (135:3cm) .. (135:5cm) arc (135:165:5cm);
	\draw[shaded] (-15:5cm) .. controls (-15:3cm) and (-45:3cm) .. (-45:5cm) arc (-45:-15:5cm);

	\draw[shaded] (0,0) -- (45:6cm)--(15:6cm)--(0,0);
	\draw[shaded] (0,0) -- (75:6cm)--(105:6cm)--(0,0);
	\draw[shaded] (0,0) -- (-105:6cm)--(-75:6cm)--(0,0);
	\draw[shaded] (0,0) -- (-165:6cm)--(-135:6cm)--(0,0);

	\node at (0,0) [Tbox] (T) {$T$};
	\node at (T.180) [left] {$\star$};
	\node at (180:5cm) [right] {$\star$};
	
	\draw [ultra thick] (0,0) circle (5cm);
\end{tikzpicture}
\; , \;
w_{2,8}=
\begin{tikzpicture}[ATLsix]
	\clip (0,0) circle (5cm);
	
	\begin{pgfonlayer}{background} 
		\fill[shaded] (0,0) circle (5cm);
	\end{pgfonlayer} 
	
	\begin{scope}[rotate=-30]
		\draw[unshaded] (165:5cm) .. controls (165:3cm) and (135:3cm) .. (135:5cm) arc (135:165:5cm);
		\draw[unshaded] (-15:5cm) .. controls (-15:3cm) and (-45:3cm) .. (-45:5cm) arc (-45:-15:5cm);

		\draw[unshaded] (0,0) -- (45:6cm)--(15:6cm)--(0,0);
		\draw[unshaded] (0,0) -- (75:6cm)--(105:6cm)--(0,0);
		\draw[unshaded] (0,0) -- (-105:6cm)--(-75:6cm)--(0,0);
		\draw[unshaded] (0,0) -- (-165:6cm)--(-135:6cm)--(0,0);
	\end{scope}

	\node at (0,0) [Tbox] (T) {$T$};
	\node at (T.180) [left] {$\star$};
	\node at (180:5cm) [right] {$\star$};
	
	\draw [ultra thick] (0,0) circle (5cm);
\end{tikzpicture}
\;  ,  \ldots , \;
w_{6,12}=
\begin{tikzpicture}[ATLsix]
	\clip (0,0) circle (5cm);
	
	\begin{pgfonlayer}{background} 
		\fill[shaded] (0,0) circle (5cm);
	\end{pgfonlayer} 
	
	\begin{scope}[rotate=30]
		\draw[unshaded] (165:5cm) .. controls (165:3cm) and (135:3cm) .. (135:5cm) arc (135:165:5cm);
		\draw[unshaded] (-15:5cm) .. controls (-15:3cm) and (-45:3cm) .. (-45:5cm) arc (-45:-15:5cm);

		\draw[unshaded] (0,0) -- (45:6cm)--(15:6cm)--(0,0);
		\draw[unshaded] (0,0) -- (75:6cm)--(105:6cm)--(0,0);
		\draw[unshaded] (0,0) -- (-105:6cm)--(-75:6cm)--(0,0);
		\draw[unshaded] (0,0) -- (-165:6cm)--(-135:6cm)--(0,0);
	\end{scope}

	\node at (0,0) [Tbox] (T) {$T$};
	\node at (T.-110) [below left] {$\star$};
	\node at (180:5cm) [right] {$\star$};
	
	\draw [ultra thick] (0,0) circle (5cm);
\end{tikzpicture}

%% file: Diagrams/ATLFiveEGs/SampleDual.tex
\begin{tikzpicture}[ATLsix]
	\clip (0,0) circle (10cm);
	\draw[shaded] (162:12cm) -- (162:5cm) .. controls (162:3cm) and (126:3cm) .. (126:5cm) -- (126:12cm);
	\draw[shaded] (0,0) -- (90:12cm)--(54:12cm)--(0,0);
	\draw[shaded] (0,0) -- (-18:12cm)--(18:12cm)--(0,0);
	\draw[shaded] (0,0) -- (-90:12cm)--(-54:12cm)--(0,0);
	\draw[shaded] (0,0) -- (-162:12cm)--(-126:12cm)--(0,0);
	\node at (0,0) [Tbox] (T) {$T$};
	\node at (T.180) [left] {$\star$};
	\node at (180:10cm) [right] {$\star$};
	\filldraw[fill=white, thick] (152:5cm) arc (-208:136:5cm) -- (136:8cm) arc (136:-208:8cm) -- (152:5cm);
	\node at (-90:6.5cm) {$f^{(10)}$};
	\draw [ultra thick] (0,0) circle (10cm);
\end{tikzpicture}

%% file: Diagrams/ATLSixEGs/SampleDual.tex
\begin{tikzpicture}[scale=.3, baseline]
	\clip (0,0) circle (9cm);
	\begin{pgfonlayer}{background} 
		\fill[shaded] (0,0) circle (9cm);
	\end{pgfonlayer} 
	\begin{scope}[rotate=-90]
		\draw[unshaded] (165:9cm )-- (165:5cm) .. controls (165:3cm) and (75:3cm) .. (75:5cm) -- (75:9cm) arc (75:105:9cm) -- (105:5cm) .. controls (105:3cm) and (135:3cm) .. (135:5cm) -- (135:9cm) arc (135:165:9cm);

		\draw[unshaded] (0,0) -- (45:10cm)--(15:10cm)--(0,0);
		\draw[unshaded] (0,0) -- (-45:10cm)--(-15:10cm)--(0,0);
		\draw[unshaded] (0,0) -- (-105:10cm)--(-75:10cm)--(0,0);
		\draw[unshaded] (0,0) -- (-165:10cm)--(-135:10cm)--(0,0);
	\end{scope}
	\filldraw[unshaded, thick] (25:5cm) --(25:8cm) arc (25:-325:8cm) -- (-325:5cm) arc (-325:25:5cm);
	\node at (0,0) [Tbox] (T) {$T$};
	\node at (T.180) [left] {$\star$};
	\node at (-90:6.5cm) {$f^{(12)}$};
	\node at (180:8.5cm) {$\star$};
	\draw [ultra thick] (0,0) circle (9cm);
\end{tikzpicture} 
\; , \;
\begin{tikzpicture}[scale=.3, baseline]
	\clip (0,0) circle (9cm);
	\begin{pgfonlayer}{background} 
		\fill[shaded] (0,0) circle (9cm);
	\end{pgfonlayer} 
	\begin{scope}[rotate=30]
		\draw[unshaded]  (165:9cm) -- (165:5cm) .. controls (165:3cm) and (135:3cm) .. (135:5cm) --(135:9cm) arc (135:165:9cm);
		\draw[unshaded] (-15:9cm) -- (-15:5cm) .. controls (-15:3cm) and (-45:3cm) .. (-45:5cm) -- (-45:9cm)  arc (-45:-15:9cm);
		\draw[unshaded] (0,0) -- (45:10cm)--(15:10cm)--(0,0);
		\draw[unshaded] (0,0) -- (75:10cm)--(105:10cm)--(0,0);
		\draw[unshaded] (0,0) -- (-105:10cm)--(-75:10cm)--(0,0);
		\draw[unshaded] (0,0) -- (-165:10cm)--(-135:10cm)--(0,0);
	\end{scope}
	\filldraw[unshaded, thick] (5:5cm) --(5:8cm) arc (5:175:8cm) -- (175:5cm) arc (175:5:5cm);
	\filldraw[unshaded, thick] (-5:5cm) --(-5:8cm) arc (-5:-175:8cm) -- (-175:5cm) arc (-175:-5:5cm);
	\node at (0,0) [Tbox] (T) {$T$};
	\node at (T.-110) [below left] {$\star$};
	\node at (90:6.5cm) {$f^{(6)}$};
	\node at (180:8.5cm) {$\star$};
	\node at (-90:6.5cm) {$f^{(6)}$};
	\draw [ultra thick] (0,0) circle (9cm);
\end{tikzpicture} 

%% file: Diagrams/TLTen/NonZeroBeta.tex
\begin{tikzpicture}[TL12]	 
	\foreach \x in {1,...,10} \draw (\x cm , -1.5cm)--(\x cm, 1.5cm);
\end{tikzpicture}
\, , \,
\begin{tikzpicture}[TL12]	 
	\foreach \x in {1,...,8} 
		\draw (\x cm , -1.5cm)--(\x cm, 1.5cm);
	\draw (9cm ,1.5cm) arc (-180:0:.5cm);
	\draw (9cm,-1.5cm) arc (180:0:.5cm);	
\end{tikzpicture}
\, , \,
\begin{tikzpicture}[TL12]	 
	\foreach \x in {3,...,10} 
		\draw (\x cm , -1.5cm)--(\x cm, 1.5cm);
	\draw (1cm,1.5cm) arc (-180:0:.5cm);
	\draw (1cm,-1.5cm) arc (180:0:.5cm);	
\end{tikzpicture}
\, , \,
\begin{tikzpicture}[TL12]	 
	\foreach \x in {3,...,10} 
		\draw (\x cm - 2cm , -1.5cm)--(\x cm, 1.5cm);
	\draw (1cm,1.5cm) arc (-180:0:.5cm);
	\draw (9cm,-1.5cm) arc (180:0:.5cm);	
\end{tikzpicture}
\, , \,
\begin{tikzpicture}[TL12]	 
	\foreach \x in {3,...,10} 
		\draw (\x cm , -1.5cm)--(\x cm - 2cm, 1.5cm);
	\draw (9cm,1.5cm) arc (-180:0:.5cm);
	\draw (1cm,-1.5cm) arc (180:0:.5cm);	
\end{tikzpicture}

%% file: Diagrams/TLTen/One.tex
\begin{tikzpicture}[TL12]	 
	\foreach \x in {1,...,10} \draw (\x cm , -1.5cm)--(\x cm, 1.5cm);
\end{tikzpicture}

%% file: Diagrams/TLTen/Two.tex
\begin{tikzpicture}[TL12]	 
	\foreach \x in {1,...,8} 
		\draw (\x cm , -1.5cm)--(\x cm, 1.5cm);
	\draw (9cm ,1.5cm) arc (-180:0:.5cm);
	\draw (9cm,-1.5cm) arc (180:0:.5cm);	
\end{tikzpicture}

%% file: Diagrams/TLTen/Three.tex
\begin{tikzpicture}[TL12]	 
	\foreach \x in {3,...,10} 
		\draw (\x cm , -1.5cm)--(\x cm, 1.5cm);
	\draw (1cm,1.5cm) arc (-180:0:.5cm);
	\draw (1cm,-1.5cm) arc (180:0:.5cm);	
\end{tikzpicture}

%% file: Diagrams/TLTen/Four.tex
\begin{tikzpicture}[TL12]	 
	\foreach \x in {3,...,10} 
		\draw (\x cm - 2cm , -1.5cm)--(\x cm, 1.5cm);
	\draw (1cm,1.5cm) arc (-180:0:.5cm);
	\draw (9cm,-1.5cm) arc (180:0:.5cm);	
\end{tikzpicture}

%% file: Diagrams/TLTen/Five.tex
\begin{tikzpicture}[TL12]	 
	\foreach \x in {3,...,10} 
		\draw (\x cm , -1.5cm)--(\x cm - 2cm, 1.5cm);
	\draw (9cm,1.5cm) arc (-180:0:.5cm);
	\draw (1cm,-1.5cm) arc (180:0:.5cm);	
\end{tikzpicture}

%% file: Tables/TLAlphaOneAlphaOne.tex
\begin{table}[ht!]
\begin{tabular}{ l | c | c }
$\beta \in B(TL_{12})$ & $\coeff{f^{(12)}}{\beta}$ & value of tangle \\
\hline
\rule{0pt}{.6cm}
\input{Diagrams/TLTwelve/TwelveThrough/One}
& $1$
& $[2]^2 [5]$
\\
\rule{0pt}{.6cm}
\input{Diagrams/TLTwelve/TenThrough/One}
& $-\frac{[11]}{[12]}$
& $[2] [5]$
\\
\rule{0pt}{.6cm}
\input{Diagrams/TLTwelve/TenThrough/Two}
& $\frac{[10]}{[12]}$
& $[5]$
\\
\rule{0pt}{.6cm}
\input{Diagrams/TLTwelve/TenThrough/Three}
& $-\frac{[2][10]}{[12]}$
& $[2][5]$
\\
\rule{0pt}{.6cm}
\input{Diagrams/TLTwelve/TenThrough/Four}
& $\frac{[2]}{[12]}$
& $-[5]$
\\
\rule{0pt}{.6cm}
\input{Diagrams/TLTwelve/TenThrough/Five}
& $- \frac{[2]^2}{[12]}$
& $-[2] [5]$
\\
\rule{0pt}{.6cm}
\input{Diagrams/TLTwelve/EightThrough/One}
& $\frac{[9][10]}{[11][12]}$
& $[5]$
\\
\rule{0pt}{.6cm}
\input{Diagrams/TLTwelve/EightThrough/Two}
& $\frac{[2]}{[11][12]}$
& $[5]$
\\
\rule{0pt}{.6cm}
\input{Diagrams/TLTwelve/EightThrough/Three}
& $-\frac{[9] + [11]}{[11][12]}$
& $-[2] [5]$
\\
\rule{0pt}{.6cm}
\input{Diagrams/TLTwelve/EightThrough/Four}
& $\frac{[8] + [10]}{[11][12]}$
& $- [5]$
\\
\rule{0pt}{.6cm}
\input{Diagrams/TLTwelve/EightThrough/Five}
& $\frac{[2]([9] + [11])}{[11][12]}$
& $- [5]$
\\
\rule{0pt}{.6cm}
\input{Diagrams/TLTwelve/EightThrough/Six}
& $\frac{[2]+ [10][11]}{[11][12]}$
& $[2]^2 [5]$
\\
\rule{0pt}{.6cm}
\input{Diagrams/TLTwelve/EightThrough/Seven}
& $-\frac{[2]^2+ [9][11]}{[11][12]}$
& $[2] [5]$
\\
\rule{0pt}{.6cm}
\input{Diagrams/TLTwelve/EightThrough/Eight}
& $\frac{[2][8]+[2]^2[11]+ [9][10]^2}{[10][11][12]}$
& $[5]$
\\
\rule[-.4cm]{0pt}{1cm}
\input{Diagrams/TLTwelve/EightThrough/Nine}
& $\frac{[2]^2[9]+[2]^2[11]+ [2][9][10][11]}{[10][11][12]}$
& $[5]$
\\
\end{tabular}

\vspace{6pt}
\caption{The elements of $TL_{12}$ which contribute to $\ip{\alpha_{i,i+1}(T), w_{i,i+1}(T)}$. }\label{TLTwelve}
\end{table}

%% file: Diagrams/TLsix/One.tex
\begin{tikzpicture}[TL12]	 
	\foreach \x in {1,...,6} \draw (\x cm , -1.5cm)--(\x cm, 1.5cm);
\end{tikzpicture}

%% file: Diagrams/TLsix/Two.tex
\begin{tikzpicture}[TL12]	 
	\foreach \x in {3,...,6} \draw (\x cm , -1.5cm)--(\x cm, 1.5cm);
	\draw (1cm,1.5cm) arc (-180:0:.5cm);
	\draw (1cm,-1.5cm) arc (180:0:.5cm);
\end{tikzpicture}

%% file: Diagrams/TLsix/Three.tex
\begin{tikzpicture}[TL12]	 
	\foreach \x in {4,...,6} \draw (\x cm , -1.5cm)--(\x cm, 1.5cm);
	\draw (2cm,1.5cm) arc (-180:0:.5cm);
	\draw (1cm,1.5cm) -- (3cm,-1.5cm);
	\draw (1cm,-1.5cm) arc (180:0:.5cm);
\end{tikzpicture}

%% file: Diagrams/TLsix/Four.tex
\begin{tikzpicture}[TL12]	 
	\foreach \x in {1,...,4} \draw (\x cm , -1.5cm)--(\x cm + 2cm, 1.5cm);
	\draw (1cm,1.5cm) arc (-180:0:.5cm);
	\draw (5cm,-1.5cm) arc (180:0:.5cm);
\end{tikzpicture}

%% file: Diagrams/TLsix/Five.tex
\begin{tikzpicture}[TL12]	 
	\foreach \x in {2,...,4} \draw (\x cm , -1.5cm)--(\x cm + 2cm, 1.5cm);
	\draw (1cm,1.5cm)--(1cm,-1.5cm);
	\draw (2cm,1.5cm) arc (-180:0:.5cm);
	\draw (5cm,-1.5cm) arc (180:0:.5cm);
\end{tikzpicture}

%% file: Diagrams/TLsix/Six.tex
\begin{tikzpicture}[TL12]	 
	\foreach \x in {3,...,4} \draw (\x cm , -1.5cm)--(\x cm + 2cm, 1.5cm);
	\draw (2cm,1.5cm) arc (-180:0:.5cm);
	\draw (1cm,1.5cm) arc (-180:0:1.5cm);
	\draw (5cm,-1.5cm) arc (180:0:.5cm);
	\draw (1cm,-1.5cm) arc (180:0:.5cm);
\end{tikzpicture}

%% file: Tables/TLAlphaOneAlphaSix.tex
\begin{table}[ht!]
\begin{tabular}{ c | c | c | c | c | c | c }
&
\rule[-.4cm]{0pt}{.8cm} \input{Diagrams/TLsix/One}
&
\input{Diagrams/TLsix/TwoReflect}
&
\input{Diagrams/TLsix/ThreeReflect}
&
\input{Diagrams/TLsix/FourReflect}
&
\input{Diagrams/TLsix/FiveReflect}
&
\input{Diagrams/TLsix/SixReflect}
\\
\hline
\rule{0pt}{.6cm} \input{Diagrams/TLsix/One} & 0 & 0 & 0 & $\frac{[5]}{[6]}$ & $\frac{-[2]^2[5]}{[6]}$ & $\frac{([3]+[5])[5]}{[5][6]}$
\\
\rule{0pt}{.6cm} \input{Diagrams/TLsix/Two} & 0 & 0 & 0 & $\frac{-[2][5]^2}{[6]^2}$ &  $\frac{[2][5]^2}{[6]^2}$ & 0 
\\
\rule{0pt}{.6cm} \input{Diagrams/TLsix/Three} & 0 & 0 & 0 & $\frac{[4][5]}{[6]^2}$ & 0 & 0 
\\
\rule{0pt}{.6cm} \input{Diagrams/TLsix/Four} & $\frac{-[5]}{[6]}$ & $\frac{[2][5]^2}{[6]^2}$ & $\frac{-[4][5]}{[6]^2}$ & 0 & 0 & 0 
\\
\rule{0pt}{.6cm} \input{Diagrams/TLsix/Five} & $\frac{[2]^2[5]}{[6]}$ & $\frac{-[2][5]^2}{[6]^2}$ & 0 & 0 & 0 & 0 
\\
\rule[-.4cm]{0pt}{1cm} \input{Diagrams/TLsix/Six} & $\frac{-([3]+[5])[5]}{[5][6]}$ & 0 & 0 & 0 & 0 & 0 
\\
\end{tabular}
\vspace{6pt}
\caption{The elements of $TL_{6}$ which contribute to $\ip{\alpha_{i,i+6}(T), w_{i,i+1}(T)}$, and their contributions. }\label{TLSixFirst}
\end{table}

%% file: Diagrams/TLsix/FourReflect.tex
\begin{tikzpicture}[TL12,xscale=-1]	 
	\foreach \x in {1,...,4} \draw (\x cm , -1.5cm)--(\x cm + 2cm, 1.5cm);
	\draw (1cm,1.5cm) arc (-180:0:.5cm);
	\draw (5cm,-1.5cm) arc (180:0:.5cm);
\end{tikzpicture}

%% file: Diagrams/TLsix/TwoReflect.tex
\begin{tikzpicture}[TL12,xscale=-1]	 
	\foreach \x in {3,...,6} \draw (\x cm , -1.5cm)--(\x cm, 1.5cm);
	\draw (1cm,1.5cm) arc (-180:0:.5cm);
	\draw (1cm,-1.5cm) arc (180:0:.5cm);
\end{tikzpicture}

%% file: Diagrams/TLsix/Seven.tex
\begin{tikzpicture}[TL12]	 
	\foreach \x in {3,4} \draw (\x cm , -1.5cm)--(\x cm, 1.5cm);
	\draw (1cm,1.5cm) arc (-180:0:.5cm);
	\draw (1cm,-1.5cm) arc (180:0:.5cm);
	\draw (5cm,1.5cm) arc (-180:0:.5cm);
	\draw (5cm,-1.5cm) arc (180:0:.5cm);
\end{tikzpicture}

%% file: Tables/TLAlphaSixAlphaSix.tex
\begin{table}[ht!]
\begin{tabular}{ c | c | c | c | c | c | c }
&
\rule[-.4cm]{0pt}{.8cm} \input{Diagrams/TLsix/One}
&
\input{Diagrams/TLsix/TwoReflect}
&
\input{Diagrams/TLsix/Four}
&
\input{Diagrams/TLsix/FourReflect}
&
\input{Diagrams/TLsix/Two}
&
\input{Diagrams/TLsix/Seven}
\\
\hline
\rule{0pt}{.6cm} \input{Diagrams/TLsix/One} & $[2]^2[5]$ & $\frac{-[2][5]^2}{[6]}$ & 0 & 0 & $\frac{-[2][5]^2}{[6]}$ & $\frac{([2]+[4][5])[5]}{[5][6]}$
\\
\rule{0pt}{.6cm} \input{Diagrams/TLsix/Two} & $\frac{-[2][5]^2}{[6]}$ & 0 & 0 & 0 &  $\frac{[5]^3}{[6]^2}$ & 0 
\\
\rule{0pt}{.6cm} \input{Diagrams/TLsix/FourReflect} & 0 & 0 & 0 & $\frac{-[5]}{[6]^2}$ & 0 & 0 
\\
\rule{0pt}{.6cm} \input{Diagrams/TLsix/Four} & 0 & 0 & $\frac{-[5]}{[6]^2}$ & 0 & 0 & 0 
\\
\rule{0pt}{.6cm} \input{Diagrams/TLsix/TwoReflect} & $\frac{-[2][5]^2}{[6]}$ & $\frac{[5]^3}{[6]^2}$ & 0 & 0 & 0 & 0 
\\
\rule[-.4cm]{0pt}{1cm} \input{Diagrams/TLsix/Seven} & $\frac{([2]+[4][5])[5]}{[5][6]}$ & 0 & 0 & 0 & 0 & 0 
\\
\end{tabular}

\vspace{6pt}
\caption{The elements of $TL_{6}$ which contribute to $\ip{\alpha_{i,i+6}(T), w_{i,i+6}(T)}$, and their contributions. }\label{TLSixSecond}
\end{table}

%% file: Diagrams/TLTwelve/PicsGivingFiveSix/One.tex
\begin{tikzpicture}[TL12]	 
	\foreach \x in {1,2,3,4} 
		\draw (\x cm , -3cm)--(\x cm + 4cm, 3cm);
	\foreach \x in {9,10,11,12} 
		\draw (\x cm , -3cm)--(\x cm, 3cm);
	\draw (1cm,3cm) arc (-180:0:1.5cm);
	\draw (2cm,3cm) arc (-180:0:.5cm);
	\draw (5cm,-3cm) arc (180:0:1.5cm);
	\draw (6cm,-3cm) arc (180:0:.5cm);
	
\end{tikzpicture}

%% file: Diagrams/TLTwelve/PicsGivingFiveSix/OneReflect.tex
\begin{tikzpicture}[TL12, xscale=-1]	 
	\foreach \x in {1,2,3,4} 
		\draw (\x cm , -3cm)--(\x cm + 4cm, 3cm);
	\foreach \x in {9,10,11,12} 
		\draw (\x cm , -3cm)--(\x cm, 3cm);
	\draw (1cm,3cm) arc (-180:0:1.5cm);
	\draw (2cm,3cm) arc (-180:0:.5cm);
	\draw (5cm,-3cm) arc (180:0:1.5cm);
	\draw (6cm,-3cm) arc (180:0:.5cm);
	
\end{tikzpicture}

%% file: Diagrams/TLTwelve/PicsGivingFiveSix/Two.tex
\begin{tikzpicture}[TL12]	 
	\foreach \x in {1,2,3,4} 
		\draw (\x cm , -3cm)--(\x cm + 2cm, 3cm);
	\foreach \x in {9,10,11,12} 
		\draw (\x cm , -3cm)--(\x cm-2cm, 3cm);
	\draw (1cm,3cm) arc (-180:0:.5cm);
	\draw (11cm,3cm) arc (-180:0:.5cm);
	\draw (5cm,-3cm) arc (180:0:1.5cm);
	\draw (6cm,-3cm) arc (180:0:.5cm);
\end{tikzpicture}

%% file: Diagrams/TLTwelve/PicsGivingFiveSix/Three.tex
\begin{tikzpicture}[TL12]	 
	\foreach \x in {2,3,4} 
		\draw (\x cm , -3cm)--(\x cm + 2cm, 3cm);
	\foreach \x in {9,10,11,12} 
		\draw (\x cm , -3cm)--(\x cm-2cm, 3cm);
	\draw (2cm,3cm) arc (-180:0:.5cm);
	\draw (11cm,3cm) arc (-180:0:.5cm);
	\draw (5cm,-3cm) arc (180:0:1.5cm);
	\draw (6cm,-3cm) arc (180:0:.5cm);
	
	\draw (1cm,-3cm)--(1cm,3cm);
\end{tikzpicture}

%% file: Diagrams/TLTwelve/PicsGivingFiveSix/ThreeReflect.tex
\begin{tikzpicture}[TL12, xscale=-1]	 
	\foreach \x in {2,3,4} 
		\draw (\x cm , -3cm)--(\x cm + 2cm, 3cm);
	\foreach \x in {9,10,11,12} 
		\draw (\x cm , -3cm)--(\x cm-2cm, 3cm);
	\draw (2cm,3cm) arc (-180:0:.5cm);
	\draw (11cm,3cm) arc (-180:0:.5cm);
	\draw (5cm,-3cm) arc (180:0:1.5cm);
	\draw (6cm,-3cm) arc (180:0:.5cm);
	
	\draw (1cm,-3cm)--(1cm,3cm);
\end{tikzpicture}

%% file: text/levels56.tex
Now we know enough about $\hat{w}_{i,j}(T)$ to prove (R5) and (R6).

\begin{proof}[Proof of (R5)]
To show 
$$\join{3}{T}{T}=  \frac{[5]}{[6]}\f{5} - i\sqrt{\frac{8(3+\sqrt{13})}{3}}\hat{w}_{10}$$
we show
$$ \Norm{ \join{3}{T}{T} }^2 = \frac{[5]^2}{[4]}$$
and
$$ \Norm{ \proj{TL_5}{ \join{3}{T}{T} } }^2= \Norm{ \frac{[5]}{[6]}\f{5} }^2= \frac{[5]^2}{[6]},$$
$$ \Norm{ \proj{ATL_5(T)}{  \join{3}{T}{T} } }^2 = \Norm{ -i\sqrt{\frac{8(3+\sqrt{13})}{3}} \hat{w}_{10} }^2 =
 \frac{8(3+\sqrt{13})}{3\sqrt{2(5+\sqrt{13})}},$$
 and then one can check that 
$$\frac{[5]^2}{[4]}=\frac{8(3+\sqrt{13})}{3\sqrt{2(5+\sqrt{13})}}+\frac{[5]^2}{[6]}.$$

First,
\begin{align*}
 \Norm{ \join{3}{T}{T} }^2 & = \input{Diagrams/InnerProducts/JoinThreeTTJoinThreeTT} \\
 & = \frac{[5]^2}{[4]^2} [4] = \frac{[5]^2}{[4]} .
 \end{align*}

Next,
\begin{align*}
 \Norm{ \proj{TL_5}{ \join{3}{T}{T} } }^2 & = \Norm{ \sum_{\beta \in B(TL_5)} \ip{ \join{3}{T}{T}, \beta} \hat{\beta} }^2 \\
 & = \Norm{ \ip{ \join{3}{T}{T}, \input{Diagrams/TLFive/One}  } \widehat{ \input{Diagrams/TLFive/One} } \; }^2
\\ 
  & =\Norm{ [5] \cdot \frac{\f{5}}{[6]} }^2= \frac{[5]^2}{[6]}.
\end{align*}

Finally, 
\begin{align*}
 \Norm{ \proj{ATL_5(T)}{ \join{3}{T}{T} } }^2 & = \Norm{ \sum_{\beta \in B(ATL_5(T))} \ip{ \join{3}{T}{T}, \beta} \hat{\beta} }^2 \\
& = \Norm{\ip{ \join{3}{T}{T}, w_5} \hat{w}_5 + \ip{ \join{3}{T}{T}, w_{10}} \hat{w}_{10} }^2 \\
& = \Norm{ 0 \hat{w}_5  - i\sqrt{\frac{8(3+\sqrt{13})}{3}} \hat{w}_{10} }^2
=
 \frac{8(3+\sqrt{13})}{3\sqrt{2(5+\sqrt{13})}}
\end{align*}

\end{proof}

To prove (R6), we need to know more about $N$ first.

\begin{Lemma}
$$
N=\join{2}{T}{T} - [5] \cdot \frac{f^{(6)}}{[7]} + i\sqrt{\frac{8(3+\sqrt{13})}{3}} \hat{w}_{6,12}
$$
and 
$$
\| N \|^2 = \frac{16}{351} (13+11 \sqrt{13}).
$$
\end{Lemma}

\begin{proof}
Recall the definition of $N$:
$$N=(1-\text{Proj}_{TL}-\text{Proj}_{ATL_6(T)}) \join{2}{T}{T}.$$
$B(TL_6)$ and $B(ATL_6(T))$ are the usual bases of $TL_6$ and $ATL_6(T)$.  We calculate
$$
\proj{TL_6}{\join{2}{T}{T}}  =
\sum_{\beta \in B(TL_6)} \ip{\join{2}{T}{T},\beta} \cdot \hat{\beta}
=[5] \cdot \frac{f^{(6)}}{[7]},
$$
since \,
\begin{tikzpicture}[scale=.2]
	\draw (0,0) --(0,2);
	\draw (1,0)--(1,2);
	\draw (2,0)--(2,2);
	\draw (3,0)--(3,2);
	\draw (4,0)--(4,2);
	\draw (5,0)--(5,2);
\end{tikzpicture} \,
is the only element of $B(TL_6)$ with which $\join{2}{T}{T}$ has non-zero inner product, and $\frac{f^{(6)}}{[7]}$ is dual to this element in $B(TL_6)$.  Also
$$
\proj{ATL_6(T)}{\join{2}{T}{T}}  =
\sum_{\gamma \in B(ATL_6(T))} \ip{\join{2}{T}{T},\gamma} \cdot \hat{\gamma}
= - i\sqrt{\frac{8(3+\sqrt{13})}{3}} \cdot\hat{w}_{6,12},
$$
since every element of $B(ATL_6(T))$ caps off one of the $T$s, except $\hat{w}_{5,6}$, $\hat{w}_{11,12}$ and $\hat{w}_{6,12}$.  When we compute the inner products of these with $\join{2}{T}{T}$, we get
$$\input{Diagrams/InnerProducts/JoinTwoTTALT}.$$
Using sphericality and multiplying by $(-1)$ when we need to rotate $T$'s star, we see that these are respectively $Z(T^3)=0$, $(-1) Z(T^3)=0$ and $(-1) Z( (\rho^{1/2} T)^3)= - i\sqrt{\frac{8(3+\sqrt{13})}{3}}$.

Next
\begin{align*}
\Norm{N}^2 & = \Norm{\join{2}{T}{T}}^2 - \frac{[5]^2}{[7]^2} \Norm{f^{(6)}}^2 - \frac{8(3+\sqrt{13})}{3} \Norm{\hat{w}_{6,12}}^2 \\
\intertext{and by  Lemma \ref{dualbasis} part (\ref{normiiplussix}) we get}
& = \frac{[5]^2}{[3]} - \frac{[5]^2}{[7]^2} [7] - \frac{8(3+\sqrt{13})}{3} \frac{3}{11-\sqrt{13}} =\frac{16}{351} (13+11 \sqrt{13})
\end{align*}
\end{proof}

\begin{proof}[Proof of (R6)]
To prove 
$$\joint{2}{T}{T} = - \frac{[5]}{[7]} \rho^{-1/2}\f{6} +\frac{ \sqrt{ 4-\sqrt{13} }}{2} N  -   i\sqrt{\frac{8(3+\sqrt{13})}{3}} \hat{w}_{4,5} +  i\sqrt{\frac{8(3+\sqrt{13})}{3}}  \hat{w}_{10,11}$$
we show 
\begin{align*}
\| \joint{2}{T}{T} \|^2 & = 2(3+\sqrt{13})\\
\intertext{ and }
\| \proj{N}{\joint{2}{T}{T}} \|^2 & =  \Norm{\frac{ \sqrt{ 4-\sqrt{13} }}{2}N }^2  = \frac{4(-91+31 \sqrt{13})}{351}, \\
 \Norm{ \proj{TL_6}{\joint{2}{T}{T}} }^2 & = \Norm{ [5] \cdot \frac{ \rho^{-1/2} \f{6} }{ [7] } }^2 = \frac{26+6 \sqrt{13}}{13} ,\\
\Norm{ \proj{ATL_6(T)}{\joint{2}{T}{T}} }^2 & = 
\Norm{ -i \sqrt{\frac{8(3+\sqrt{13})}{3}} \hat{w}_{4,5}  + i \sqrt{\frac{8(3+\sqrt{13})}{3}} \hat{w}_{10,11} }^2 \\
& = \frac{8(17+4\sqrt{13})}{27} \\
\end{align*}
and then one can check that
$$
2(3+\sqrt{13}) =  \frac{4(-91+31 \sqrt{13})}{351} + \frac{26+6 \sqrt{13}}{13} +  \frac{8(17+4\sqrt{13})}{27}.
$$
  
First,
 \begin{align*}
\| \joint{2}{T}{T} \|^2 & = 
 \input{Diagrams/InnerProducts/JointTwoTTJointTwoTT} \\
 & = \frac{[5]^2}{[3]} = 2(3+\sqrt{13})
  \end{align*}
  
Second,
\begin{align*}
\| \proj{N}{\joint{2}{T}{T}} \|^2 & = \Norm{ \frac{\ip{\joint{2}{T}{T},N}}{\ip{N,N}}N }^2  \\ 
= & 
\| \frac{1}{\ip{N,N}} 
(
	\ip{\joint{2}{T}{T},\join{2}{T}{T}}
	-\frac{[5]}{[7]} \ip{\joint{2}{T}{T},f^{(6)}} 
\\
&	 - i \sqrt{\frac{8(3+\sqrt{13})}{3}} \ip{\joint{2}{T}{T},\hat{w}_{6,12}} ) N \|^2
\\
\intertext{and by evaluating pictures and Lemma \ref{dualbasis} part (\ref{decompiiplussix}), which tells us the coefficients of $w_{4,5}$ and $w_{10,11}$ (the only two elements of $B(ATL_6(T))$ with non-zero inner product with $\joint{2}{T}{T}$) in $\hat{w}_{6,12}$, we get }
 = &
 \|  \frac{351}{16 (13+11\sqrt{13})}  ( \frac{[5]^2}{[3] [4]} +\frac{[5]}{[7]} \cdot \frac{[5]}{[6]} 
\\
& -i \sqrt{\frac{8(3+\sqrt{13})}{3}} ( i\sqrt{\frac{8(3+\sqrt{13})}{3}} \cdot \frac{-6}{[6] (11-\sqrt{13})})) N \|^2
\\
 = & \|  \frac{ \sqrt{ 4-\sqrt{13} }}{2}N \|^2 =  \frac{4(-91+31 \sqrt{13})}{351}
\end{align*}

Next,
\begin{align*}
\| \proj{TL_6}{\joint{2}{T}{T}} \|^2 & = \| \sum_{\beta \in B(TL_6)} \ip{\joint{2}{T}{T},\beta} \cdot \hat{\beta} \|^2 \\
& = \Norm{\ip{\joint{2}{T}{T}, \input{Diagrams/RotateTLId} } \cdot \widehat{\input{Diagrams/RotateTLId}} }^2 \\
& =\Norm{ [5] \cdot \frac{ \rho^{-1/2} \f{6} }{ [7] } }^2 \\
& = \frac{[5]^2}{[7]} = \frac{26+6 \sqrt{13}}{13}
\end{align*}

Finally, 
\begin{align*}
\| \proj{ATL_6(T)}{\joint{2}{T}{T}} \|^2 & = \| \sum_{\gamma \in B(ATL_6(T)) } \ip{\joint{2}{T}{T},\gamma} \cdot \hat{\gamma} \|^2 \\
& = \| \ip{\joint{2}{T}{T},w_{4,5}} \hat{w}_{4,5}  + \ip{\joint{2}{T}{T},w_{10,11}} \hat{w}_{10,11} \|^2 \\
& = \| -i \sqrt{\frac{8(3+\sqrt{13})}{3}} \hat{w}_{4,5}  + i \sqrt{\frac{8(3+\sqrt{13})}{3}} \hat{w}_{10,11} \|^2 \\
& = \frac{8(3+\sqrt{13})}{3} (\| \hat{w}_{4,5} \|^2 + \| \hat{w}_{10,11}  \|^2 -\ip{\hat{w}_{4,5} ,\hat{w}_{10,11,} } - \ip{\hat{w}_{10,11} , \hat{w}_{4,5} }) \\
\intertext{ and by Lemma \ref{dualbasis} parts (\ref{normiiplusone}) and (\ref{decompiiplusone}), we get }
& = \frac{8(17+4\sqrt{13})}{27}
\end{align*}
\end{proof}

%% file: Diagrams/TLFive/One.tex
\begin{tikzpicture}[TL12,baseline=-.1cm]	 
	\foreach \x in {1,...,6} \draw (\x cm , -1.5cm)--(\x cm, 1.5cm);
\end{tikzpicture}

%% file: Diagrams/RotateTLId.tex
\begin{tikzpicture}[scale=.2]
	\draw (0,0) arc (180:0:.5cm);
	\draw (2,0)--(0,2);
	\draw (3,0)--(1,2);
	\draw (4,0)--(2,2);
	\draw (5,0)--(3,2);
	\draw (4,2) arc (-180:0:.5cm);
\end{tikzpicture}

%% file: text/PATisHaagerup.tex
Recall that we are trying to show that the planar algebra generated by $T$, abbreviated $PA(T)$, is a subfactor planar algebra.  All the required properties are inherited from $PABG(H)$, except one.  We need to check that the zero-box space of $PA(T)$ is one-dimensional.

In Section \ref{Quadratic} we proved four quadratic skein relations; in this section, we show that these relations (together with the annular relations that defined $T$) are sufficient to evaluate any $0$-box.  Specifically, we describe how to take any diagram of $T$s joined among themselves with no output strands, and turn this diagram into an element of Temperley-Lieb.   This procedure is inductive; to show it can be used to completely evaluate any closed diagram, it is enough to show that we can always reduce the number of $T$s in a diagram.  Relations (R4), (R4'), and (R5) accomplish this immediately, but relation (R6) does not; it preserves the number of $T$s, but changes the way they're connected.  Thus, the main thing we need to show in this section is that we can always use relation (R6) to change the connectivity of a diagram so that one of relations (R4), (R4'), or (R5) can be applied.

\begin{figure}[!ht]
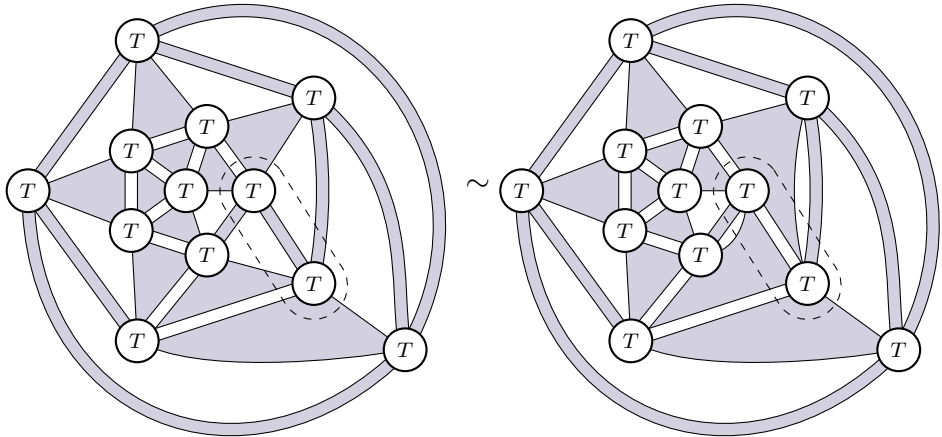

$\input{Diagrams/IcosahedronNetwork/Original} \sim \input{Diagrams/IcosahedronNetwork/Fixed}$
\caption{A diagram which cannot be reduced, but can be reduced after an application of (R6) in the indicated region.}\label{icosahedron}
\end{figure}

Figure \ref{icosahedron} is an example of a diagram where (R4), (R4'), and (R5) cannot be applied immediately, but (R5) can be applied after (R6) has been applied in a certain place.  In what follows, we use the $\sim$ symbol to mean that two diagrams are similar, i.e. they are non-zero multiples of each other, up to some lower-order terms (ordered by the number of inputs).  For example, (R6) says that $\join{2}{T}{T} \sim \joint{2}{T}{T}$.   When working with similarity instead of equality, we can forget the placement of stars in unshaded regions around each $T$, because changing the position of a star only changes the value of the diagram by a factor of $-1$.

\begin{figure}[!ht]
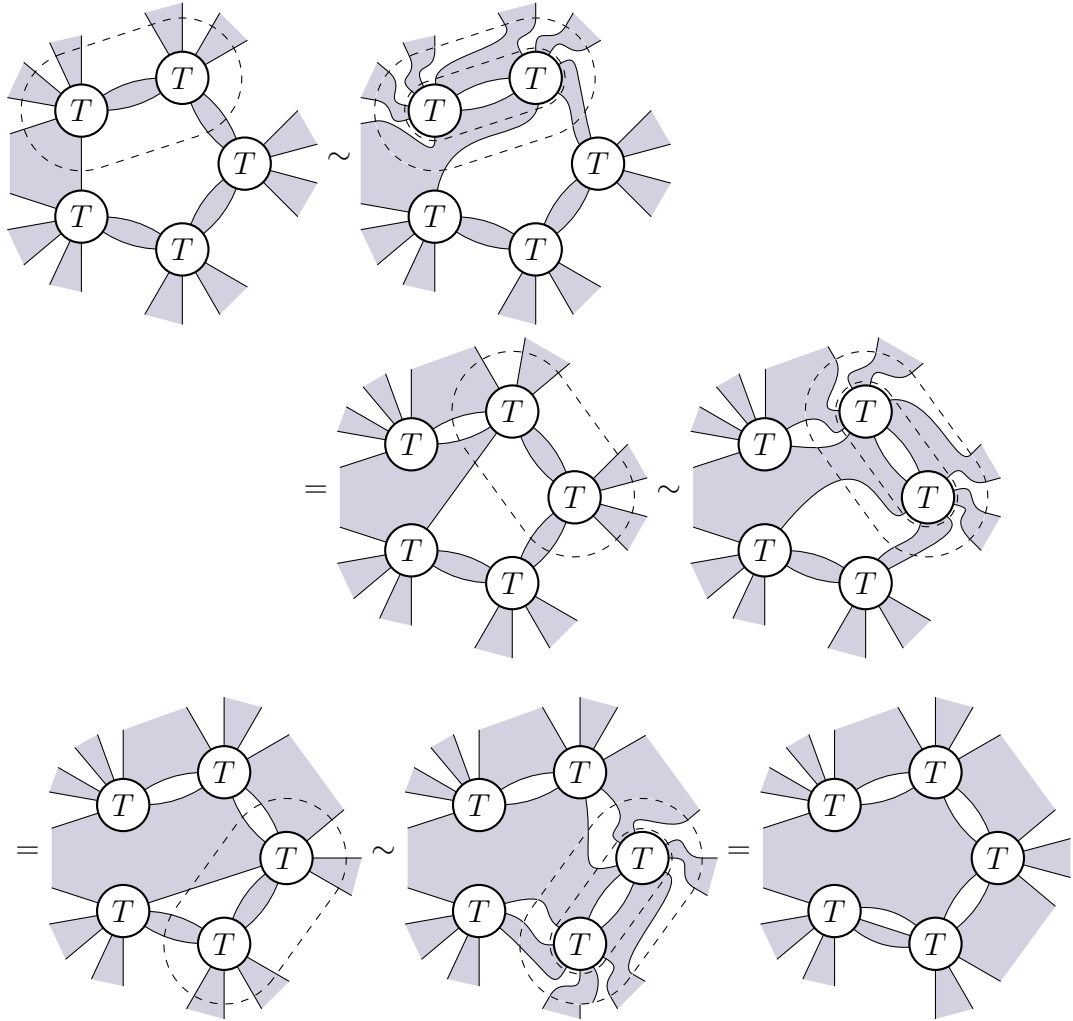

\begin{multline*}
\input{Diagrams/MakeTripleEdge/Pentagon}
\sim 
\input{Diagrams/MakeTripleEdge/SquiglySquare} 
 \\ 
=
 \input{Diagrams/MakeTripleEdge/Square} 
 \sim
\input{Diagrams/MakeTripleEdge/SquiglyTriangle} 
\end{multline*}
\begin{multline*}
 =
\input{Diagrams/MakeTripleEdge/Triangle} 
 \sim
\input{Diagrams/MakeTripleEdge/SquiglyTripleEdge}
 =
\input{Diagrams/MakeTripleEdge/TripleEdge}
\end{multline*}
\caption{Relation (R6) is used to pass an edge across a pentagon so that relation (R5) can be applied.}\label{reducepentagon}
\end{figure}

We reduced this diagram by passing an edge across a triangular face to create a triple edge.  This trick is generally applicable if we have an $n$-gon in which $(n-1)$ of the edges are double edges.  (By {\em double edge}, we mean an adjacent pair of edges between the same two vertices; similarly for {\em triple edge}.)  See Figure \ref{reducepentagon} for a more general example.

Using this trick, we can always reduce the number of $T$s in a diagram.  To show this, we must show that any diagram of $T$s contains an $n$-gon in which $(n-1)$ of the edges are doubled.  This is a fortunate consequence of working in a graph where all vertices have degree 8.  The argument proving this is somewhat involved, but uses nothing more sophisticated than the Euler characteristic.

\begin{Lemma}[Reduction Lemma]\label{ngon} If $G$ is a planar (multi-)graph with vertices of degree 8 and no edges connecting a vertex to itself, then $G$ contains either a triple edge, or an $n$-gonal face which has $(n-1)$ doubled edges.
\end{Lemma}

\begin{proof}
Suppose that $G$ has no triple edges.  To simplify the following argument, we switch from thinking about single and double edges to thin and thick edges.

Consider a graph with only triangular faces; if we could show that more than a third of the edges were thick, we would be done, because some face would (by the pigeonhole principle) have two thick edges.  If our graph had all square faces, we could proceed by showing that more than half of all edges were thick; and the percentage of edges which must be thick goes up as the number of sides per face goes up.
Fortunately, we are working in a graph with a fixed number of vertices and $d_{\text{thin}}(v)+2 d_{\text{thick}}(v)=8$ for all $v$ (where $d_{-}$ is the degree of a type of edge).  The Euler characteristic tells us that in this kinds of graph, there is a trade-off between how many sides the average face has and how many edges have to be thick.  For instance, if all the faces are triangles then three-fifths of the edges are thick; if all the faces are squares, then all edges are thick (of course, this is impossible since it gives the square tiling of the plane and our graph is finite.)
So as the number of non-triangle faces in our graph goes up, so does the proportion of edges which have to be thick.  

To make the above paragraph more precise, we first note some relations among edges and vertices.  Let $f_k$ be the number of $k$-gonal faces in $G$,
let $e$ be the number of edges, and let $v$ be the number of vertices.  Say $d$ is the average degree of the graph; $d$ must be between $4$ and $8$.  Then we have the relations $e = \frac{3}{2} f_3 + 2 f_4 + \frac{5}{2} f_5 + \cdots$,
and $d * v = 2 e$.   
The Euler characteristic tells us 
\begin{align*}
2= & v-e+f \\
  = & \frac{2e}{d} - e + (f_3 + f_4 + f_5 + \cdots ) \\
  = &  (\frac{2}{d}-1)e+(  f_3 + f_4 + f_5 + \cdots )\\
  = & (\frac{2}{d}-1)(\frac{3}{2} f_3 + 2 f_4 + \frac{5}{2} f_5 + \cdots)+( f_2 + f_3 + 
         f_4 + f_5 + \cdots ),\\
  \end{align*}
  which implies
  \begin{align*}
 2 = & ( \frac{3}{d}-\frac{1}{2})f_3 + (\frac{4}{d}-1)f_4 + (\frac{5}{d}-\frac{3}{2})f_5 + \cdots .
  \end{align*}
  
Knowing that the right hand side of this equation is non-negative is very informative.  We immediately see that $d < 6$, as $d \geq 6$ would make all coefficients of $f_i$'s zero or negative.   

In fact, since $d \geq 4$, we know that the coefficients of all the $f_i$ except for $f_3$ are zero or negative; this means that the squares, pentagons, and general $n$-gons in this graph all have to be balanced out by a large number of triangles.  More precisely,
\begin{align*}
f_3 & = \frac{ 2 -  (\frac{4}{d}-1)f_4 - (\frac{5}{d}-\frac{3}{2})f_5 - \cdots }{\frac{3}{d}-\frac{1}{2}} \\
& = \frac{ 4d - (8-2d) f_4 - (10 - 3 d) f_5 - \cdots }{6-d}
\end{align*}

Now suppose that our graph has no $n$-gon with $(n-1)$-thick edges, which is to say that at most $(n-2)$ edges of each $n$-gon are thick.   Then the total number of thick edges in the graph is
\begin{align*}
2 \cdot \text{thick edges} &  \leq f_3+2 f_4 + 3 f_5 + \cdots \\
& = \frac{ 4d - (8-2d) f_4 - (10 - 3 d) f_5 - \cdots }{6-d} +2 f_4 + 3 f_5 + \cdots \\
& = \frac{ 4d + 4 f_4 +8 f_5 + \cdots }{6-d} \\
\intertext{On the other hand, since each vertex in our original graph had degree 8,  we know that the fraction of edges in $G$ which  are thick is $\frac{8-d}{d}$.  So we also have }
2 \cdot  \text{thick edges} & = \frac{8-d}{d}(3 f_3+4 f_4 + 5 f_5 + \cdots )\\
& =  \frac{8-d}{d} (3 \frac{4d -(8-2d) f_4 - (10 - 3 d) f_5 - \cdots }{6-d} +4 f_4 + 5 f_5 + \cdots )\\
& = \frac{8-d}{6 - d}(12 + 2 f_4 + 4 f_5 + \cdots ) \\
& > \frac{2}{6-d}( 12 + 2 f_4 + 4 f_5 + \cdots )=\frac{24 + 4 f_4 + 8 f_5 + \cdots }{6-d}\\
\end{align*}
since $d <6$ implies $8-d > 2$.

But
$$\frac{24 + 4 f_4 + 8 f_5 + \cdots }{6-d} < 2 \cdot \text{thick edges} \leq \frac{4d + 4 f_4 + 8 f_5 + \cdots }{6-d}
$$ 
implies $6<d$, which is a contradiction.
Therefore the graph must contain an $n$-gon with $(n-1)$ thick edges.
\end{proof}

 \begin{Theorem} $PA(T)$ is a subfactor planar algebra. 
 \end{Theorem}\label{IsSubfactor}

\begin{proof}
Because $PA(T)$ is a sub-planar algebra of $PABG(H)$, all the properties needed to be a subfactor planar algebra are inherited, except for $\dim(PA(T)_{0,\pm})=1$.  Thus we must show that any closed diagram evaluates to a multiple of the single Temperley-Lieb zero-box.  We do this by induction on the number of $T$s in a diagram.  The base cases are no $T$s or one $T$.  By definition, a diagram with no $T$s in it is in Temperley-Lieb.  In a diagram with one $T$, which is planarly joined to itself, there is an innermost cap somewhere, which kills $T$; so this evaluates to zero.  

For the inductive step, consider any diagram $\tau$.  If $\tau$ has two $T$s connected by at least three adjacent strings, rotate the stars on these $T$s to make them look like the left hand side of (R5).  Then apply (R5).  The result is a sum of diagrams, each of which has at least one less $T$ than $\tau$ has.  

If $\tau$ does not have two $T$s connected by at least three strings, then by Lemma \ref{ngon} we know that somewhere in $\tau$ we can find an $n$-gon face with $(n-1)$ double edges.  We repeatedly apply (R6) to the double edges of this face.  Each time we apply the relation, we cut a double edge out of the face; finally we get down to a triangle which has two double edges.  When we apply (R6) again we create a new string between two $T$s already joined by two strings; in other words, a triple edge.  Then we apply (R5) to reduce the number of $T$s.
\end{proof}
  
 \begin{Theorem} $PA(T)$ is irreducible, i.e., $\dim(PA(T)_1)=1$.
 \end{Theorem}
 \begin{proof}
 
 We need to show that $PA(T)_1=TL_1$.
 We proceed very similarly to the above proof.  If we have a diagram with two output strings, we will inductively reduce the number of $T$s inside of it until we have none left.  One to none is easy; since the diagram only has two output strings, some of the strings of $T$ connect back to $T$, so $T$ has an innermost cap somewhere that kills it.  
 
 To go from $k$ to $k-1$, let $\tau$ be a diagram with $k$ $T$s connected up among themselves so that there are only two output strings.  Consider the inner product $\ip{\tau,\tau}=
\begin{tikzpicture}[baseline, scale=.5]
	\draw (0,0) circle (1.1cm);
	\node at (-1,0) [Tbox] {$\tau$};
	\node at (1,0) [Tbox] {$\tau^*$};	
\end{tikzpicture}. $
It's not actually important that this is the inner product diagram.  What matters is that it's a $0$-box, so by Lemma \ref{ngon}, it has either two $T$s connected by three edges, or an $n$-gon with $n-1$ thick sides.  Since there are only two edges between $\tau$ and $\tau^*$, all the vertices involved in the three-connected pair or the $n$-gon lie solely in $\tau$ or $\tau^*$; but $\tau^*$ is a reflection of $\tau$, so $\tau$ has some piece which can be reduced by (R5), or repeated applications of (R6), followed by (R5).  Thus $\tau$ is equal to a sum of diagrams, each of which has at most $(k-1)$ $T$s.  So by induction, $\tau$ is in Temperley-Lieb.
\end{proof}

\begin{Theorem}  $PA(T)$ has principal graph $H$.
\end{Theorem}

\begin{proof}  $PA(T)$ is an irreducible subfactor planar algebra with parameter $\delta=\sqrt{\frac{5+\sqrt{13}}{2}}$.   It follows that the associated subfactor is also irreducible, with index   $\frac{5+\sqrt{13}}{2}$. 
Since $PA(T)$ has a new element at level $4$, the principal graph has a branch and cannot be $A_{\infty}$.  By Haagerup's analysis in \cite{MR1317352}, therefore, $PA(T)$ must have one of the Haagerup graphs as its principal graph.  We know that the principal graph is $H$ and not $H'$ by \cite[equation 4.3.5]{quadratic}.
\end{proof}

\begin{Theorem}  The Haagerup subfactor exists
\end{Theorem}
\begin{proof}
By any of \cite{MR1334479, 0712.2904, 0807.4146}, the existence of a planar algebra with principal graph $H$ implies the existence a subfactor with principal graph $H$.
\end{proof}

%% file: text/OtherGraphs.tex
Having constructed the Haagerup planar algebra in the graph planar algebra of one graph, one naturally wants to know if the Haagerup planar algebra might appear in the graph planar algebras of other bipartite graphs.  In this section we show that many graph planar algebras do not contain a Haagerup planar algebra, and (perhaps more surprisingly) one additional graph planar algebra which does contain the Haagerup planar algebra.  

We start by proving some general theorems about graphs whose graph planar algebras may contain $n$-excess 1 planar algebras (defined below).  We then apply these to the question of finding the Haagerup planar algebra inside the graph planar algebra of particular graphs.  We go on to construct the Haagerup planar algebra inside the graph planar algebra of \begin{tikzpicture}[baseline=-1mm, scale=.7]
	\filldraw (0,0) circle (1mm);
	\filldraw (1,0) circle (1mm);
	\filldraw (2,0) circle (1mm);
	\filldraw (2.7,.7) circle (1mm);
	\filldraw (2.7,-.7) circle (1mm);
	\filldraw (3,0) circle (1mm);
	
	\draw (0,0) -- (3,0);
	\draw (2,0)--(2.7,.7);
	\draw (2,0)--(2.7,-.7);
\end{tikzpicture} .
Many of the results in this section were first seen experimentaly, using a {\tt Mathematica} package called FusionAtlas \cite{fusionatlas}.

In the rest of this section, $n$ will always be an integer greater than $0$.  Recall some notation from \cite{quadratic}:

\begin{Defn} If the codimension of $TL_n$ in the planar algebra $\PP_n$ is $k$, we say that  $\PP$ has $n$-excess $k$.
\end{Defn}

We are mostly interested in the case where a subfactor planar algebra has $n$-excess 1.  This means that its principal graph up to depth $n$ looks like
$$
\begin{tikzpicture}[scale=1.2]
	\filldraw[black] (0,0) circle (1mm);
	\node at (0,-1.2) {$0$};
	\filldraw[black] (1,0) circle (1mm);
	\node at (1,-1.2) {$1$};
	\node at (2,0) [above] {$\ldots$};
	\node at (2,-1.2) {\ldots};	
	\filldraw[black] (3,0) circle (1mm);
	\node at (3,-1.2) {$n-2$};
	\filldraw[black] (4,0) circle (1mm);
	\node at (4,-1.2)  {$n-1$};

	\filldraw[black] (5,.7) circle (1mm);
	\filldraw[black] (5,-.7) circle (1mm);
	\node at (5,-1.2) {$n$};

	\draw (0,0) -- (4,0) -- (5,.7);
	\draw (4,0) -- (5,-.7);

	\draw[dashed] (5,-1) -- (5,1);
\end{tikzpicture}
$$

\begin{Example*}
The Haagerup planar algebra has $4$-excess $1$.
\end{Example*}

\begin{Defn}
Let $\PP$ be an $n$-excess 1 subfactor planar algebra and name the two depth-$n$ minimal idempotents  $e$ and $f$ in such a way that $r=\frac{\tr{f}}{\tr{e}} \geq 1$.  We  define
$$T=-r e + f,
$$
which generates $TL_n^\perp \subset \PP_n$ and has $\tr{T^2}=r[n+1]$ and $\tr{\rho^{1/2}(T)^2}=-r[n+1]$ (see \cite{quadratic}).  
\end{Defn}

\begin{Remark}
If $T$ is orthogonal to Temperley-Lieb, so is $\rho(T)$.  So in an $n$-excess 1 subfactor,  $T$ must be an eigenvalue of the rotation operator.  Thus $T$ generates an irreducible annular Temperley-Lieb module.
\end{Remark}

\begin{Defn}
Let $\PP$ be an $n$-excess 1 subfactor planar algebra, and let $\omega$ be the root of unity such that $\rho(T)=\omega T$.  We call $\omega$ the {\em chirality} of T.  
\end{Defn}

\begin{Example*}
The Haagerup planar algebra has chirality $-1$.
\end{Example*}

Recall from Section \ref{pabg} our notation for paths and loops in a graph planar algebra:

\begin{Notation}
All graphs in this section are simply laced.  Therefore paths or loops can, and will, be entirely described by the vertices they pass through. 
When we concatenate two paths written this way, we have to drop a vertex:
$$a b \ldots c d \sqcup d e \ldots f g = a b \ldots c d e \ldots  f g.$$
When a loop is described in this notation, we might forget to notate the last vertex.

Paths and loops are written in boldface, and indices indicate the position of a vertex in a path or loop: $\pth{p} = p_1 p_2 \ldots p_k$.
\end{Notation}

\begin{Lemma}\label{linearrelns}
If $S$ is an element of $PABG(G)_{n,\pm}$ such that $\epsilon_i(S)=0$ for all $i$, the coefficients in $S$ of loops which are identical except at one vertex are subject to the linear relation
$$
0  = \sum_{ w \text{ adjacent to } x_{i} } \frac{\sqrt{\lambda(w)}}{\sqrt{\lambda(x_{i})}} \cdot \coeff{S}{\pth x w \pth y}  $$
where $\pth x$ is a length $i$ path and $\pth y$ is a length $2n-2-i$ path such that $\pth y$ starts where $\pth x$ ends and $\pth x$ starts where $\pth y$ ends (i.e., $x_i = y_1$, $y_{2n-2-i}=x_1$).
\end{Lemma}

\begin{proof} In the graph planar algebra, we calculate
$$0  = \coeff{ \epsilon_{i} (S)}{\pth x \sqcup \pth y}   
 = \sum_{ w \text{ adjacent to } x_{i} } \frac{\sqrt{\lambda(w)}}{\sqrt{\lambda(x_{i})}} \cdot \coeff{S}{\pth x w \pth y}  
$$
\end{proof}

\begin{Corollary}\label{degreeone}
Suppose $S$ is an element of $PABG(G)_{n,\pm}$, such that $\epsilon_i(S)=0$.  If $v$ is a degree-one vertex, and a loop $\lp p$ ``stays at'' $v$ for two steps (i.e., $p_i=p_{i+2}=v$), then 
$$\coeff{S}{\lp p}=0.$$
\end{Corollary}

\begin{Example*}[$\lp p$ which ``stays at'' a degree-one vertex for two steps]
$$\begin{tikzpicture}[scale=.7]
	\filldraw (0,.5) circle (1mm);
	\filldraw (1,.5) circle (1mm);
	\filldraw (2,.5) circle (1mm);
	
	\draw (0,.5)--(2.5,.5);
	\draw (2.6,.5)--(2.7,.5);
	\draw (2.8,.5)--(2.9,.5);
	
	\draw[rounded corners = .5mm,->] (2,1) -- (0,1)--(0,1.2)--(1,1.2) -- (1,1.4)--(0,1.4)--(0,1.6)--(2,1.6);
	
	\node at (0,.5) [below] {$v$};
	\node at (1,.5) [below] {$w$};
	\node at (2,.5) [below] {$x$};

	\node at (3,1.3) [right] {$\lp p=x w v w v w x$};
\end{tikzpicture}$$
\end{Example*}

\begin{Corollary}\label{flips}
Suppose $S$ is an element of $PABG(G)_{n,\pm}$, such that $\epsilon_i(S)=0$.  If loops $\lp p$, $\lp q$ are related by a series of flips through degree-two vertices, their coefficients in $S$ are subject to a linear relation:
$$\coeff{S}{\lp p} = k \cdot \coeff{S}{\lp q},$$
where $k$ can be calculated from the sequence of flips taking $\lp p$ to $\lp q$.
\end{Corollary}

\begin{Example*}[$\lp p$ and $\lp q$ related by two flips through $w$, then a flip through $x$.]
$$\begin{tikzpicture}[scale=.6, baseline=8mm]
	\filldraw (0,.5) circle (1mm);
	\filldraw (1,.5) circle (1mm);
	\filldraw (2,.5) circle (1mm);
	\filldraw (3,.5) circle (1mm);
	
	\draw (0,.5)--(3.5,.5);
	\draw (3.6,.5)--(3.7,.5);
	\draw (3.8,.5)--(3.9,.5);
	
	\draw[rounded corners = .5mm,->] (3,1) -- (0,1)--(0,1.2)--(1,1.2) -- (1,1.4)--(0,1.4)--(0,1.6)--(3,1.6);

	\node at (3,1.3) [right] {$=\lp p$};
	
	\node at (1,.5) [below] {$w$};
	\node at (2,.5) [below] {$x$};
\end{tikzpicture}
\;\; \rightsquigarrow \;\;
\begin{tikzpicture}[scale=.6,  baseline=8mm]
	\filldraw (0,.5) circle (1mm);
	\filldraw (1,.5) circle (1mm);
	\filldraw (2,.5) circle (1mm);
	\filldraw (3,.5) circle (1mm);
	
	\draw (0,.5)--(3.5,.5);
	\draw (3.6,.5)--(3.7,.5);
	\draw (3.8,.5)--(3.9,.5);
	
	\draw[rounded corners = .5mm,->] (3,1) -- (1,1)--(1,1.2)--(2,1.2) -- (2,1.4)--(0,1.4)--(0,1.6)--(3,1.6);
	
	\node at (1,.5) [below] {$w$};
	\node at (2,.5) [below] {$x$};
\end{tikzpicture}
\;\; \rightsquigarrow \;\;
\begin{tikzpicture}[scale=.6,  baseline=8mm]
	\filldraw (0,.5) circle (1mm);
	\filldraw (1,.5) circle (1mm);
	\filldraw (2,.5) circle (1mm);
	\filldraw (3,.5) circle (1mm);
	
	\draw (0,.5)--(3.5,.5);
	\draw (3.6,.5)--(3.7,.5);
	\draw (3.8,.5)--(3.9,.5);
	
	\draw[rounded corners = .5mm,->] (3,1) -- (1,1)--(1,1.2)--(2,1.2) -- (2,1.4)--(1,1.4)--(1,1.6)--(2,1.6)--(2,1.8)--(1,1.8)--(1,2)--(3,2);	
	
	\node at (1,.5) [below] {$w$};
	\node at (2,.5) [below] {$x$};
\end{tikzpicture}
\;\; \rightsquigarrow \;\;
\begin{tikzpicture}[scale=.6,  baseline=8mm]
	\filldraw (0,.5) circle (1mm);
	\filldraw (1,.5) circle (1mm);
	\filldraw (2,.5) circle (1mm);
	\filldraw (3,.5) circle (1mm);
	
	\draw (0,.5)--(3.5,.5);
	\draw (3.6,.5)--(3.7,.5);
	\draw (3.8,.5)--(3.9,.5);
	
	\draw[rounded corners = .5mm,->] (3,1) -- (1,1)--(1,1.2)--(3,1.2) -- (3,1.4)--(1,1.4)--(1,1.6)--(3,1.6);	

	\node at (1,.5) [below] {$w$};
	\node at (2,.5) [below] {$x$};
	
	\node at (3,1.3) [right] {$=\lp q $};
\end{tikzpicture}
$$
\end{Example*}

\begin{proof}
It is enough to show this for $\lp p$ and $\lp q$ related by one flip.  In this case, it follows immediately from Lemma \ref{linearrelns}.
\end{proof}

\begin{Lemma}\label{traceSsq}
Let $\PP$ be an $n$-excess $1$, chirality $\omega$, $\delta>2$ subfactor planar algebra.  Suppose that for all $S \in PABG(G)$ such that $\epsilon_i(S)=0$ and $\rho(S)=\omega S$, there is some vertex $v$ such that $\coeff{\tr{S^2}}{v}=0$ (if $v$ is even) or $\coeff{\tr{\rho^{1/2}(S)^2}}{v}=0$ (if $v$ is odd).  Then $\PP \not\subset PABG(G)$.
\end{Lemma}

\begin{proof}
Let $T \in \PP$ be the element described in quadratic tangles (and above) which is in $TL_n^{\perp} \subset \PP_n$ and appropriately normalized so that $\tr{T^2} = r[n+1]$, $\tr{\rho^{1/2}(T)^2}  =-r[n+1]$ (where $r$ is a ratio of traces coming from $\PP$'s principal graph).

Suppose for the sake of contradiction that there is an injective planar algebra homomorphism
$$f:\PP \rightarrow PABG(G)$$
and let $S=f(T)$.  Then relations on $T$ must hold on $S$, and so we know that 
\begin{align*}
\epsilon_i(S)& =0, & \rho(S) & =\omega S, \\
\tr{S^2} &= r[n+1] \cdot  \tikz[baseline=-1mm] \draw[thick] (0,0) circle (2mm);  
& \tr{\rho^{1/2}(S)^2} & = \check{r} [n+1] \cdot  \tikz[baseline=-1mm] \filldraw[thick, shaded] (0,0) circle (2mm);
\end{align*}

Recall that as $PABG(G)_{0, \pm}$ is most likely more-than-one-dimensional, the trace tangle is not necessarily a map to $\mathbb{C}$.  However, the claim above is that the trace tangle applied to $S^2$ or $\rho^{1/2}(S)^2$ produces a particular multiple of the empty diagram, which is the only basis element of $TL_{0,\pm}$.  (In particular, we mean $\tr$ and not $Z$ in the formulas above.)

On the one hand, we have
$$\tr{S^2}= r[n+1]  \cdot \tikz[baseline=-1mm] \draw[thick] (0,0) circle (2mm); =  r[n+1] \cdot \sum_{\text{even vertices } w} w.$$
But on the other hand,  $\coeff{\tr{S^2}}{v}=0$.  As $\delta>2$, $[n+1] \neq 0$ and $r$ is never $0$, so this is a contradiction.  Thus there is no injective homomorphism from $\PP$ into $PABG(G)$.

If $v$ is an odd vertex, the same argument works using $\rho^{1/2}(S)$ instead of $S$.
\end{proof}

\begin{Theorem}\label{TooLong}
Let $\PP$ be an $n$-excess $1$, chirality $\omega \neq 1$, $\delta>2$ subfactor planar algebra.  
  If $G$ is a bipartite graph with a vertex $w$ such that
$$d(v,w) \geq n \qquad \text{ for all } v \text{ such that } \operatorname{degree}(v)\geq 3,$$
then $\PP \not\subset PABG(G)$.
\end{Theorem}

\begin{proof}
Let $S$ be an element of $PABG(G)_{n,\pm}$ such that $\epsilon_i(S)=0$ for all $i$, and $\rho(S) = \omega S$.
Suppose $w$ is an even vertex, and let $\lp{p}$ be any length-$2n$ loop in $G$ based at $w$.  
Pick a vertex adjacent to $w$, call it $x$, and let $\lp q$ be the loop, based at $w$, which goes from $w$ to $x$ and back to $w$, $k$ times in total.  

The essential feature of $\lp q$ is that  $\rho( \lp q)=\lp q$.  Because $\rho(S)=\omega S$ , we have $\omega \coeff{S}{\lp q} = \coeff{\omega S}{\lp q}= \coeff{\rho(S)}{\lp q}=\coeff{S}{\rho^{-1}(\lp q)}=\coeff{S}{\lp q}$.  But we also know $\omega \neq 1$, and so $\coeff{S}{\lp q} = 0$.  Because $\lp p$ only passes through degree-two (or degree-one) vertices, it is related by a series of flips to $\lp q$.   Now Lemma \ref{flips} implies that $\coeff{S}{\lp p }=k \cdot \coeff{S}{\lp q}=0$ for all $\lp p$.

Thus, for any loop $\lp p$ based at $w$,
$$\coeff{S^2}{\lp p} 
= \sum_{
\begin{subarray}{l} 
 \pth{s} \text{ length } k \text{ path } p_1 \rightarrow p_{k+1} \\ 
 \pth{r} \text{ length } k \text{ path } p_{k+1} \rightarrow p_{1} 
\end{subarray}} \coeff{S}{p_1 \ldots p_k p_{k+1} \sqcup \pth{r} } \cdot \coeff{S}{\pth{s} \sqcup p_{k+1} \ldots p_{2k} p_1}
=0$$
and therefore 
$$\coeff{\tr{S^2}}{w}=\sum_{\lp p \in PABG(G)_{n,+}, \text{ based at } w} \coeff{S^2}{\lp p} \cdot \tr{\lp p}=0.$$
Now Lemma \ref{traceSsq} implies that $\PP \not\subset PABG(G)$.

If $w$ is an odd vertex, the same argument works using $\rho^{1/2}(S)$ instead of $S$.
\end{proof}

\begin{Defn}[A {\em pruning} of a graph]
A pruning of a graph is a subgraph which preserves degrees, except at one vertex.  Thus 
$$
\begin{tikzpicture}[scale=.7]
	\filldraw (0,0) circle (1mm);
	\filldraw (1,0) circle (1mm);
	\filldraw (2,0) circle (1mm);
	\filldraw (2,1) circle (1mm);
	\filldraw (3,0) circle (1mm);
	
	\draw (0,0) -- (4,0);
	\draw (2,0)--(2,1);
	
	\filldraw[fill=white] (4,0) circle (1mm);
\end{tikzpicture}
$$
is a pruning of
$$
\begin{tikzpicture}[scale=.7]
	\filldraw (0,0) circle (1mm);
	\filldraw (1,0) circle (1mm);
	\filldraw (2,0) circle (1mm);
	\filldraw (2,1) circle (1mm);
	\filldraw (3,0) circle (1mm);
	\filldraw (4,0) circle (1mm);
	\filldraw (4,1) circle (1mm);
	\filldraw (5,0) circle (1mm);
	\filldraw (6,0) circle (1mm);
	
	\draw (0,0) -- (6,0);
	\draw (2,0)--(2,1);
	\draw (4,0)--(4,1);
\end{tikzpicture}
$$
but not of
$$
\begin{tikzpicture}[scale=.7]
	\filldraw (-1,0) circle (1mm);
	\filldraw (0,0) circle (1mm);
	\filldraw (1,0) circle (1mm);
	\filldraw (2,0) circle (1mm);
	\filldraw (2,1) circle (1mm);
	\filldraw (3,0) circle (1mm);
	\filldraw (4,0) circle (1mm);
	\filldraw (5,0) circle (1mm);

	\draw (-1,0) -- (5,0);
	\draw (2,0)--(2,1);
\end{tikzpicture}
$$
\end{Defn}

\begin{Theorem}\label{ForbiddenPruning}
Let $\PP$ be an $n$-excess $1$, $\delta>2$ subfactor planar algebra.  If $G$ is a bipartite graph which has a pruning
$$
\begin{tikzpicture}[scale=.7]
	\filldraw (-1,0) circle (1mm);
	\node at (0,0) [below] {$\cdots$};
	\filldraw (1,0) circle (1mm);
	\draw[snake=brace] (1,-.5) -- node[below] {$n-2$} (-1,-.5);
	\filldraw (2,0) circle (1mm);
	\filldraw (2,1) circle (1mm);
	\filldraw (3,0) circle (1mm);
	\filldraw[fill=white] (4,0) circle (1mm);
	
	\draw (-1,0) -- (4,0);
	\draw (2,0)--(2,1);
	
	\filldraw[fill=white] (4,0) circle (1mm);
\end{tikzpicture}
$$
then $\PP \not\subset PABG(G)$.
\end{Theorem}

\begin{proof}
We need some notation: name the vertices in the pruning $a_1$, \ldots , $a_{n-2}$, $c$, $d$, $e$ and $f$:
$$
\begin{tikzpicture}[scale=.7]
	\filldraw (-1,0) node[below] {$a_1$} circle (1mm);
	\node at (0,0) [below] {$\cdots$};
	\filldraw (1,0) node[below] {$a_{n-2}$} circle (1mm);
	\filldraw (2,0) node[below] {$c$} circle (1mm);
	\filldraw (2,1) node[right] {$d$} circle (1mm);
	\filldraw (3,0) node[below] {$e$} circle (1mm);
	\filldraw[fill=white] (4,0) node[below] {$f$} circle (1mm);
	
	\draw (-1,0) -- (4,0);
	\draw (2,0)--(2,1);
	
	\filldraw[fill=white] (4,0) circle (1mm);
\end{tikzpicture}
$$
If $a_1$ is an even vertex, we are going to show that $\coeff{\tr{S^2}}{a_1}=0$ by showing (as in Theorem \ref{TooLong}) that for any path $\lp p$ based at $a_1$, $\coeff{S}{\lp p}=0$.

Of the loops of length $2n$ in $G$ based at $a_1$, we can immediately show that all but $4$ have coefficient $0$.  
Loops which stay at $a_1$ or $d$ for two steps automatically have coefficient $0$ in $S$ (see Lemma \ref{degreeone}).
Further, loops which stay at $a_i$ for two steps are related by a series of flips to $\ldots a_1 a_2 a_1 \ldots$, and so by Lemma \ref{flips} also have coefficient $0$ in $S$.  So we only need to show the following:
\begin{align*}
\coeff{S}{a_1 \ldots a_{n-2} cdceca_{n-2} \ldots a_1}&=0\\
\coeff{S}{a_1 \ldots a_{n-2}cefeca_{n-2} \ldots a_1}&=0\\
\coeff{S}{a_1 \ldots a_{n-2}cecdca_{n-2} \ldots a_1}&=0\\
\coeff{S}{a_1 \ldots a_{n-2}cececa_{n-2} \ldots a_1}&=0
\end{align*}

Lemma \ref{linearrelns} tells us that  we have a linear relation among $\pth x a_{n-2} \pth y$, $\pth x d \pth y$ and $\pth x e \pth y$ (for $\pth x$ starting at $a_1$ and ending at $c$, $\pth y$ starting at $c$ and ending at $a_1$).    This means that
\begin{multline*}
\coeff{S}{a_1 \ldots a_{n-2}cdceca_{n-2} \ldots a_1}\\
=k_1 \cdot \coeff{S}{a_1 \ldots a_{n-2}cdcbca_{n-2} \ldots a_1} \\
+ k_2 \cdot \coeff{S}{a_1 \ldots a_{n-2}cdcdca_{n-2} \ldots a_1}  =0
\end{multline*}
\begin{multline*}
\coeff{S}{a_1 \ldots a_{n-2}cecdca_{n-2} \ldots a_1}\\
=k_1 \cdot\coeff{S}{ a_1 \ldots a_{n-2}cbcdca_{n-2} \ldots a_1} \\
+ k_2 \cdot \coeff{S}{a_1 \ldots a_{n-2}cdcdca_{n-2} \ldots a_1} = 0
\end{multline*}
\begin{multline*}
\coeff{S}{a_1 \ldots a_{n-2}cececa_{n-2} \ldots a_1}\\
=k_1 \cdot \coeff{S}{a_1 \ldots a_{n-2}cbceca_{n-2} \ldots a_1}  \\
+ k_2 \cdot \coeff{S}{a_1 \ldots a_{n-2}cdceca_{n-2} \ldots a_1} = 0 
\end{multline*}
Now by Lemma \ref{flips}, there is a linear relation
\begin{multline*}
\coeff{S}{a_1 \ldots a_{n-2}cefeca_{n-2} \ldots a_1} \\\
= k \cdot  \coeff{S}{a_1 \ldots a_{n-2}cececa_{n-2} \ldots a_1}=0
\end{multline*}
and so for any $\lp p$, we have $\coeff{S}{\lp p}=0$.  Therefore 
$$\coeff{\tr{S^2}}{a_1}=\sum_{\lp p \in PABG(G)_{n,+}, \text{ based at } a_1} \coeff{S^2}{\lp p} \cdot \tr{\lp p}=0.$$
Now Lemma \ref{traceSsq} implies that $\PP \not\subset PABG(G)$.

If $a_1$ is an odd vertex, the same argument works using $\rho^{1/2}(S)$ instead of $S$.

\end{proof}

\begin{figure}[ht]
\input{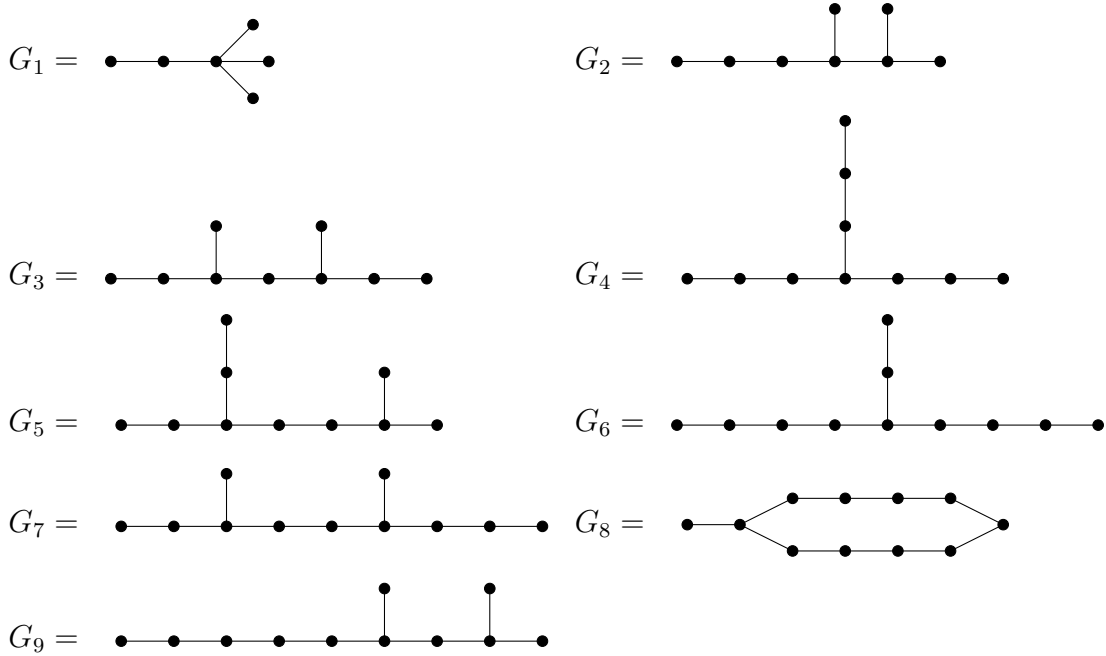}
\caption{Graphs with norm $\sqrt{\frac{5+\sqrt{13}}{2}}$ and fewer than $12$ vertices}\label{upto12}
\end{figure}

\begin{figure}[ht]
\input{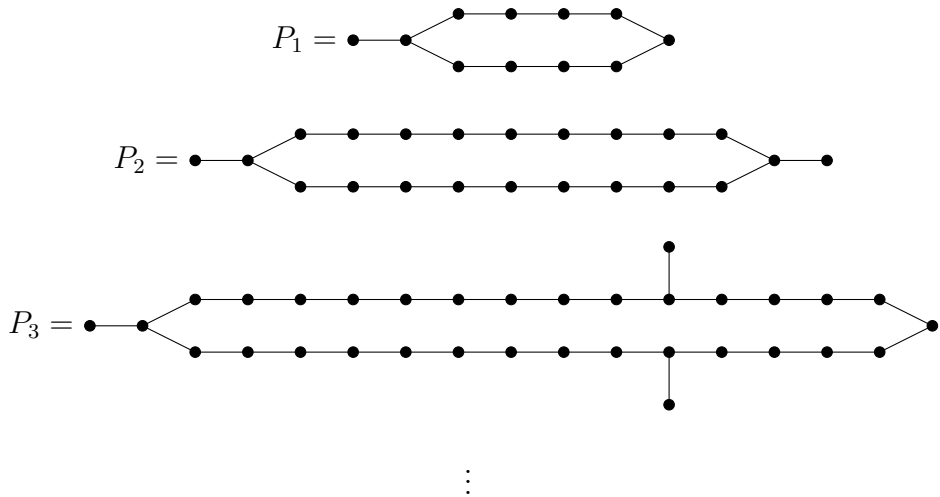}
\caption{An infinite family of graphs with norm $\sqrt{\frac{5+\sqrt{13}}{2}}$.}\label{10gons}
\end{figure}

From now on, $\PP$ is the Haagerup planar algebra.

\begin{Theorem}  The Haagerup planar algebra makes a surprise appearance:
$$\PP \subset PABG( B =
\begin{tikzpicture}[baseline=-1mm, scale=.7]
	\filldraw (0,0) circle (1mm);
	\filldraw (1,0) circle (1mm);
	\filldraw (2,0) circle (1mm);
	\filldraw (2.7,.7) circle (1mm);
	\filldraw (2.7,-.7) circle (1mm);
	\filldraw (3,0) circle (1mm);
	
	\draw (0,0) -- (3,0);
	\draw (2,0)--(2.7,.7);
	\draw (2,0)--(2.7,-.7);
\end{tikzpicture}
)
$$
\end{Theorem}

\begin{proof}
It is sufficient to give an element $T_B \in PABG(B)_{4,+}$ such that 
\begin{align*}
\epsilon_i(T_B)&=0  \text{ for all } i, & 
 Z(T_B^2) & = 3+\sqrt{13},&
Z(\rho^{1/2}(T_B)^2) & = -3 - \sqrt{13},\\
\rho(T_B)& =-T_B, &
Z(T_B^3) & = 0,&
Z(\rho^{1/2}(T_B)^3) & = i\sqrt{\frac{8(3+\sqrt{13})}{3}},\\
 T_B^* & =T_B ,& 
 Z(T_B^4) & = 3+\sqrt{13}, &
 Z(\rho^{1/2}(T_B)^4) & = \frac{17+3\sqrt{13}}{3};
\end{align*}
once we find such a $T_B$, the calculations of Sections \ref{Quadratic} and \ref{PATisHaagerup} all hold.  

Name the vertices of $B$ as
$$
\begin{tikzpicture}[baseline=-1mm, scale=.7]
	\filldraw (0,0) node [above] {$a$} circle (1mm);
	\filldraw (1,0) node [above] {$b$}  circle (1mm);
	\filldraw (2,0) node [above] {$c$}  circle (1mm);
	\filldraw (2.7,.7) node [right] {$d_2$}  circle (1mm);
	\filldraw (2.7,-.7) node [right] {$d_4$}  circle (1mm);
	\filldraw (3,0) node [right] {$d_3$}  circle (1mm);
	
	\draw (0,0) -- (3,0);
	\draw (2,0)--(2.7,.7);
	\draw (2,0)--(2.7,-.7);
\end{tikzpicture};
$$
The  Perron-Frobenius data is $\delta=\sqrt \frac{5+ \sqrt{13}}{2}$ and 
\begin{align*}
\lambda(a) & =\sqrt{\frac{13-3\sqrt{13}}{26}}, &
\lambda(b) & =\sqrt{\frac{13-\sqrt{13}}{26}}, \\
\lambda(c) & =\sqrt{\frac{13+3\sqrt{13}}{26}}, &
\lambda(d_i) & =\sqrt{\frac{13+\sqrt{13}}{78}}
\end{align*}

Now we can define the element $T_B$ which will generate $\PP$.
\begin{align*} 
T_B =  t \cdot \big(1 & -\rho + \rho^2 - \rho^3 \big) \cdot \\
\Big( 
& k \cdot \big( 1 -  (23) -  (24) -  (34) + (234) +  (243) \big)\big(1+f_2 \big) 
\big( cbcbc d_2 c d_3 \big)\\
& + \big( 1 + (234) + (243) \big) \big( cbcd_2 c b c d_3 \big)\\
& -  \big( 1 + (234) + (243) \big) \big( cbcd_2 c b c d_2 \big)\\
& - k \cdot \big( 1 +  (23)  +  (24) + (34) + (234) + (243) \big) \big( cbcd_2 c d_3 c d_4 \big)\\
& + k ( 1+k) \cdot \big( 1 + (234) + (243) \big) \big( cbcd_2 c d_3 c d_2 \big)\\
& + k ( 1-k) \cdot \big( 1 + (234) + (243) \big) \big( cbcd_3 c d_2 c d_3 \big)\\
& + k^2 \cdot  \big( 1 + (234) + (243) \big) \big( c d_2 c d_3 c d_2 c d_4 \big)\\
& + k^3 \cdot  \big( 1  + (234) + (243) \big) \big( c d_2 c d_3 c d_2 c d_3 \big) \Big)
\end{align*}
for $t=\frac{i}{3} \sqrt{-7+2 \sqrt{13}} \approx 0.153153 i$ and $k=\sqrt{\frac{-1+\sqrt{13}}{2}} \approx 1.14139$.

Note that the coefficients in $T_B$ are all purely imaginary.

From this presentation, it is easy to see that $T_B$ is self-adjoint and has rotational eigenvalue $-1$; it is also straightforward to show that is has $\epsilon_i(T_B)=0$ for all $i$.

To verify that the traces of the powers are as desired, we convert $T_B$ into matrix form in Appendix \ref{traces}.  Once we know these traces, all the calculations of Sections \ref{Quadratic} and \ref{PATisHaagerup} carry through and we know that $T_B$ generates the Haagerup planar algebra.
\end{proof}

\begin{Theorem}
$\PP$ is not contained in $PABG(G_i)$ for $G_i$ from Figure \ref{upto12} except for $G_1$, $G_2$ and $G_4$.
$\PP$ is not contained in $PABG(P_i)$ for all $P_i$ from Figure \ref{10gons}. 
\end{Theorem}

\begin{proof}
Graphs $G_6$, $G_8$, $G_9$ and all $P_i$ are ruled out by Theorem \ref{TooLong}.  Graphs $G_3$ and $G_7$ are ruled out by Theorem \ref{ForbiddenPruning}.  We know of no way to rule out $G_5$ except by direct computation (performed in FusionAtlas \cite{fusionatlas}):  One finds all the elements $S \in PABG(G_5)_{4,\pm}$ satisfying $\epsilon_i(S)=0$, $\rho(S)=-S$ and $S^*=S$ and shows that none of them satisfy $\tr{S^2} \in TL_{0,\pm}$.
\end{proof}

%% file: text/Traces.tex
Because $PABG_{k,\pm}$ is a finite-dimensional algebra with a positive definite inner product, it is isomorphic (by the Artin-Wedderburn theorem) to a direct sum of matrix algebras.  The minimal central idempotents are projection onto the space of loops starting at $v$ having midpoint $w$.  The matrix corresponding to loops based at $v$ with midpoint $w$ has rows labelled by length-$k$ paths from $v$ to $w$, and columns labelled by the same paths, traversed backwards.
$$PABG_{k,\pm} \simeq \bigoplus_{( v, w) \in U_\pm \times U_\pm} M_{i_{ v, w} }(\mathbb{C})$$

  Given $m \in PABG_{k,+}$, let $\tilde{m}=\oplus \tilde{m}_{ v, w }$ be its image in $\bigoplus M_{i_{ v, w} }(\mathbb{C})$.  Then 
$$Z(m)=\sum_{( v, w) } \lambda(v) \lambda(w) \tr {\tilde{m}_{v, w }},$$ where $\tr$ is the usual matrix trace.  For this reason, working in $\bigoplus M_{i_{ v, w} }(\mathbb{C})$ is sometimes computationally convenient.

For the purposes of constructing the Haagerup subfactor, we work in $PABG(H)_{4,+}$.  We order  the set of pairs of vertices lexicographically, with the ``alphabet'' ordered $b_0 < z_0 < b_1 < z_1 < b_2 < z_2$.

The rows of $M_{i_{b_0,v}}(\mathbb{C})$ are labelled
\begin{align*}
& \left(
\begin{array}{c}
b_0 a_0 z_0 a_0 b_0 \\
b_0 a_0 b_0 a_0 b_0 \\
b_0 c b_0 a_0 b_0 \\
b_0 a_0 b_0 c b_0 \\
b_0 c b_0 c b_0 \\
b_0 c b_1 c b_0 \\
b_0 c b_2 c b_0 \\
\end{array} \right),
 \left(
\begin{array}{c}
b_0 a_0 z_0 a_0 z_0 \\
b_0 a_0 b_0 a_0 z_0 \\
b_0 c b_0 a_0 z_0 \\
\end{array} \right),
\left(
\begin{array}{c}
b_0 a_0 b_0 c b_1 \\
b_0 c b_0 c b_1 \\
b_0 c b_1 c b_1 \\
b_0 c b_2 c b_1 \\
b_0 c b_1 a_1 b_1 \\
\end{array} \right),
\left(
\begin{array}{c}
b_0 c b_1 a_1 z_1 \\
\end{array} \right),
\\
&
\left(
\begin{array}{c}
b_0 a_0 b_0 c b_2 \\
b_0 c b_0 c b_2 \\
b_0 c b_1 c b_2 \\
b_0 c b_2 c b_2 \\
b_0 c b_1 a_1 b_2 \\
\end{array} \right),
\left(
\begin{array}{c}
b_0 c b_1 a_1 z_2 \\
\end{array} \right)
\end{align*}
and the rows of $M_{i_{z_0,v}}(\mathbb{C})$ are labelled
\begin{align*}
\left(
\begin{array}{c}
z_0 a_0 z_0 a_0 b_0 \\
z_0 a_0 b_0 a_0 b_0 \\
z_0 a_0 b_0 c b_0 \\
\end{array} \right),
\left(
\begin{array}{c}
z_0 a_0 z_0 a_0 z_0 \\
z_0 a_0 b_0 a_0 z_0 \\
\end{array} \right),
\left(
\begin{array}{c}
z_0 a_0 b_0 c b_1 \\
\end{array} \right),
\left(
\begin{array}{c}
z_0 a_0 b_0 c b_2 \\
\end{array} \right)
\end{align*}

The labels of the rows of $M_{i_{z_1,v}}(\mathbb{C})$ and $M_{i_{b_1,v}}(\mathbb{C})$ are as above, but with the permutation $(012)$ applied to the indices of $z$, $a$ and $b$; the labels of the rows of $M_{i_{z_2,v}}(\mathbb{C})$ and $M_{i_{b_2,v}}(\mathbb{C})$ are also as above, but with the permutation $(021)$ applied to the indices.

Then, if we set
$$
\theta = \sqrt{\frac{1+\sqrt{13}}{2}}, \eta =  \sqrt{\frac{3+\sqrt{13}}{2}}, t =\frac{4-\sqrt{13}}{3}
$$ 
we have
\begin{align*}
\tilde{T}= &
\biggl(
\left(
\begin{array}{ccccccc}
 0 & 0 & 0 & 0 & 0 & -t \eta  \theta ^2 & t \eta  \theta ^2 \\
 0 & 0 & 0 & 0 & 0 & t \theta ^2 & -t \theta ^2 \\
 0 & 0 & 0 & 0 & 0 & -t \theta  & t \theta  \\
 0 & 0 & 0 & 0 & 0 & -t \theta  & t \theta  \\
 0 & 0 & 0 & 0 & 0 & t & -t \\
 -t \eta  \theta ^2 & t \theta ^2 & -t \theta  & -t \theta  & t & -t & 0 \\
 t \eta  \theta ^2 & -t \theta ^2 & t \theta  & t \theta  & -t & 0 & t
\end{array}
\right)
\oplus
\left(
\begin{array}{ccc}
 0 & 0 & 0 \\
 0 & 0 & 0 \\
 0 & 0 & 0
\end{array}
\right)
\\
&  \oplus 
\left(
\begin{array}{ccccc}
 -t \theta ^2 & t \theta  & -i t \theta ^2 & t (-1+i \theta ) \theta  & i t \theta ^3 \\
 t \theta  & -t & i t \theta  & t (1-i \theta ) & -i t \theta ^2 \\
 i t \theta ^2 & -i t \theta  & t & t (-1+i \theta ) & -t \theta  \\
 t (-1-i \theta ) \theta  & t (1+i \theta ) & t (-1-i \theta ) & 0 & t (1+i \theta ) \theta  \\
 -i t \theta ^3 & i t \theta ^2 & -t \theta  & t (1-i \theta ) \theta  & t \theta ^2
\end{array}
\right) 
\oplus 
\left(
\begin{array}{c}
 -t \eta ^2 \theta ^2
\end{array}
\right)
\\
& \oplus 
\left(
\begin{array}{ccccc}
 t \theta ^2 & -t \theta  & t (1+i \theta ) \theta  & -i t \theta ^2 & i t \theta ^3 \\
 -t \theta  & t & t (-1-i \theta ) & i t \theta  & -i t \theta ^2 \\
 t (1-i \theta ) \theta  & t (-1+i \theta ) & 0 & t (1-i \theta ) & t (-1+i \theta ) \theta  \\
 i t \theta ^2 & -i t \theta  & t (1+i \theta ) & -t & t \theta  \\
 -i t \theta ^3 & i t \theta ^2 & t (-1-i \theta ) \theta  & t \theta  & -t \theta ^2
\end{array}
\right) 
 \oplus 
 \left(
\begin{array}{c}
 t \eta ^2 \theta ^2
\end{array}
\right)
\\
& \oplus 
\left(
\begin{array}{ccc}
 0 & 0 & 0 \\
 0 & 0 & 0 \\
 0 & 0 & 0
\end{array}
\right)
\oplus
\left(
\begin{array}{cc}
 0 & 0 \\
 0& 0
\end{array}
\right) 
 \oplus 
 \left(
\begin{array}{c}
 t \eta ^2 \theta ^2
\end{array}
\right) 
\oplus 
\left(
\begin{array}{c}
 -t \eta ^2 \theta ^2
\end{array}
\right) 
\biggr)^{\otimes 3}
\end{align*}

From this representation of $T$, it is easy to check that 
\begin{align*}
Z(T^2) & =3+\sqrt{13} \\ 
Z(T^3) & =0 \\
Z(T^4) & =3+\sqrt{13}.
\end{align*}

For $PABG(H)_{4,-}$, we can similarly represent $\rho^{1/2} T$ as a direct sum of matrices; from which it follows that 

\begin{align*}
Z((\rho^{1/2} T)^2) & = -3 - \sqrt{13} \\
Z((\rho^{1/2} T)^3)& =i\sqrt{\frac{8(3+\sqrt{13})}{3}} \\
Z((\rho^{1/2}T)^4) & =\frac{17+3\sqrt{13}}{3} \\ 
\end{align*}

To construct the Haagerup subfactor as a sub-planar algebra of $PABG(B)$, we need to similarly present its generator $T_B$.

The even vertices of $B$ are ordered $a < c$; the rows of $M_{i_{a,v}}(\mathbb{C})$ are labelled 
\begin{align*}
\left(
	\begin{array}{c}
	ababa \\
	abcba \\
	\end{array} \right),
\left(
	\begin{array}{c}
	ababc \\
	abcbc \\
	abcd_2c \\
	abcd_3c \\
	abcd_3c \\
	\end{array} \right),
\end{align*}
and the rows of $M_{i_{c,v}}(\mathbb{C})$ are ordered
\begin{align*}
\left(
	\begin{array}{c}
	cbaba \\
	cbcba \\
	cd_2cba \\
	cd_3cba \\
	cd_4cba
	\end{array} 
\right),
\left(
	\begin{array}{c}
	cbabc \\
	cbcbc \\
	cbcd_2c \\
	cbcd_3c \\
	cbcd_4c \\
	cd_2cbc \\
	cd_2cd_2c \\
	cd_2cd_3c \\
	cd_2cd_4c \\
	cd_3cbc \\
	cd_3cd_2c \\
	cd_3cd_3c \\
	cd_3cd_4c \\
	cd_4cbc \\
	cd_4cd_2c \\
	cd_4cd_3c \\
	cd_4cd_4c 
	\end{array} 
\right).
\end{align*}

Then, if we set 
\begin{align*}
t & =\frac{i}{3} \sqrt{-7+2 \sqrt{13}} \\
k & =\sqrt{\frac{\lambda(b)}{\lambda(d)}}=\sqrt{\frac{-1+\sqrt{13}}{2}} \\
m &=\sqrt{\frac{\lambda(c)}{\lambda(a)}}=\sqrt{\frac{3+\sqrt{13}}{2}} \\
n & =k(1+k) \\
n' & =k(1-k)
\end{align*}
we have
{\small \begin{multline*}
\tilde{T}_B= t \Biggl(
\left(
	\begin{array}{cc}
	 0 & 0 \\
	 0 & 0 
	 \end{array}
\right) 
\oplus
\left(
	\begin{array}{ccccc}
	 0 & 0 & 0 & 0 & 0 \\
	 0 & 0 & 0 & 0 & 0 \\
	 0 & 0 & 0 & k  m ^2 & -k  m ^2 \\
	 0 & 0 & -k  m ^2 & 0 & k  m ^2 \\
	 0 & 0 & k  m ^2 & -k  m ^2 & 0
	\end{array}
\right)
\oplus
\left(
	\begin{array}{ccccc}
	 0 & 0 & 0 & 0 & 0 \\
	 0 & 0 & 0 & 0 & 0 \\
	 0 & 0 & 0 & -k  m^2  & k  m^2  \\
	 0 & 0 & k  m^2  & 0 & -k  m^2  \\
	 0 & 0 & -k  m^2  & k  m^2  & 0
\end{array}
\right) \oplus \\
\left(
	\begin{array}{ccccccccccccccccc}
	 0 & 0 & 0 & 0 & 0 & 0 & 0 & k  m  & -k  m  & 0 & -k  m  & 0 & k  m  & 0 & k  m  & -k  m  & 0 \\
	 0 & 0 & 0 & 0 & 0 & 0 & 0 & -k  & k  & 0 & k  & 0 & -k  & 0 & -k  & k  & 0 \\
	 0 & 0 & 0 & -k  & k  & -2 & 0 & n  & n ' & 1 & 0 & 0 & -k  & 1 & 0 & -k  & 0 \\
	 0 & 0 & k  & 0 & -k  & 1 & 0 & 0 & -k  & -2 & n ' & 0 & n  & 1 & -k  & 0 & 0 \\
	 0 & 0 & -k  & k  & 0 & 1 & 0 & -k  & 0 & 1 & -k  & 0 & 0 & -2 & n  & n ' & 0 \\
	 0 & 0 & 2 & -1 & -1 & 0 & 0 & 0 & 0 & k  & -n  & 0 & k  & -k  & -n ' & k  & 0 \\
	 0 & 0 & 0 & 0 & 0 & 0 & 0 & 0 & 0 & 0 & 0 & 0 & 0 & 0 & 0 & 0 & 0 \\
	 -k  m  & k  & -n  & 0 & k  & 0 & 0 & 0 & 0 & n ' & k ^3 & 0 & -k ^2 & -k  & k ^2 & 0 & 0 \\
	 k  m  & -k  & -n ' & k  & 0 & 0 & 0 & 0 & 0 & -k  & k ^2 & 0 & 0 & n  & -k ^3 & -k ^2 & 0 \\
	 0 & 0 & -1 & 2 & -1 & -k  & 0 & -n ' & k  & 0 & 0 & 0 & 0 & k  & k  & -n  & 0 \\
	 k  m  & -k  & 0 & -n ' & k  & n  & 0 & -k ^3 & -k ^2 & 0 & 0 & 0 & 0 & -k  & 0 & k ^2 & 0 \\
	 0 & 0 & 0 & 0 & 0 & 0 & 0 & 0 & 0 & 0 & 0 & 0 & 0 & 0 & 0 & 0 & 0 \\
	 -k  m  & k  & k  & -n  & 0 & -k  & 0 & k ^2 & 0 & 0 & 0 & 0 & 0 & n ' & -k ^2 & k ^3 & 0 \\
	 0 & 0 & -1 & -1 & 2 & k  & 0 & k  & -n  & -k  & k  & 0 & -n ' & 0 & 0 & 0 & 0 \\
	 -k  m  & k  & 0 & k  & -n  & n ' & 0 & -k ^2 & k ^3 & -k  & 0 & 0 & k ^2 & 0 & 0 & 0 & 0 \\
	 k  m  & -k  & k  & 0 & -n ' & -k  & 0 & 0 & k ^2 & n  & -k ^2 & 0 & -k ^3 & 0 & 0 & 0 & 0 \\
	 0 & 0 & 0 & 0 & 0 & 0 & 0 & 0 & 0 & 0 & 0 & 0 & 0 & 0 & 0 & 0 & 0
	\end{array}
\right) \Biggr)
\end{multline*} }

Again, from this representation of $T_B$ it is easy to check that 
\begin{align*}
Z(T^2) & =3+\sqrt{13} \\ 
Z(T^3) & =0 \\
Z(T^4) & =3+\sqrt{13}.
\end{align*}
and for $PABG(B)_{4,-}$, we can similarly represent $\rho^{1/2} T$ as a direct sum of matrices; from which it follows that 
\begin{align*}
Z((\rho^{1/2} T)^2) & = -3 - \sqrt{13} \\
Z((\rho^{1/2} T)^3)& =i\sqrt{\frac{8(3+\sqrt{13})}{3}} \\
Z((\rho^{1/2}T)^4) & =\frac{17+3\sqrt{13}}{3}. \\ 
\end{align*}